\documentclass[11pt,oneside]{article}

\usepackage{amsthm}
\usepackage{amsmath,amssymb,amsfonts}

\usepackage{graphicx}

\usepackage{color}
\usepackage[dvipsnames]{xcolor}

\def\mym{ {m} }

\def\me{\mathrm e}

\def\dif{\mathrm d}

\def\T{ {\mathrm{\scriptscriptstyle T}} }

\def\argmin{\mathrm{argmin}}
\def\TV{\mbox{TV}}

\def\R{{\mathbb{R}}}

\newenvironment{prf}
{\noindent \textbf{Proof.}}{\hfill $\Box$ \vspace{.1in}}

\newtheorem{thm}{Theorem}
\newtheorem{lem}{Lemma}
\newtheorem{pro}{Proposition}
\newtheorem{cor}{Corollary}
\newtheorem{ass}{Assumption}

\theoremstyle{definition}
\newtheorem{eg}{Example}

\theoremstyle{definition}
\newtheorem{rem}{Remark}

\begin{document}

\begin{titlepage}

\begin{center}
%{\bf On penalized estimation \chz{in} functional additive regression with high-dimensional data}
{\Large\sc Penalized Estimation in Additive Regression\\ with High-Dimensional Data}

\vspace{.1in} Zhiqiang Tan\footnotemark[1] \& Cun-Hui Zhang \footnotemark[1]

\vspace{.1in}
\today
\end{center}

\footnotetext[1]{Department of Statistics \& Biostatistics, Rutgers University. Address: 110 Frelinghuysen Road,
Piscataway, NJ 08854. E-mail: ztan@stat.rutgers.edu, czhang@stat.rutgers.edu. The research of Z.~Tan was supported in part by PCORI grant ME-1511-32740.
The research of C.-H. Zhang was supported in part by NSF grants DMS-1513378, IIS-1250985, and IIS-1407939.}

\paragraph{Abstract.} Additive regression provides an extension of linear regression by modeling
the signal of a response as a sum of functions of covariates of relatively low complexity.
We study penalized estimation in high-dimensional nonparametric
additive regression where functional semi-norms are
used to induce smoothness of component functions and the empirical $L_2$ norm is used to induce sparsity.
The functional semi-norms can be of Sobolev or bounded variation types and are allowed to be different
amongst individual component functions.
We establish new oracle inequalities for the predictive performance of such methods
under three simple technical conditions: a sub-gaussian condition on the noise,
a compatibility condition on the design and the functional classes under consideration,
and an entropy condition on the functional classes.
For random designs, the sample compatibility condition can be replaced by its population version
under an additional condition to ensure suitable convergence of empirical norms.
In homogeneous settings where the complexities of the component functions are of the same order,
our results provide a spectrum of explicit convergence rates, from the so-called slow rate without
requiring the compatibility condition to the fast rate under the hard sparsity or certain $L_q$ sparsity
to allow many small components in the true regression function.
These results significantly broadens and sharpens existing ones in the literature.

%\vspace{-.2in}

\end{titlepage}

\section{Introduction} \label{sect:intro}

Additive regression is an extension of linear regression where the signal of a response can be
written as a sum of functions of covariates of relatively low complexity.
Let $(Y_i,X_i)$, $i=1,\ldots, n$, be a set of $n$ independent (possibly non-identically distributed) observations,
where $Y_i\in \R$ is a response variable and $X_i\in \R^d$ is a covariate (or design) vector.
Consider an additive regression model, $Y_i = g^*(X_i) + \varepsilon_i$ with
\begin{align}
g^*(x) = \hbox{$\sum_{j=1}^p$ } g^*_j( x^{(j)}),\label{additive-reg}
\end{align}
where $\varepsilon_i$ is a noise with mean 0 given $X_i$,
$x^{(j)}$ is a vector composed of a small subset of the components of $x\in \R^d$,
and $g^*_j$ belongs to a certain functional class $\mathcal G_j$.
That is, $g^*(x)$ lies in the space of additive functions
$\mathcal G = \{\sum_{j=1}^p g_j (x^{(j)}): g_j \in \mathcal G_j, j=1,\ldots,p\}$.
A function $g \in \mathcal G$ may admit the decomposition
$g(x) = \sum_{j=1}^p g_j( x^{(j)})$ for multiple choices of $(g_1,\ldots,g_p)$.
In what follows, such choices are considered equivalent but a favorite decomposition
can be used to evaluate properties of the components of $g \in \mathcal G$.

In a classical setting (e.g., Stone 1985), each $g_j^*$ is a univariate function and $x^{(j)}$ is the $j$th
component of $x\in [0,1]^d$, so that $p=d$.
We take a broad view of additive regression and our analysis will accommodate the general setting
where $g^*_j$ can be multivariate with $X_i^{(j)}$ being a block of covariates,
possibly overlapping across different $j$ as in functional ANOVA (e.g., Gu 2002).
However, most concrete examples will be given in the classical setting.

Additive modeling has been well studied in the setting where the number of components $p$ is fixed.
%the sample size $n$ is much larger than the number of functions $p$.
See Hastie \& Tibshirani (1990) and references therein.
Recently, building upon related works in penalized linear regression,
there have been considerable progresses in the development of theory and methods for sparse additive regression
in high-dimensional settings where $p$ can be of greater order than the sample size $n$
but the number of significant components is still smaller than $n$.
See, for example, Lin \& Zhang (2006), Meier et~al. (2009), Ravikumar et~al.~(2009), Huang et~al.~(2010), Koltchinskii \& Yuan (2010), Raskutti et~al.~(2012),
Suzuki \& Sugiyama (2013),
Petersen et~al.~(2016), and Yuan \& Zhou (2016).

In this article, we study a penalized estimator $\hat g$ with a specific associated decomposition
$\hat g = \sum_{j=1}^p \hat g_j$ defined as a minimizer of
a penalized loss
\begin{align*}
\| Y - g \|_n^2/2 +  %A_0
\hbox{$\sum_{j=1}^p$}\big( \rho_{nj} \|g_j\|_{F,j} + \lambda_{nj} \|g_j\|_n \big)
\end{align*}
over $g\in \mathcal G$ and decompositions $g=\sum_{j=1}^p g_j$,
where %$A_0>1$ is a technical constant (see Section~\ref{sect:main-results}),
$(\lambda_{nj}, \rho_{nj})$ are tuning parameters,
$\|\cdot\|_n$ is the empirical $L_2$ norm based on the data points,
e.g. $\|Y-g\|_n^2 = n^{-1} \sum_{i=1}^n \{Y_i-g(X_i)\}^2$,
%and $\| g_j \|_n^2 = n^{-1} \sum_{i=1}^n g_j^2(X_i^{(j)})$,
and $\|g_j\|_{F,j}$ is a semi-norm describing the complexity of $g_j\in \mathcal G_j$.
For simplicity, the association of $\|g_j\|_n$ and $\|g_j\|_{F,j}$ with $X_i^{(j)}$ is typically suppressed.
%Our analysis holds for the decomposition $\hat g = \sum_{j=1}^p \hat g_j\in \mathcal G$
%derived from the minimization of $K_n(g)$.

In the above penalty function, the primary role of the empirical norm $\|\cdot\|_n$ is to
induce sparsity, whereas the primary role of the functional semi-norm $\|\cdot\|_{F,j}$ is to induce
smoothness of the estimated regression function.
%Previously, similar penalties were used for reproducing kernel Hilbert spaces, but with the semi-norm $\|\cdot\|_{F,j}$ replaced by the Hilbert norm (Koltchinskii \& Yuan 2010).
For example, $\|g_j\|_{F,j} = \{\int_0^1 (g_j^{(\mym)})^2 \,\dif z \}^{1/2}$ when
$\mathcal G_j$ is the $L_2$-Sobolev space $\mathcal W^m_2$ on $[0,1]$, where $g_j^{(\mym)}$ denotes
the $\mym$th derivative of $g_j$.

%The penalty above combines the semi-norm $\|\cdot\|_{F,j}$ (inducing smoothness) and
%the empirical norm $\|\cdot\|_n$ (inducing sparsity).
%Previously, similar penalties were used for reproducing kernel Hilbert spaces, but with the semi-norm $\|\cdot\|_{F,j}$ replaced by the Hilbert norm (Koltchinskii \& Yuan 2010).
%For $\mathcal G_j$ as the Sobolev Hilbert space $\mathcal W^\mym_2$ on $[0,1]$, the semi-norm is
%$\|g_j\|_{F,j} = \{\int_0^1 (g_j^{(\mym)})^2 \,\dif z \}^{1/2}$, where $g_j^{(\mym)}$ denotes $\mym$th derivative of $g_j$.
%In general, we allow each class $\mathcal G_j$ to be a Sobolev space $\mathcal W_r^\mym$ or a bounded variation space $\mathcal V^\mym$.
%See Section~\ref{sect:function-class} for a review of univariate functional spaces.

We consider both fixed and random designs and establish oracle inequalities
for the predictive performance of $\hat g$ under three simple technical conditions:
a sub-gaussian condition on noises, a compatibility condition on the design and
the functional classes $\mathcal G_j$, and an entropy condition on $\mathcal G_j$.
The compatibility condition is similar to the restricted eigenvalue condition used in analysis of Lasso,
and for random designs, the empirical compatibility condition can be replaced by its population version
under an additional condition to ensue suitable convergence of empirical norms.
For the Sobolev and bounded variation classes, the entropy condition on $\mathcal G_j$
follows from standard results in the literature (e.g., Lorentz et~al.~1996).

%We consider both fixed and random designs and establish new oracle inequalities for the predictive performance of
%$\hat g$ under three simple technical conditions: sub-gaussian noises,
%empirical or theoretical compatibility condition (similar to the restricted eigenvalue condition used in analysis of Lasso),
%and, for random-design results, a rate condition related to the convergence of empirical norms,
%in addition to an entropy condition that can be satisfied using standard estimates of entropies for Sobolev and bounded variation spaces.
%See Theorems~\ref{main-thm} and~\ref{main2-thm}.

The implications of our oracle inequalities can be highlighted in the classical homogeneous setting
where $X_i^{(j)}$ is the $j$th component of $X_i$ and $\mathcal G_j=\mathcal G_0$ for all $j$,
where $\mathcal G_0$ is either an $L_r$-Sobolev space $\mathcal W_r^\mym$
or a bounded variation space $\mathcal V^\mym$ of univariate functions on $[0,1]$,
where $r\ge 1$ and $\mym \ge 1$ are shape and smoothness indices of the space,
and $r=1$ for $\mathcal V^\mym$.
%$\mathcal G_1,\ldots,\mathcal G_p$ are identical classes of univariate functions on $[0,1]$
%(denoted by $\mathcal G_0$)
In this setting, it is natural to set $(\lambda_{nj},\rho_{nj}) = (\lambda_n,\rho_n)$ for all $j$. %$j=1,\ldots,p$.
Consider random designs, and suppose that for some choice of $(g^*_1,\ldots,g^*_p)$ satisfying (\ref{additive-reg}),
\begin{align}
\sum_{j=1}^p \|g^*_j\|_{F,j} \le C_1 M_F, \quad \sum_{j=1}^p \| g^*_j\|_Q^q  \le C_1^q M_q, \label{Lq-ball}
\end{align}
where $\|f\|_Q^2 =n^{-1} \sum_{i=1}^n E\{f^2(X_i)\}$
for a function $f(x)$,
$C_1>0$ is a scaling constant depending only on the moments of $\varepsilon_i$, and
$0 \le q \le 1$, $M_q>0$ and $M_F>0$ are allowed to depend on $(n,p)$. In the case of hard sparsity, $q=0$,
$M_0 = \#\{j: g^*_j\neq 0\}$.
As a summary, the following result can be easily deduced from Proposition~\ref{cor-slow}, \ref{cor-fast0}, and \ref{cor-medium0}.

Let $\beta_0=1/m$ and define
\begin{align*}
& w^*_n(q) = \max \left\{ n^{\frac{-1}{2+\beta_0(1-q)}}, \left(\log(p)/n\right)^{\frac{1-q}{2}} \right\} ,\\
& \gamma^*_n(q) =\min \left\{ n^{\frac{-1}{2+\beta_0(1-q)}}, n^{-1/2} \left(\log(p)/n\right)^{\frac{-(1-q)\beta_0}{4}} \right\}.
\end{align*}
For simplicity, we restrict to the case where $1 \le r \le 2$.
For $rm>1$, we assume that the average marginal density of
$(X^{(j)}_1,\ldots,X^{(j)}_n)$ are uniformly bounded away from 0 and, if $q\not=1$, also uniformly bounded from above for all $j=1,\ldots,p$.
The assumption of marginal densities bounded from above, as well as the restriction $1 \le r \le 2$, can be relaxed
under slightly different technical conditions (see Propositions~\ref{cor-slow}, \ref{cor-fast}, and \ref{cor-medium}).
For $r=m=1$, neither the lower bound nor the upper bound of marginal densities need to be assumed.

\begin{pro} \label{pro-summary}
Let $\mathcal G_0$ be a Sobolev space $\mathcal W_r^\mym$ with $1 \le r \le 2$ and $\mym \ge 1$
or a bounded variation space $\mathcal V^\mym$ with $r=1$ and $\mym \ge 1$.
Suppose that
the noises are sub-gaussian and $\log(p)=o(n)$.
Let $\tau_0 = 1/(2\mym+1-2/r)$, $\Gamma_n=1$ for $r \mym>1$ and $\Gamma_n =\sqrt{\log n}$ for $r=\mym=1$.

(i) Let $q=1$ and $\lambda_n = \rho_n = A_0\{\log(p)/n\}^{1/2}$ for a sufficiently large constant $A_0$. If $p\to\infty$, then
\begin{align}
\| \hat g - g^*\|_Q^2 = O_p(1) C_1^2(M_F^2 + M_1^2) \left\{n^{-1/2} \Gamma_n+ \sqrt{\log(p)/n} \right\}.  \label{intro-result1}
\end{align}

(ii) Let $q=0$, $\lambda_n =A_0[\gamma^*_n(0) + \{\log(p)/n\}^{1/2}]$
and $\rho_n = \lambda_n w^*_n(0)$.
Suppose that
\begin{align}
\left\{ {w^*_n(0)}^{-\tau_0} \sqrt{\log(np)/n}\right\} (1 + M_F + M_0) = o(1) \label{intro-condition2}
\end{align}
and a population compatibility condition (Assumption~\ref{compat-condition-Q}) holds. Then,
\begin{align}
\| \hat g - g^*\|_Q^2 = O_p(1)C_1^2 ( M_F + M_0) \left\{ n^{\frac{-1}{2+\beta_0} }  + \sqrt{\log(p)/n}\right\}^2.  \label{intro-result2}
\end{align}

(iii) Let $0<q<1$, $\lambda_n = A_0[\gamma^*_n(q) + \{\log(p)/n\}^{1/2}]$
and $\rho_n = \lambda_n w^*_n(q)$.
Suppose that
\begin{align*}
\left\{{w^*_n(q)}^{-\tau_0} \left(\log(np)/n\right)^{\frac{1-q}{2}} \right\} (1 + M_F + M_q) = O(1)
\end{align*}
and a population compatibility condition (Assumption~\ref{compat-mono}) holds. Then,
\begin{align}
\| \hat g - g^*\|_Q^2 = O_p(1)C_1^2 ( M_F + M_q)  \left\{ n^{\frac{-1}{2+\beta_0(1-q)} }  + \sqrt{\log(p)/n}\right\}^{2-q}.  \label{intro-result3}
\end{align}
\end{pro}

There are several important features achieved by the foregoing result, distinct from existing results.
First, our results are established for additive regression with
Sobolev spaces of general shape  and bounded variation spaces.
An important innovation in our proofs involves a delicate application of
maximal inequalities based on the metric entropy of a particular choice of bounded subsets of $\mathcal G_0$ (see Lemma~\ref{max-ineq2}).
All previous results seem to be limited to the $L_2$-Sobolev spaces or similar reproducing kernel Hilbert spaces,
except for Petersen et~al.~(2016), who studied additive regression with
the bounded variation space $\mathcal V^1$ and obtained the rate $\{\log(np)/n\}^{1/2}$ for in-sample prediction under assumption (\ref{Lq-ball}) with $q=1$.
In contrast, our analysis in the case of $q=1$ yields the sharper, yet standard, rate $\{\log(p)/n\}^{1/2}$ for in-sample prediction (see Proposition~\ref{cor-slow}),
whereas $\{\log(np)/n\}^{1/2}$ for out-of-sample prediction by (\ref{intro-result1}).

Second, the restricted parameter set (\ref{Lq-ball}) represents
an $L_1$ ball in $\|\cdot\|_F$ semi-norm (inducing smoothness) but an $L_q$ ball in $\|\cdot\|_Q$ norm (inducing sparsity)
for the component functions $(g^*_1, \ldots, g^*_p)$.
That is, the parameter set (\ref{Lq-ball}) decouples conditions for sparsity and smoothness
in additive regression: it can encourage sparsity at different levels $0\le q\le 1$
while enforcing smoothness only to a limited extent.
Accordingly,
our result leads to a spectrum of convergence rates, which are easily seen to slow down as $q$ increases from 0 to 1,
corresponding to weaker sparsity assumptions. While most of previous results are obtained under exact sparsity ($q=0$), Yuan \& Zhou (2016) studied
additive regression with reproducing kernel Hilbert spaces under an $L_q$ ball in the Hilbert norm $\|\cdot\|_H$: $\sum_{j=1}^p \|g^*_j\|_H^q \le M_q$.
This parameter set induces smoothness and sparsity simultaneously and is in general more restrictive than (\ref{Lq-ball}).
As a result, the minimax rate of estimation obtained by Yuan \& Zhou (2016), based on constrained least squares with {\it known} $M_q$ instead of penalized estimation,
is faster than (\ref{intro-result3}), in the form
$n^{-2/(2+\beta_0)}  +  \{\log(p)/n\}^{(2-q)/2}$, unless $q=0$ or 1.

Third, in the case of $q=1$, our result (\ref{intro-result1}) shows that the rate $\{\log(p)/n\}^{1/2}$, with an additional
$\{\log (n)/n\}^{1/2}$ term for the bounded variation space $\mathcal V^1$,
can be achieved via penalized estimation without requiring a compatibility condition.
This generalizes a slow-rate result for constrained least-squares (instead of penalization) with {\it known} $(M_1,M_F)$ in additive regression with the Sobolev Hilbert space in Ravikumar et~al.~(2009).
Both are related to earlier results for linear regression (Greenhstein \& Ritov 2004; Bunea et~al.~2007).

Finally, compared with previous results giving the same rate of convergence (\ref{intro-result2}) under exact sparsity ($q=0$) for Hilbert spaces, our results
are stronger in requiring much weaker technical conditions.
The penalized estimation procedures in Koltchinskii \& Yuan (2010) and Raskutti et~al.~(2012), while
minimizing a similar criterion as $K_n(g)$, involve additional constraints:
Koltchinskii \& Yuan (2010) assumed that the sup-norm of possible $g^*$ is bounded by a known constant,
where as Raskutti et~al.~(2012) assumed $\max_j \|g_j\|_H$ is bounded by a known constant.
Moreover, Raskutti et~al.~(2012) assumed that the covariates $(X_i^{(1)},\ldots,X_i^{(p)})$ are independent
of each other.
These restrictions were relaxed in Suzuki \& Sugiyama (2013), but only explicitly under the assumption that the noises $\varepsilon_i$ are uniformly bounded by a constant.
Moreover, our rate condition (\ref{intro-condition2}) about the sizes of $(M_0,M_F)$ is much weaker
than in Suzuki \& Sugiyama (2013), due to
improved analysis of convergence of empirical norms and the more careful choices $(\lambda_n,\rho_n)$.
For example, if $(M_0,M_F)$ are bounded, then condition (\ref{intro-condition2}) holds whenever $\log(p)/n=o(1)$ for Sobolev Hilbert spaces, but
the condition previously required amounts to $\log(p) n^{-1/2}=o(1)$.
Finally, the seemingly faster rate in Suzuki \& Sugiyama (2013) can be deduced from our results
when $(\lambda_n,\rho_n)$ is allowed to depend on $(M_0,M_F)$. See Remarks~\ref{rem-main2-cor5}
and \ref{rem:cor-fast1}--\ref{rem:cor-fast3} for relevant discussion.

The rest of the article is organized as follows. Section~\ref{sect:function-class} gives a review of univariate functional classes and entropies.
Section~\ref{sect:main-results} presents general results for fixed designs (Section~\ref{sect:fixed-design}) and random designs (Section~\ref{sect:random-design}),
and then provides specific results with Sobolev and bounded variation spaces (Section~\ref{sect:specific-result}) after a study of convergence of empirical norms (Section~\ref{sect:conv-norm}).
Section~\ref{sect:conclusion} concludes with a discussion. For space limitation, all proofs are collected in Section~\ref{sect:proofs} and technical tools are stated in Section~\ref{sect:tech-tools} of
the Supplementary Material.

\section{Functional classes and entropies} \label{sect:function-class}

As a building block of additive regression, we discuss two broad choices for the function space
$\mathcal G_j$ and the associated semi-norm $\| g_j \|_{F,j}$ in the context of univariate regression.
%For simplicity, assumption (\ref{centered}) is temporarily removed about $g_j$, that is, the intercept $\mu$ is absorbed into $g_j$.
For concreteness, we consider a fixed function space, say $\mathcal G_1$,
although our discussion is applicable to $\mathcal G_j$ for $j=1,\ldots,p$.
For $r \ge 1$, the $L_r$ norm of a function $f$ on $[0,1]$ is defined as $\|f \|_{L_r}=\{\int_0^1 |f(z)|^r \,\dif z\}^{1/r}$.

\begin{eg}[Sobolev spaces] \label{eg:sspline}
For $r \ge 1$ and $\mym \ge 1$, let $\mathcal W_r^\mym=\mathcal W_r^\mym([0,1])$ be the Sobolev space
of all functions, $g_1: [0,1] \to \mathbb R$,
such that $g_1^{(\mym-1)}$ is absolutely continuous and the norm
$\|g_1 \|_{\mathcal W_r^\mym} = \|g_1\|_{L_r} + \| g_1^{(\mym)} \|_{L_r}$ is finite, where $g_1^{(\mym)}$ denotes the $\mym$th (weak) derivative of $g_1$.
To describe the smoothness, a semi-norm
$\|g_1\|_{F,1} = \| g_1^{(\mym)} \|_{L_r}$ is often used for $g_1 \in \mathcal W_r^\mym$.

In the statistical literature, a major example of Sobolev spaces is $\mathcal W_{2}^{m}=\{ g_1: \|g_1\|_{L_2}+\|g_1^{(\mym)} \|_{L_2} < \infty\}$, which is a
reproducing kernel Hilbert space (e.g., Gu 2002).
Consider a univariate regression model
\begin{align}
Y_i = g_1( X_i^{(1)}) + \varepsilon_i, \quad i=1,\ldots, n. \label{uni-regression}
\end{align}
The Sobolev space $\mathcal W_2^\mym$ is known to lead to polynomial smoothing splines through penalized estimation: there exists a unique solution, in the form of
a spline of order $(2\mym-1)$, when minimizing over $g_1 \in \mathcal W^\mym_2$ the following criterion
\begin{align}
\frac{1}{2n} \sum_{i=1}^n \left\{ Y_i - g_1(X_i^{(1)}) \right\}^2 +  \rho_{n1} \| g_1 \|_{F,1}. \label{uni-pen-est}
\end{align}
This solution can be made equivalent to the standard derivation of smoothing splines, where the penalty in (\ref{uni-pen-est}) is $\rho^\prime_{n1} \| g_1 \|_{F,1}^2$ for a different tuning parameter $\rho^\prime_{n1}$.
Particularly, cubic smoothing splines are obtained with the choice $\mym=2$.
%For additive modeling, the use of Hilbert spaces $W^\mym_2$ has also been studied in both low- and high-dimensional settings (e.g., Chen et~al.~1989; Meier et~al.~2009).
\end{eg}

\begin{eg}[Bounded variation spaces] \label{eg:TVspline}
For a function $f$ on $[0,1]$, the total variation (TV) of $f$ is defined as
\begin{align*}
\TV (f) = \sup \left\{ \sum_{i=1}^k |f( z_{i} ) - f( z_{i-1}) |:  z_0 < z_1 < \ldots < z_k \mbox{ is any partition of } [0,1] \right\}.
\end{align*}
If $f$ is differentiable, then $\TV(f) = \int_0^1 |f^{(1)}(z)| \,\dif z$.
For $\mym \ge 1$, let $\mathcal V^\mym=\mathcal V^\mym([0,1])$ be the bounded variation space that consists of all functions, $g_1: [0,1] \to \mathbb R$,
such that $g_1^{(\mym-2)}$, if $\mym \ge 2$, is absolutely continuous and the norm
$\|g_1 \|_{\mathcal V^\mym} = \|g_1\|_{L_1} + \TV(g_1^{(\mym-1)})$ is finite.
For $g_1 \in \mathcal V^\mym$, the semi-norm $\|g_1\|_{F,1} = \TV(g_1^{(\mym-1)})$ is often used to describe
smoothness.
The bounded variation space $\mathcal V^\mym$ includes as a strict subset the Sobolev space $\mathcal W_1^\mym$, where
the semi-norms also agree: $\TV(g_1^{(\mym-1)}) = \|g_1^{(\mym)} \|_{L_1} $ for $g_1 \in \mathcal W_1^\mym$.

For univariate regression (\ref{uni-regression}) with bounded variation spaces, TV semi-norms can be used as penalties in (\ref{uni-pen-est}) for penalized estimation.
This leads to a class of TV splines, which are shown to adapt well to spatial inhomogeneous smoothness (Mammen \& van de Geer 1997).
For $\mym=1$ or $2$, a minimizer of (\ref{uni-pen-est}) over $g_1 \in \mathcal V^\mym$ can always be chosen as a spline of order $\mym$, with the knots in the set of design points
$\{X_i^{(1)}: i=1,\ldots,n\}$. But, as a complication, this is in general not true for $\mym\ge 3$.

Recently, there is another smoothing method related to TV splines, called trend filtering (Kim et~al.~2009),
where (\ref{uni-pen-est}) is minimized over all possible values $\{g_1(X^{(1)}_i): i=1,\ldots,n\}$ with $\|g_1\|_{F,1}$ replaced by $L_1$ norm
of $\mym$th-order differences of these values. This method is equivalent to TV splines only for $\mym=1$ or $2$. But when the design points are evenly spaced,
it achieves the minimax rate of convergence over functions of bounded variation for general $m\ge 1$, similarly as TV splines (Tibshirani 2014).
\end{eg}

The complexity of a functional class can be described by its metric entropy,
which plays an important role
in the study of empirical processes (van der Vaart \& Wellner 1996).
For a subset $\mathcal F$ in a metric space $\overline{\mathcal F}$ endowed with norm $\|\cdot\|$, the covering number $N( \delta, \mathcal F, \|\cdot \|)$ is defined as the
smallest number of balls of radius $\delta$ in the $\|\cdot\|$-metric needed to cover $\mathcal F$, i.e., the smallest value of $N$ such that
there exist $f_1, \ldots, f_N \in \overline{\mathcal F}$, satisfying $\min_{j=1,\ldots,N} \| f- f_j \| \le \delta$ for any $f \in \mathcal F$.
The entropy of $(\mathcal F, \|\cdot\|)$ is defined as $H(\delta, \mathcal F, \|\cdot\|) = \log N( \delta, \mathcal F, \|\cdot\|)$.
%It is easily shown that $H(\delta, \mathcal F, \|\cdot\|)$ is non-increasing in $\delta >0$.

For analysis of regression models, our approach involves using entropies of functional classes
for empirical norms based on design points, for example,
$\{X^{(1)}_i: i=1,\ldots,n\}$ for subsets of $\mathcal G_1$.
One type of such norms is the empirical $L_2$ norm, %$\|\cdot\|_n $,
$\|g_1\|_n = \{n^{-1} \sum_{i=1}^n g_1^2( X^{(1)}_i )\}^{1/2}$.
Another is the empirical supremum norm, %$\|\cdot\|_{n,\infty}$,
$ \|g_1 \|_{n,\infty} = \max_{i=1,\ldots,n} g_1(X^{(1)}_i )$.
If $\mathcal F$ is the unit ball in the Sobolev space $\mathcal W_r^\mym$ or
the bounded variation space $\mathcal V^\mym$ on $[0,1]$, the general picture is
$H(\delta, \mathcal F, \|\cdot\|)\lesssim \delta^{-1/m}$ for commonly used norms.
See Section~\ref{sect:entropy} for more.

\section{Main results} \label{sect:main-results}

As in Section~\ref{sect:intro}, consider the estimator
\begin{align}\label{estimator}
\hat g = \argmin_{g\in\mathcal G} K_n(g), \quad K_n(g) = \| Y -g\|_n^2/2 +  A_0 R_n( g) ,
\end{align}
where $A_0>1$ is a constant, ${\mathcal G}=\{g=\sum_{j=1}^p g_j:g_j\in {\mathcal G}_j\}$ and
the penalty is of the form
\begin{align*}
R_n(g) = \sum_{j=1}^p R_{nj}(g_j) = \sum_{j=1}^p \Big(\rho_{nj} \|g_j \|_{F,j} + \lambda_{nj} \| g_j \|_n\Big)
\end{align*}
for any decomposition $g=\sum_{j=1}^p g_j$ with $g_j\in {\mathcal G}_j$,
with certain functional penalties $\|f_j\|_{F,j}$ and the empirical $L_2$ penalty $\|f_j\|_n$.
Here the regularization parameters $(\lambda_{nj}, \rho_{nj})$ are of the form
\begin{align*}
\rho_{nj} = \lambda_{nj} w_{nj} ,\quad
\lambda_{nj} = C_1 \left\{ \gamma_{nj} + \sqrt{ \log(p/\epsilon)/n} \right\} ,
\end{align*}
where $C_1>0$ is a noise level depending only on parameters in Assumption~\ref{sub-gaussian-error} below,
$0<\epsilon<1$ is a tail probability for the validity of error bounds, $0<w_{nj}\le 1$ is a rate parameter, and
\begin{align}
\gamma_{nj} & = n^{-1/2} \psi_{nj} ( w_{nj} ) / w_{nj} \label{gamma-def}
\end{align}
for a function $\psi_{nj} (\cdot)$ depending on the entropy of the unit ball of the space
$\mathcal G_j$ under the associated functional penalty.
See Assumption~\ref{ass:fixed-entropy} or \ref{ass:random-entropy} below.

Before theoretical analysis, we briefly comment on computation of $\hat g$.
By standard properties of norms and semi-norms,
the objective function $K_n(g)$ is convex in $g$.
Moreover, there are at least two situations where the infinitely-dimensional problem of minimizing
$K_n( g)$ can be reduced to a finite-dimensional one.
First, if each class $\mathcal G_j$ is a reproducing kernel Hilbert space such as $\mathcal W^\mym_2$,
then a solution $\hat g = \sum_{j=1}^p \hat g_j$ can be obtained such that each $\hat g_j$ is a smoothing
spline with knots in the design points $\{ X^{(j)}_i: i=1,\ldots,n\}$ (e.g., Meier et~al.~2009).
Second, by the following proposition, the optimization problem can be also
reduced to a finite-dimensional one when each class $\mathcal G_j$ is the bounded variation space
$\mathcal V^1$ or $\mathcal V^2$.
As a result, the algorithm in Petersen et~al.~(2016) can be directly used to find $\hat g$ when all classes $(\mathcal G_1,\ldots,\mathcal G_p)$ are $\mathcal V^1$.

\begin{pro} \label{pro:BV-comput}
Suppose that the functional class $\mathcal G_j$ is $\mathcal V^\mym$ for some $1\le j\le p$ and $\mym=1$ or 2. Then a solution $\hat g = \sum_{j=1}^p \hat g_j$
can be chosen such that $\hat g_j$ is piecewise constant with jump points only in $\{X^{(j)}_i:i=1,\ldots,n\}$ if $\mym=1$, or $\hat g_j$  is continuous and piecewise linear
with break points only in $\{X^{(j)}_i:i=1,\ldots,n\}$ if $\mym=2$.
\end{pro}

%By the additive structure of $g \in\mathcal G$ and the penalty, the finite-dimensional reduction is also valid when
%some classes $\mathcal G_j$ are $\mathcal W^\mym_2$ and others are $\mathcal V^1$ or $\mathcal V^2$.
By Example~\ref{eg:TVspline},
it can be challenging to compute $\hat g$ when some classes $\mathcal G_j$ are $\mathcal V^\mym$ with $\mym \ge 3$.
However, this issue may be tackled using trend filtering (Kim et al.~2009) as an approximation.
%We are investigating computational methods and will report elsewhere.

\subsection{Fixed designs} \label{sect:fixed-design}

For fixed designs, the covariates $(X_1,\ldots,X_n)$ are fixed as observed,
whereas $(\varepsilon_1,\ldots,\varepsilon_n)$ and hence $(Y_1,\ldots,Y_n)$ are independent random variables.
The responses are to be predicted when new observations are drawn with covariates from the sample $(X_1,\ldots,X_n)$.
The predictive performance of $\hat g$ is measured by
$\| \hat g - g^*\|_n^2$.

Consider the following three assumptions.
First, we assume sub-gaussian tails for the noises.
This condition can be relaxed, but with increasing technical complexity and possible modification of
the estimators, which we will not pursue here.

\begin{ass}[Sub-gaussian noises] \label{sub-gaussian-error}
Assume that the noises $(\varepsilon_1, \ldots, \varepsilon_n)$ are mutually independent and
uniformly sub-Gaussian: For some constants $D_0>0$ and $D_1 >0$,
\begin{align*}
\max_{i=1,\ldots,n}  D_0\, E\exp ( \varepsilon_i^2 /D_0)\le D_1 .
\end{align*}
We will also impose this assumption for random designs with the interpretation that
the above probability and expectation are taken conditionally on $(X_1,\ldots,X_n)$.
\end{ass}

Second, we impose an entropy condition which describes the relationship between the function
$\psi_{nj} (\cdot)$ in the definition of $\gamma_{nj}$ and the complexity of bounded subsets in $\mathcal G_j$.
Although entropy conditions are widely used to analyze nonparametric regression
(e.g., Section 10.1, van de Geer 2000),
the subset $\mathcal G_j(\delta)$ in our entropy condition below is carefully aligned with
the penalty
$R_{nj}(g_j)=\lambda_{nj}(w_{nj} \|g_j\|_{F,j} + \|g_j\|_n)$.
This leads to a delicate use of maximal inequalities
so as to relax and in some cased remove some restrictions in previous studies of additive models.
See Lemma~\ref{max-ineq2} in the Supplement %Section~\ref{sect:proof-main-thm}
and Raskutti et~al.~(2012, Lemma 1).

\begin{ass}[Entropy condition for fixed designs] \label{ass:fixed-entropy}
For $j=1,\ldots,p$, let $\mathcal G_j (\delta) = \{f_j \in \mathcal G_j: \|f_j\|_{F,j} + \|f_j\|_n/\delta \le 1\}$
and $\psi_{nj} (\delta)$ be an upper bound of the entropy integral as follows:
\begin{align}
\psi_{nj} (\delta) \ge \int_0^{\delta} H^{1/2}(u, \mathcal G_j(\delta), \|\cdot\|_n )\,
\dif u, \quad 0< \delta \le 1. \label{fixed-entropy-ineq}
\end{align}
In general, $\mathcal G_j (\delta)$ and
the entropy $H(\cdot, \mathcal G_j(\delta), \|\cdot\|_n )$
may depend on the design points $\{X_i^{(j)}\}$.
\end{ass}

The third assumption is a compatibility condition, which resembles the restricted eigenvalue condition used in high-dimensional analysis of Lasso in linear regression (Bickel et~al.~2009).
Similar compatibility conditions were used by Meier et~al.~(2009) and Koltchinskii \& Yuan (2010)
in their analysis of penalized estimation in high-dimensional additive regression.

\begin{ass}[Empirical compatibility condition] \label{compat-condition}
For certain subset $S\subset \{1,2,\ldots,p\}$ and constants $\kappa_0>0$ and $\xi_0  > 1$,
assume that
\begin{align*}
\kappa_0^2 \left(\sum_{j\in S} \lambda_{nj} \|f_j\|_n \right)^2 \le \left( \sum_{j\in S} \lambda_{nj}^2 \right) \| f\|^2_n
%\label{compat-ineq2}
\end{align*}
for any functions $\{f_j \in \mathcal G_j: j=1,\ldots, p\}$ and $f=\sum_{j=1}^p f_j \in \mathcal G$ satisfying
\begin{align*}
\sum_{j=1}^p \lambda_{nj} w_{nj} \|f_j\|_{F,j} + \sum_{j \in S^c} \lambda_{nj} \| f_j \|_n \le \xi_0 \sum_{j\in S} \lambda_{nj} \|f_j\|_n. %\label{compat-ineq1}
\end{align*}
\end{ass}

\vspace{.1in}
\begin{rem} \label{emptyset}
The subset $S$ can be different from $\{1\le j\le p: g^*_j\neq 0\}$. In fact, $S$ is arbitrary in the sense that a larger
$S$ leads to a smaller compatibility coefficient $\kappa_0$ which appears as a factor in the denominator
of the ``noise'' term in the prediction error bound below, whereas a smaller $S$ leads to a larger ``bias" term.
Assumption~\ref{compat-condition} is automatically satisfied for the choice $S=\emptyset$.
In this case, it is possible to take $\xi_0 =\infty$ and any $\kappa_0 >0$, provided that we treat summation over an empty set as 0
and $\infty\times 0$ as 0.
\end{rem}

Our main result for fixed designs is an oracle inequality stated in Theorem~\ref{main-thm} below,
where $\bar g =\sum_{j=1}^p \bar g_j\in \mathcal G$ as an estimation target
is an additive function but the true regression function $g^*$ may not be additive.
Denote as a penalized prediction loss
\begin{align*}
\mathcal D_n(\hat g, \bar g) = \frac{1}{2} \|\hat g -g^*\|_n^2 + \frac{1}{2} \|\hat g -\bar g \|_n^2 + (A_0-1) R_n(\hat g - \bar g).
\end{align*}
For a subset $S \subset \{1,2,\ldots,p\}$, write as a bias term for the target $\bar g$
\begin{align*}
\Delta_n(\bar g, S) =  \frac{1}{2} \|\bar g - g^*\|_n^2
+ 2 A_0\left( \hbox{$\sum$}_{j=1}^p  \rho_{nj}  \|\bar g_j\|_{F,j} + \hbox{$\sum$}_{j\in S^c} \lambda_{nj} \| \bar g_j \|_n\right) .
\end{align*}
%In particular, taking $S = \emptyset$ gives
%\begin{align*}
%\Delta_n(\bar g, g^*, \emptyset)= \frac{1}{2} \|\bar g - g^*\|_n^2 + \lambda^2_{n0} + 2 A_0 R_n(\bar g).
%\end{align*}
The bias term is small when $\bar g$ is smooth and sparse and predicts $g^*$ well.

\vspace{.1in}
\begin{thm} \label{main-thm}
Suppose that Assumptions~\ref{sub-gaussian-error}, \ref{ass:fixed-entropy}, and \ref{compat-condition} hold. Then for any $A_0 > (\xi_0+1)/(\xi_0-1)$, we have
with probability at least $1-\epsilon$,
\begin{align}
\mathcal D_n(\hat g, \bar g)
\le \xi_1^{-1} \Delta_n(\bar g, S)
+ 2\xi_2^2  \kappa_0^{-2} \left( \hbox{$\sum$}_{j\in S} \lambda_{nj}^2 \right). \label{main-thm-ii}
\end{align}
where $\xi_1 = 1- 2A_0/\{(\xi_0+1)(A_0-1)\} \in (0,1]$ and $\xi_2 = (\xi_0+1) (A_0-1)$.
\end{thm}

\begin{rem} \label{ANOVA-fixed}
As seen from our proofs, Theorem~\ref{main-thm} and subsequent corollaries are directly applicable
to functional ANOVA modeling, where each function $g_j$ may depend on $X_i^{(j)}$,
a block of covariates,
and the variable blocks are allowed to overlap across different $j$.
The entropy associated with the functional class $\mathcal G_j$
need to be determined accordingly. %This extension also holds for the subsequent corollaries.
\end{rem}

\begin{rem} \label{rem:Bellec}
Using ideas from Bellec and Tsybakov (2016), it is possible to refine the oracle inequality for $\hat g$, such that the scaling parameter $\epsilon$ is fixed, for example, $\epsilon=1/2$
in the definition of $\hat g$ in (\ref{estimator}),
but at any level $0<\tilde\epsilon<1$, (\ref{main-thm-ii}) holds with probability
$1-\tilde\epsilon$ when an additional term of the form $\log(1/\tilde\epsilon)/n$ on the right-hand side.
\end{rem}

Taking $S=\emptyset$ and $\xi_0=\infty$ leads to the following corollary, which explicitly does not require the
compatibility condition (Assumption~\ref{compat-condition}).

\begin{cor} \label{main-cor}
Suppose that Assumptions~\ref{sub-gaussian-error} and \ref{ass:fixed-entropy} hold. Then for any $A_0>1$, we have
with probability at least $1-\epsilon$,
\begin{align}
\mathcal D_n(\hat g, \bar g)
\le \Delta_n(\bar g,\emptyset) = \frac{1}{2} \|\bar g - g^*\|_n^2 + 2 A_0 R_n(\bar g). \label{main-thm-i}
\end{align}
\end{cor}

The following result can be derived from Theorem~\ref{main-thm} through the choice $S = \{1 \le j \le p: \|\bar g_j\|_n > C_0 \lambda_{nj}\}$ for some constant $C_0 >0$.

\begin{cor} \label{main-cor2}
Suppose that Assumptions~\ref{sub-gaussian-error}, \ref{ass:fixed-entropy}, and \ref{compat-condition} hold with
$S = \{1 \le j \le p: \|\bar g_j\|_n > C_0 \lambda_{nj}\}$ for some constant $C_0 >0$. Then for any $0\le q \le 1$ and $A_0> (\xi_0+1)/(\xi_0-1)$, we have
with probability at least $1-\epsilon$,
\begin{align*}
\mathcal D_n(\hat g, \bar g)
%\le \xi_1^{-1} \big\{ \lambda^2_{n0} + \frac{1}{2} \|\bar g - g^*\|_n^2 \big\} +  O(1) \sum_{j=1}^p  \Big\{ \rho_{nj}  \|\bar g_j\|_{F,j} + \min(\lambda_{nj} \| \bar g_j \|_n, C_0 \lambda_{nj}^2) \Big\} \\
\le  O(1) \left\{ \|\bar g - g^*\|_n^2 + \sum_{j=1}^p  \Big( \rho_{nj}  \|\bar g_j\|_{F,j} + \lambda_{nj}^{2-q} \|\bar g_j\|_n^q \Big) \right\},
\end{align*}
where $O(1)$ depends only on $(q,A_0, C_0,\xi_0,\kappa_0)$.
\end{cor}

It is instructive to examine the implications of Corollary~\ref{main-cor2} in a homogenous situation where
for some constants $B_0>0$ and $0<\beta_0<2$,
\begin{align}
\max_{j=1,\ldots,p} \int_0^\delta H^{1/2}(u, \mathcal G_j(\delta), \|\cdot\|_n ) \,\dif u \le B_0 \delta^{1-\beta_0/2}, \quad 0<\delta \le 1. \label{fixed-entropy}
\end{align}
That is, we assume $\psi_{nj}(\delta) = B_0 \delta^{1-\beta_0/2}$ in (\ref{fixed-entropy-ineq}).
For $j=1,\ldots,p$, let
\begin{align}
w_{nj}=w_n(q)= \{\gamma_n(q)\}^{1-q}, \quad \gamma_{nj} = \gamma_n (q)= B_0^{\frac{2}{2+\beta_0(1-q)}} n^{ \frac{-1}{2+\beta_0(1-q)} } , \label{fixed-rate}
\end{align}
which are determined by balancing the two rates $\rho_{nj} = \lambda_{nj}^{2-q}$, that is, $w_{nj} = \lambda_{nj}^{1-q}$, along
with the definition $\gamma_{nj} = B_0 n^{-1/2} w_{nj}^{-\beta_0/2}$ by (\ref{gamma-def}).
For $g=\sum_{j=1}^p g_j\in\mathcal G$, denote $\|g\|_{F,1} = \sum_{j=1}^p \|g_j\|_{F,j}$ and $\|g\|_{n,q} = \sum_{j=1}^p \|g_j\|_n^q$.
For simplicity, we also assume that $g^*$ is an additive function and set $\bar g = g^*$ for Corollary~\ref{main-cor3}.

\begin{cor} \label{main-cor3}
Assume that (\ref{additive-reg}) holds and $\|g^*\|_{F,1} \le C_1 M_F$ and $\|g^*\|_{n,q} \le C_1^q M_q$ for $0\le q \le 1$, $M_q>0$, and $M_F>0$, possibly depending on $(n,p)$.
In addition, suppose that (\ref{fixed-entropy}) and (\ref{fixed-rate}) hold,
and Assumptions~\ref{sub-gaussian-error} and \ref{compat-condition} are satisifed with
$S = \{1 \le j \le p: \|g_j^*\|_n > C_0 \lambda_{nj}\}$ for some constant $C_0 >0$.
If $0< w_n(q) \le 1$ for sufficiently large $n$, then for any $A_0> (\xi_0+1)/(\xi_0-1)$, we have
with probability at least $1-\epsilon$,
\begin{align}
& \mathcal D_n(\hat g, g^*) =  \|\hat g -g^*\|_n^2  + (A_0-1) R_n(\hat g - g^*) \nonumber \\
& \le O(1) C_1^2 (M_F+M_q) \left\{ \gamma_n(q) + \sqrt{\log(p/\epsilon)/n}\right\}^{2-q} , \label{main-cor3-equ}
\end{align}
where $O(1)$ depends only on $(q, A_0, C_0,\xi_0,\kappa_0)$.
\end{cor}

\begin{rem} \label{rem-main-cor3}
There are several interesting features in the convergence rate (\ref{main-cor3-equ}).
First, (\ref{main-cor3-equ}) presents a spectrum of convergence rates in the form
$$
\left\{ n^{\frac{-1}{2+\beta_0(1-q)}} + \sqrt{\log(p)/n}\right\}^{2-q},
$$
which are easily shown to become slower as $q$ increases from 0 to 1,
that is, the exponent $(2-q)/\{2+\beta_0(1-q)\}$ is decreasing in $q$ for $0<\beta_0<2$.
The rate (\ref{main-cor3-equ}) gives the slow rate $\{\log(p)/n\}^{1/2}$ for $q=1$,
or the fast rate $n^{\frac{-2}{2+\beta_0}} +\log(p)/n$ for $q=0$,
as previously obtained for additive regression with reproducing kernel Hilbert spaces.
We defer to Section~\ref{sect:specific-result} the comparison
with existing results in random designs.
Second, the rate (\ref{main-cor3-equ}) is in general at least as fast as
$$
\left\{ n^{\frac{-1}{2+\beta_0}} + \sqrt{\log(p)/n}\right\}^{2-q} .
$$
Therefore, weaker sparsity (larger $q$) leads to a slower rate of convergence, but not as slow as the fast rate $\{n^{\frac{-2}{2+\beta_0}} +\log(p)/n\}$ raised to the power of $(2-q)/2$.
This is in contrast with previous results on penalized estimation over $L_q$ sparsity balls, for example,
the rate $\{k/n +\log(p)/n\}^{(2-q)/2}$ obtained for group Lasso estimation in linear regression (Neghaban et~al.~2012), where $k$ is the group size.
Third, the rate (\ref{main-cor3-equ}) is in general not as fast as the following rate (unless $q=0$ or 1)
$$
n^{\frac{-2}{2+\beta_0}} + \{\log(p)/n\}^{(2-q)/2} ,
$$
which was obtained by Yuan \& Zhou (2016) using constrained least squares
for additive regression with reproducing kernel Hilbert spaces under an $L_q$ ball in the Hilbert norm:
$\sum_{j=1}^p \|g^*_j\|_H^q \le M_q$.
This difference can be explained by the fact that an $L_q$ ball in $\|\cdot\|_H$ norm is more restrictive than in $\|\cdot\|_n$ or $\|\cdot\|_Q$ norm for our results.
\end{rem}

\subsection{Random designs} \label{sect:random-design}

For random designs, prediction of the responses can be sought when new observations are randomly drawn with covariates
from the distributions of $(X_1,\ldots,X_n)$, instead of within the sample $(X_1,\ldots,X_n)$ as in Section~\ref{sect:fixed-design}.
For such out-of-sample prediction, the performance of $\hat g$ is measured by
$\| \hat g - g^*\|_Q^2$, where $\|\cdot\|_Q$ denotes the theoretical norm: $\|f\|_Q^2 =n^{-1} \sum_{i=1}^n E\{f^2(X_i)\}$ for a function $f(x)$.

Consider the following two extensions of Assumptions~\ref{ass:fixed-entropy} and \ref{compat-condition},
such that dependency on the empirical norm $\|\cdot\|_n$ and hence on $(X_1,\ldots,X_n)$ are removed.

\begin{ass}[Entropy condition for random designs] \label{ass:random-entropy}
For some constant $0<\eta_0<1$ and $j=1,\ldots,p$, let $\psi_{nj} (\delta)$ be an upper bound of the entropy integral, independent of the realizations
$\{X_i^{(j)}: i=1,\ldots,n\}$, as follows:
\begin{align}
\psi_{nj} (\delta) \ge \int_0^{\delta} H^{*1/2}( (1-\eta_0)u, \mathcal G^*_j(\delta), \|\cdot\|_n ) \,\dif u, \quad 0<\delta\le 1,\label{random-entropy-ineq}
\end{align}
where $\mathcal G^*_j (\delta) = \{f_j \in \mathcal G_j: \|f_j\|_{F,j} + \|f_j\|_Q /\delta \le 1\}$ and
\begin{align*}
H^{*}( u, \mathcal G^*_j(\delta), \|\cdot\|_n ) = \sup_{(X_1^{(j)},\ldots,X_n^{(j)})} H(u, \mathcal G^*_j (\delta),\|\cdot\|_n).
\end{align*}
\end{ass}

\begin{ass}[Theoretical compatibility condition] \label{compat-condition-Q}
For some subset $S\subset \{1,2,\ldots,p\}$ and constants $\kappa^*_0>0$ and $\xi^*_0 > 1$, assume that for
any functions $\{f_j \in \mathcal G_j: j=1,\ldots, p\}$ and $f=\sum_{j=1}^p f_j\in\mathcal G$, if
\begin{align}
\sum_{j=1}^p \lambda_{nj} w_{nj} \|f_j\|_{F,j} + \sum_{j \in S^c} \lambda_{nj} \| f_j \|_Q \le \xi^*_0 \sum_{j\in S} \lambda_{nj} \|f_j\|_Q, \label{compat-Q-equ1}
\end{align}
then
\begin{align}
\kappa_0^{*2} \left(\sum_{j\in S} \lambda_{nj} \|f_j\|_Q \right)^2 \le \left( \sum_{j\in S} \lambda_{nj}^2 \right) \| f\|^2_Q  . \label{compat-Q-equ2}
\end{align}
\end{ass}

\begin{rem} \label{emptyset-Q}
Similarly as in Remark~\ref{emptyset} about the empirical compatibility condition,
Assumption~\ref{compat-condition-Q} is also automatically satisfied for the choice $S=\emptyset$,
in which case it is possible to take $\xi^*_0 =\infty$ and any $\kappa^*_0 >0$.
\end{rem}

To tackle random designs, our approach relies on
establishing appropriate convergence of empirical norms $\|\cdot\|_n$ to $\|\cdot\|_Q$ uniformly over the
space of additive functions $\mathcal G$, similarly as in Meier et~al.~(2009) and Koltchinskii \& Yuan (2010).
For clarity, we postulate the following assumption on the rate of such convergence to develop general analysis of $\hat g$.
We will study convergence of empirical norms specifically for Sobolev and bounded variation spaces in Section~\ref{sect:conv-norm},
and then provide corresponding results on the performance of  $\hat g$ in Section~\ref{sect:specific-result}.
For $g=\sum_{j=1}^p g_j \in \mathcal G$, denote
\begin{align*}
& R^*_n(g) = \sum_{j=1}^p  R^*_{nj} ( g_j), \quad
R^*_{nj} ( g_j) = \lambda_{nj} (w_{nj}\|g_j\|_{F,j} + \|g_j\|_Q ) ,
\end{align*}
as the population version of the penalty $R_n(g)$, with $\|g_j\|_Q$ in place of $\|g_j\|_n$.

\begin{ass}[Convergence of empirical norms] \label{rate-assumption}
Assume that
\begin{align}
P\left\{ \sup_{g\in \mathcal G} \frac{ \left| \|g\|_n^2 - \|g \|_Q^2 \right|}{ R^{*2}_n (g) } > \phi_n \right\} \le \pi , \label{norm-conv-prob}
\end{align}
where $0 < \pi <1$ and $\phi_n >0$ such that for sufficiently large $n$, one or both of the following conditions are valid.
\begin{itemize}
\item[(i)] $\phi_n (\max_{j=1,\ldots,p}\lambda^2_{nj} ) \le \eta_0^2$, where $\eta_0$ is from Assumption~\ref{ass:random-entropy}.

\item[(ii)] For some constant $0 \le \eta_1 < 1$, we have
\begin{align}
\phi_n\, (\xi_0^*+1)^2 \kappa_0^{*-2} \left(\sum_{j\in S} \lambda_{nj}^2 \right) \le \eta_1^2 ,\label{rate-assumption-ii}
\end{align}
where $S$ is the subset of $\{1,2,\ldots,p\}$ used in Assumption~\ref{compat-condition-Q}.
\end{itemize}
\end{ass}

Our main result, Theorem~\ref{main2-thm}, gives an oracle inequality for random designs, where the predictive performance of $\hat g$ is compared with that of an arbitrary additive function
$\bar g = \sum_{j=1}^p \bar g_j \in \mathcal G$, but the true regression function $g^*$ may not be additive, similarly as in Theorem~\ref{main-thm} for fixed designs.
For a subset $S \subset \{1,2,\ldots,p\}$, denote
\begin{align*}
\Delta^*_n(\bar g, S) = \frac{1}{2} \|\bar g - g^*\|_n^2 + 2 A_0(1-\eta_0) \left( \sum_{j=1}^p  \rho_{nj}  \|\bar g_j\|_{F,j} + \sum_{j\in S^c} \lambda_{nj} \| \bar g_j \|_Q \right) ,
\end{align*}
which, unlike $\Delta_n(\bar g, S)$, involves $\| \bar g_j\|_Q$ and $\eta_0$ from Assumptions~\ref{ass:random-entropy} and \ref{rate-assumption}(i).
%In particular, taking $S = \emptyset$ gives
%\begin{align*}
%\Delta^*_n(\bar g, g^*, \emptyset)= \frac{1}{2} \|\bar g - g^*\|_n^2 + \lambda^2_{n0} + 2 A_0 R(\bar g).
%\end{align*}

\begin{thm} \label{main2-thm}
Suppose that Assumptions~\ref{sub-gaussian-error}, \ref{ass:random-entropy}, \ref{compat-condition-Q} and \ref{rate-assumption}(i)--(ii) hold with $0 < \eta_0 < (\xi_0^*-1)/(\xi^*_0+1)$.
Let $A(\xi^*_0,\eta_0) =  \{\xi^*_0+1 +\eta_0(\xi^*_0 +1)\} / \{\xi^*_0-1 - \eta_0(\xi^*_0+1)\} >(1+\eta_0)/(1-\eta_0)$.
Then for any $A_0 >A(\xi^*_0,\eta_0)$, we have with probability at least $1-\epsilon-\pi$,
 \begin{align}
& \frac{1}{2} \|\hat g -g^*\|_n^2 + \frac{1}{2} \|\hat g -\bar g \|_n^2 + (1-\eta_1) A_1 R^*_n(\hat g - \bar g) \nonumber \\
& \le \xi_1^{*-1} \Delta^*_n(\bar g, S) +  2 \xi_2^{*2} \kappa_0^{*-2} \left( \sum_{j\in S} \lambda_{nj}^2 \right) , \label{main2-thm-ii1}
\end{align}
where $A_1=(A_0-1)-\eta_0(A_0+1)>0$, $\xi_1^* = 1- 2 A_0 / \{(\xi_0^*+1)A_1\} \in (0,1]$ and $ \xi_2^* = (\xi_0^*+1) A_1$.
Moreover, we have with probability at least $1-\epsilon-\pi$,
\begin{align}
& \mathcal D^*_n (\hat g, \bar g) := \frac{1}{2} \|\hat g -g^*\|_n^2 + \frac{1}{2} \|\hat g - \bar g\|_Q^2 + A_2 R^*_n(\hat g - \bar g) \nonumber \\
& \le  \xi_3^{*-1} \Delta^*_n(\bar g, S) +  2 \xi_4^{*2} \kappa_0^{*-2} \left( \sum_{j\in S} \lambda_{nj}^2 \right)  +
\frac{\phi_n}{ 2A_1 A_2} \xi_3^{*-2} \Delta^{*2}_n(\bar g, g^*, S), \label{main2-thm-ii2}
\end{align}
where $ A_2 = A_1 / (1 - \eta_1^2)$, $ \xi^*_3 = \xi^*_1 ( 1- \eta_1^2) $, and $ \xi^*_4 = \xi^*_2 /( 1- \eta_1^2) $.
\end{thm}

\vspace{.1in}
\begin{rem} \label{ANOVA-random}
Similarly as in Remark~\ref{ANOVA-fixed}, we emphasize that Theorem~\ref{main2-thm} and subsequent corollaries are also applicable to functional ANOVA modeling (e.g., Gu 2002).
For example, consider model (\ref{additive-reg}) studied in Yang \& Tokdar (2015), where each $g^*_j$ is assumed to depend only on $d_0$ of a total of $d$ covariates
and lie in a H\"{o}lder space with smoothness level $\alpha_0$. Then $p=\binom{d}{d_0}$, and the entropy condition (\ref{random-entropy}) holds with $\beta_0 = d_0/\alpha_0$.
Under certain additional conditions, Corollary~\ref{main2-cor3} with $q=0$ shows that
penalized estimation studied here achieves a convergence rate $M_0 n^{\frac{-2}{2+\beta_0}} + M_0 \log(p)/n$ under exact sparsity of size $M_0$, where $n^{\frac{-2}{2+\beta_0}}$ is
the rate for estimation of a single regression function in the H\"{o}lder class in dimension $d_0$ with smoothness $\beta_0^{-1}$, and $\log(p)/n \asymp d_0 \log(d/d_0)/n$ is the term associated with handling $p$ regressors.
This result agrees with the minimax rate derived in Yang \& Tokdar (2015), but can be applied when more general functional classes are used such as multi-dimensional Sobolev spaces.
%In addition, the analysis in Yang \& Tokdar (2015) seems to assume that all the $d$ covariates are independent of each other (as seen from their Theorem 6.4),
%whereas the theoretical compatibility condition used in our analysis is in general much weaker.
In addition, Yang \& Todkar (2015) considered adaptive Bayes estimators which are nearly minimax with some extra logarithmic factor in $n$.
\end{rem}

Taking $S=\emptyset$, $\xi^*_0=\infty$, and $\eta_1=0$  leads to the following corollary,
which explicitly does not require the theoretical compatibility condition (Assumption~\ref{compat-condition-Q}) or the rate condition, Assumption~\ref{rate-assumption}(ii),
for convergence of empirical norms.

\begin{cor} \label{main2-cor}
Suppose that Assumptions~\ref{sub-gaussian-error}, \ref{ass:random-entropy}, and \ref{rate-assumption}(i) hold.
Then for any $A_0> (1+\eta_0)/(1-\eta_0)$, we have with probability at least $1-\epsilon-\pi$,
\begin{align}
& \frac{1}{2} \|\hat g -g^*\|_n^2 + \frac{1}{2} \|\hat g -\bar g \|_n^2 + A_1 R^*_n(\hat g - \bar g) \nonumber \\
& \le \Delta^*_n(\bar g, \emptyset)= \lambda^2_{n0} + \frac{1}{2} \|\bar g - g^*\|_n^2 + 2 A_0 R^*_n(\bar g). \label{main2-thm-i1}
\end{align}
Moreover, we have with probability at least $1-\epsilon-\pi$,
\begin{align}
& \frac{1}{2} \|\hat g -g^*\|_n^2 + \frac{1}{2} \|\hat g -\bar g \|_Q^2 + A_1 R^*_n(\hat g - \bar g) \nonumber \\
& \le \Delta^*_n(\bar g, \emptyset) + \frac{\phi_n}{2A_1^2} \Delta^{*2}_n(\bar g, \emptyset). \label{main2-thm-i2}
\end{align}
\end{cor}

\vspace{.1in}
The preceding results deal with both in-sample and out-of-sample prediction.
For space limitation, except in Proposition~\ref{cor-slow},
we hereafter focus on the more challenging out-of-sample prediction.
Under some rate condition about $\phi_n$ in (\ref{norm-conv-prob}), the additional term involving $\phi_n \Delta_n^{*2}(\bar g, S)$ can be absorbed into the first term,
as shown in the following corollary. Two possible scenarios are accommodated.
On one hand, taking $\bar g= g^*$ directly gives high-probability bounds on the prediction error $\|\hat g - g^*\|_Q^2$
provided that $g^*$ is additive, that is, model (\ref{additive-reg}) is correctly specified.
On the other hand, the error $\|\hat g - g^*\|_Q^2$ can also be bounded, albeit in probability, in terms of an arbitrary additive function $\bar g \in \mathcal G$,
while allowing $g^*$ to be non-additive.

\begin{cor} \label{main2-cor2}
Suppose that the conditions of Theorem~\ref{main2-thm} hold with $S = \{1 \le j \le p: \|\bar g_j\|_Q > C^*_0 \lambda_{nj}\}$ for some constant $C^*_0>0$,
and (\ref{norm-conv-prob}) holds with $\phi_n>0$ also satisfying
\begin{align}
\phi_n \left( \sum_{j=1}^p  \rho_{nj}  \|\bar g_j\|_{F,j} + \sum_{j\in S^c} \lambda_{nj} \| \bar g_j \|_Q \right) \le \eta_2 ,  \label{main2-cor2-equ}
\end{align}
for some constant $\eta_2>0$. Then for any $0\le q\le 1$ and $A_0 > A(\xi^*_0,\eta_0)$, we have with probability at least $1-\epsilon-\pi$,
\begin{align*}
& \mathcal D^*_n (\hat g, \bar g) \le \{O(1)+\phi_n \|\bar g- g^*\|_n^2\} \left\{ \|\bar g- g^*\|_n^2 + \sum_{j=1}^p \Big( \rho_{nj}  \|\bar g_j\|_{F,j} + \lambda_{nj}^{2-q} \|\bar g_j\|_Q^q \Big) \right\},
\end{align*}
where $O(1)$ depends only on $(q, A^*_0, C^*_0,\xi^*_0,\kappa^*_0,\eta_0,\eta_1, \eta_2)$. In addition, suppose that
$\phi_n \|\bar g - g^*\|_Q^2 $ is bounded by a constant and $\epsilon=\epsilon(n,p)$ tends to $0$ in the definition of $\hat g$ in (\ref{estimator}).
Then for any $0 \le q \le 1$, we have
\begin{align*}
& \|\hat g - g^* \|_Q^2
\le O_p(1)\left\{ \| \bar g- g^*\|_Q^2 +\sum_{j=1}^p \Big( \rho_{nj}  \|\bar g_j\|_{F,j} + \lambda_{nj}^{2-q} \|\bar g_j\|_Q^q \Big)\right\} .
\end{align*}
\end{cor}

\vspace{.1in}
Similarly as Corollary~\ref{main-cor3}, it is useful to deduce the following result in a homogeneous situation where
for some constants $B^*_0>0$ and $0<\beta_0<2$,
\begin{align}
\max_{j=1,\ldots,p} \int_0^\delta H^{*1/2}( (1-\eta_0) u, \mathcal G^*_j(\delta), \|\cdot\|_n ) \,\dif u \le B^*_0 \delta^{1-\beta_0/2}, \quad 0<\delta \le 1. \label{random-entropy}
\end{align}
That is, we assume $\psi_{nj}(\delta) = B^*_0 \delta^{1-\beta_0/2}$ in (\ref{random-entropy-ineq}). For $j=1,\ldots,p$, let
\begin{align}
& w_{nj}=w^*_n(q) =\max\{\gamma_n(q)^{1-q},\, \nu_n^{1-q}\}, \label{random-rate1} \\
& \gamma_{nj} = \gamma^*_n(q) = \min\{\gamma_n(q),\, B_0^* n^{-1/2} \nu_n^{-(1-q)\beta_0/2}\},\label{random-rate2}
\end{align}
where $\nu_n=\{\log(p/\epsilon)/n\}^{1/2}$, and $w_n(q) = \gamma_n(q)^{1-q}$ and
\begin{align*}
\gamma_n(q) = {B_0^*}^{\frac{2}{2+\beta_0(1-q)}} n^{ \frac{-1}{2+\beta_0(1-q)} } \asymp n^{ \frac{-1}{2+\beta_0(1-q)} }
\end{align*}
are determined from the relationship (\ref{gamma-def}), that is, $\gamma_n(q) = B^*_0 n^{-1/2} w_n(q)^{-\beta_0/2}$.
The reason why $(w^*_n(q), \gamma^*_n(q))$ are used instead of the simpler choices $(w_n(q), \gamma_n(q))$ is
that the rate condition (\ref{main2-cor3-equ}) needed below would become stronger if $\gamma^*_n(q)$ were replaced by $\gamma_n(q)$.
The rate of convergence, however, remains the same even if $\gamma^*_n(q)$ is substituted for $\gamma_n(q)$ in (\ref{main2-cor3-rate}).
See Remark~\ref{rem:cor-fast3} for further discussion.
For $g=\sum_{j=1}^p g_j\in\mathcal G$, denote $\|g\|_{F,1} = \sum_{j=1}^p \|g_j\|_{F,j}$ and $\|g\|_{Q,q} = \sum_{j=1}^p \|g_j\|_Q^q$.

\begin{cor} \label{main2-cor3}
Assume that (\ref{additive-reg}) holds and $\|g^*\|_{F,1} \le C_1 M_F$ and $\|g^*\|_{Q,q} \le C_1^q M_q$ for $0\le q \le 1$, $M_q>0$, and $M_F>0$, possibly depending on $(n,p)$.
In addition, suppose that (\ref{random-entropy}), (\ref{random-rate1}), and (\ref{random-rate2}) hold,
Assumptions~\ref{sub-gaussian-error}, \ref{compat-condition-Q}, and \ref{rate-assumption}(i)
are satisfied with $0 < \eta_0 < (\xi_0^*-1)/(\xi^*_0+1)$ and $S = \{1 \le j \le p: \|g^*_j\|_Q > C^*_0 \lambda_{nj}\}$ for some constant $C^*_0>0$,
and (\ref{norm-conv-prob}) holds with $\phi_n>0$ satisfying
\begin{align}
\phi_n C_1^2 (M_F+M_q) \left\{ \gamma^*_n(q) + \sqrt{\log(p/\epsilon)/n}\right\}^{2-q} = o(1). \label{main2-cor3-equ}
\end{align}
Then for sufficiently large $n$, depending on $(M_F, M_q)$ only through the convergence rate in (\ref{main2-cor3-equ}), and any $A_0> A(\xi^*_0,\eta_0)$, we have with probability at least $1-\epsilon-\pi$,
\begin{align}
\mathcal D^*_n (\hat g, g^*)
\le O(1) C_1^2 (M_F+M_q) \left\{ \gamma_n(q) + \sqrt{\log(p/\epsilon)/n}\right\}^{2-q} , \label{main2-cor3-rate}
\end{align}
where $O(1)$ depends only on $(q, A^*_0, C^*_0,\xi^*_0,\kappa^*_0,\eta_0)$.
\end{cor}

In the case of $q\not=0$, Corollary~\ref{main2-cor3} can be improved by relaxing the rate condition (\ref{main2-cor3-equ})
but requiring the following compatibility condition, which is seemingly stronger than Assumption~\ref{compat-condition-Q},
and also more aligned with those used in related works on additive regression (Meier et~al.~2009; Koltchinskii \& Yuan 2010).

\begin{ass}[Monotone compatibility condition] \label{compat-mono}
For some subset $S\subset \{1,2,\ldots,p\}$ and constants $\kappa^*_0>0$ and $\xi^*_0 > 1$, assume that for
any functions $\{f_j \in \mathcal G_j: j=1,\ldots, p\}$ and $f=\sum_{j=1}^p f_j\in\mathcal G$, if (\ref{compat-Q-equ1}) holds then
\begin{align}
\kappa_0^{*2} \sum_{j\in S} \|f_j\|_Q^2 \le \| f\|^2_Q  . \label{compat-Q-equ3}
\end{align}
\end{ass}

\begin{rem} \label{mono-S}
By the Cauchy--Schwartz inequality,
(\ref{compat-Q-equ3}) implies (\ref{compat-Q-equ2}), and hence
Asssumption~\ref{compat-mono} is stronger than Assumption~\ref{compat-condition-Q}.
However, there is a monotonicity in $S$ for the validity of Assumption~\ref{compat-mono} with (\ref{compat-Q-equ3}) used.
In fact, for any subset $S^\prime \subset S$ and
any functions $\{f^\prime_j \in \mathcal G_j: j=1,\ldots, p\}$ and $f^\prime=\sum_{j=1}^p f^\prime_j\in\mathcal G$, if
\begin{align*}
\sum_{j=1}^p \lambda_{nj} w_{nj} \|f^\prime_j\|_{F,j} + \sum_{j \in S^{\prime c}} \lambda_{nj} \| f^\prime_j \|_Q \le \xi^*_0 \sum_{j\in S^\prime} \lambda_{nj} \|f^\prime_j\|_Q,
\end{align*}
then (\ref{compat-Q-equ1}) holds with $f_j=f^\prime_j$, $j=1,\ldots,p$, and hence, via (\ref{compat-Q-equ3}),  implies
\begin{align*}
\| f^\prime \|^2_Q \ge \kappa_0^{*2} \sum_{j\in S} \|f^\prime_j\|_Q^2 \ge \kappa_0^{*2} \sum_{j\in S^\prime} \|f^\prime_j\|_Q^2.
\end{align*}
Therefore, if Assumption~\ref{compat-mono} holds for a subset $S$, then it also holds for any subset $S^\prime \subset S$ with the same constants $(\xi^*_0, \kappa^*)$.
\end{rem}

\begin{cor} \label{main2-cor4}
Suppose that the conditions of Corollary~\ref{main2-cor3} are satisfied with $0<q \le 1$ (excluding $q=0$),
Assumption~\ref{compat-mono} holds instead of Assumption~\ref{compat-condition-Q}, and the following condition holds instead of (\ref{main2-cor3-equ}),
\begin{align}
\phi_n C_1^2 (M_F+M_q) \left\{ \gamma^*_n(q) + \sqrt{\log(p/\epsilon)/n}\right\}^{2-q} \le \eta_3, \label{main2-cor4-equ}
\end{align}
for some constant $\eta_3>0$.
If $0<w^*_n(q) \le 1$ for sufficiently large $n$, then for any $A_0> A(\xi^*_0,\eta_0)$, inequality (\ref{main2-cor3-rate}) holds with probability at least $1-\epsilon-\pi$,
where $O(1)$ depends only on $(q, A^*_0, C^*_0,\xi^*_0,\kappa^*_0,\eta_0,\eta_3)$.
\end{cor}

To demonstrate the flexibility of our approach and compare with related results, notably Suzuki \& Sugiyama (2013),
we provide another result in the context of Corolloary~\ref{main2-cor3} with $(w_{nj},\gamma_{nj})$ allowed to
depend on $(M_F,M_q)$, in contrast with the choices (\ref{random-rate1})--(\ref{random-rate2}) independent of $(M_F,M_q)$.
For $j=1,\ldots,p$, let
\begin{align}
& w_{nj}=w^\dag_n(q) =\max\{w^\prime_n(q),\, \nu_n^{1-q} (M_q/M_F)\}, \label{random-rate3} \\
& \gamma_{nj} = \gamma^\dag_n(q) = \min\{\gamma^\prime_n(q),\, B_0^* n^{-1/2} \nu_n^{-(1-q)\beta_0/2}(M_q/M_F)^{-\beta_0/2} \}, \label{random-rate4}
\end{align}
where $w^\prime_n(q) = \gamma^\prime_n(q)^{1-q} (M_q/M_F)$ and $\gamma^\prime_n(q) = {B_0^*}^{\frac{2}{2+\beta_0(1-q)}} n^{\frac{-1}{2+\beta_0(1-q)}} (M_q / M_F)^{\frac{-\beta_0}{2+\beta_0(1-q)}}$
are determined along with the relationship $\gamma^\prime_n(q) = B^*_0 n^{-1/2} w^\prime_n(q)^{-\beta_0/2}$ by (\ref{gamma-def}).
These choices are picked to balance the two rates: $\lambda_n w_n M_F$ and  $\lambda_n^{2-q}  M_q$,
where $w_n$ and $\lambda_n$  denote the common values of $w_{nj}$ and $\lambda_{nj}$ for $j=1,\ldots,p$.

\begin{cor} \label{main2-cor5}
Suppose that the conditions of Corollary~\ref{main2-cor3} are satisfied
except that $(w_{nj},\gamma_{nj})$ are defined by (\ref{random-rate3})--(\ref{random-rate4}),
and the following condition holds instead of (\ref{main2-cor3-equ}),
\begin{align}
\phi_n C_1^2 M_q \left\{ \gamma^\dag_n(q) + \sqrt{\log(p/\epsilon)/n}\right\}^{2-q} = o(1). \label{main2-cor5-equ}
\end{align}
Then for sufficiently large $n$, depending on $(M_F, M_q)$ only through the convergence rate in (\ref{main2-cor5-equ}), and any $A_0> A(\xi^*_0,\eta_0)$, we have with probability at least $1-\epsilon-\pi$,
\begin{align}
\mathcal D^*_n (\hat g, g^*)
\le O(1) C_1^2 \left\{ M_q^{\frac{2-\beta_0}{2+\beta_0(1-q)}} M_F^{\frac{(2-q)\beta_0}{2+\beta_0(1-q)}} n^{\frac{-(2-q)}{2+\beta_0(1-q)}} + M_q \nu_n^{2-q} \right\} , \label{main2-cor5-rate}
\end{align}
where $O(1)$ depends only on $(q, B^*_0, A^*_0, C^*_0,\xi^*_0,\kappa^*_0,\eta_0)$.
\end{cor}

\begin{rem} \label{rem-main2-cor5}
In the special case of $q=0$ (exact sparsity), the convergence rate (\ref{main2-cor5-rate}) reduces to
$M_0^{\frac{2-\beta_0}{2+\beta_0}} M_F^{\frac{2\beta_0}{2+\beta_0}} n^{\frac{-2}{2+\beta_0}} + M_0\nu_n^{2} $.
The same rate was obtained in Suzuki \& Sugiyama (2013) under
\begin{align}
\sum_{j=1}^p \|g^*_j\|_Q^0 \le M_0, \quad \sum_{j=1}^p \|g^*_j \|_H  \le M_F\le cM_0, \label{SuSu}
\end{align}
with a constant $c$ for additive regression with reproducing kernel Hilbert spaces, where $\|g^*_j\|_H$ is the Hilbert norm.
As one of their main points, this rate was argued to be faster than $(M_0+M_F)n^{\frac{-2}{2+\beta_0}} + M_0 \nu_n^{2}$, that is,
the rate (\ref{main2-cor3-rate}) with $q=0$ under (\ref{SuSu}).
Our analysis sheds new light on the relationship between the rates (\ref{main2-cor3-rate}) and (\ref{main2-cor5-rate}): their difference
mainly lies in whether the tuning parameters $(w_{nj},\gamma_{nj})$ are chosen independently of $(M_F,M_0)$ or depending on $(M_F,M_0)$.
\end{rem}

\subsection{Convergence of empirical norms} \label{sect:conv-norm}

We provide two explicit results on the convergence of empirical norms as needed for Assumption~\ref{rate-assumption}.
These results can also be useful for other applications.

Our first result, Theorem~\ref{thm-Sobolev}, is applicable (but not limited) to Sobolev and bounded variation spaces in general.
For clarity, we postulate another entropy condition, similar to Assumption~\ref{ass:random-entropy} but with the empirical supremum norms.

\begin{ass}[Entropy condition in supremum norms] \label{ass:supre-entropy}
For $j=1,\ldots,p$, let $\psi_{nj,\infty} (\cdot,\delta)$ be an upper envelope of the entropy integral, independent of the realizations
$\{X_i^{(j)}: i=1,\ldots,n\}$, as follows:
\begin{align*}
\psi_{nj,\infty} (z, \delta) \ge \int_0^{z} H^{*1/2}(u/2, \mathcal G^*_j (\delta), \|\cdot\|_{n,\infty} ) \,\dif u, \quad z>0,\, 0<\delta\le 1,
\end{align*}
where $\mathcal G^*_j(\delta) = \{f_j \in \mathcal G_j: \|f_j\|_{F,j} + \|f_j\|_Q /\delta \le 1\}$ as in Assumption~\ref{ass:random-entropy} and
\begin{align*}
H^{*}( u, \mathcal G^*_j(\delta), \|\cdot\|_{n,\infty} ) = \sup_{(X_1^{(j)},\ldots,X_n^{(j)})} H(u, \mathcal G^*_j (\delta),\|\cdot\|_{n,\infty}).
\end{align*}
\end{ass}

We also make use of the following two conditions about metric entropies and sup-norms.
Suppose that for $j=1,\ldots,p$, $\psi_{nj}(\delta)$ and
$\psi_{nj,\infty}(z,\delta)$ in Assumptions~\ref{ass:random-entropy} and \ref{ass:supre-entropy} are in the polynomial forms
\begin{align}
& \psi_{nj}(\delta) = B_{nj} \delta^{1-\beta_j/2}, \quad 0 < \delta \le 1 , \label{entropy-cond1} \\
& \psi_{nj,\infty}(z,\delta) = B_{nj,\infty} z^{1-\beta_j/2}, \quad z>0,\, 0<\delta\le 1, \label{entropy-cond2}
\end{align}
where $0<\beta_j<2$ is a constant, and $B_{nj} >0$ and $B_{nj,\infty} >0$ are constants, possibly depending on $n$.
Denote $\Gamma_{n} = \max_{j=1,\ldots,p} (B_{nj,\infty}/B_{nj})$.
In addition, suppose that for $j=1,\ldots,p$,
\begin{align}
\| g_j\|_\infty \le (C_{4,j}/2) \big( \|g_j\|_{F,j} + \|g_j\|_{Q} \big)^{\tau_j} \|g_j\|_{Q}^{1-\tau_j}, \quad  g_j\in\mathcal G_j,\label{sup-norm-cond}
\end{align}
where $C_{4,j} \ge 1$ and $0<\tau_j \le (2/\beta_j-1)^{-1}$ are constants.
Let $\gamma_{nj} = n^{-1/2} \psi_{nj}(w_{nj}) /w_{nj} = n^{-1/2} B_{nj} w_{nj}^{-\beta_j/2}$ by (\ref{gamma-def}) and
%$\gamma_{nj,\infty} = n^{-1/2} \psi_{nj,\infty}(w_{nj},w_{nj}) /w_{nj}$ for $0<w_{nj}\le 1$.
$\tilde \gamma_{nj} = n^{-1/2} w_{nj}^{-\tau_j}$ for $j=1,\ldots,p$.
As a function of $w_{nj}$, the quantity $\tilde \gamma_{nj}$ in general differs from $\gamma_{nj}$ even up to a multiplicative constant unless $\tau_j = \beta_j /2$ as in
the case where $\mathcal G_j$ is an $L_2$-Sobolev space; see (\ref{tau-sobolev}) below.

%For Theorem~\ref{thm-Sobolev}, take each class $\mathcal G_j$, $j=1,\ldots,p$, to be a Sobolev space $\mathcal W_{r_j}^{\mym_j}$ for $r_j \ge 1$ and $\mym_j\ge 1$,
%or a bounded variation space $\mathcal V^{\mym_j}$ for $r_j=1$ and $\mym_j \ge 1$, on $[0,1]$.
%Then, by the entropy estimates in Lemmas~\ref{sobolev-entropy3} and \ref{TV-entropy} and the interpolation inequalities (Lemma~\ref{sobolev-interpolation}),
%Assumptions~\ref{ass:random-entropy} and \ref{ass:supre-entropy} are satisfied,
%with $\psi_{nj}(\delta)$ and
%$\psi_{nj,\infty}(z)$ in the polynomial forms $\delta^{1-\beta_j/2}$ and $z^{1-\beta_j/2}$, up to some multiplicative constants possibly depending on $n$, where $\beta_j=1/m_j$.
%For technical convenience, we take $\psi_{nj}(\delta) = B^*_j \delta^{1-\beta_j/2}$ for a constant $B^*_j \ge 1$,
%$\psi_{nj,\infty}(\delta) =\psi_{nj}(\delta) $ if $r_j \mym_j >1$ or
%$\psi_{nj,\infty}(\delta) =\{1+\log(n)\} \psi_{nj}(\delta)$ if $r_j=\mym_j =1$.
%Let $\gamma_{nj} = n^{-1/2} \psi_{nj}(w_{nj}) /w_{nj}$ by (\ref{gamma-def}) and $\gamma_{nj,\infty} = n^{-1/2} \psi_{nj,\infty}(w_{nj}) /w_{nj}$ for $0<w_{nj}\le 1$.
%See Section~\ref{prf:thm-Sobolev} for details.

\vspace{.1in}
\begin{thm} \label{thm-Sobolev}
Suppose that Assumptions~\ref{ass:random-entropy} and \ref{ass:supre-entropy} hold with $\psi_{nj}(\delta)$ and
$\psi_{nj,\infty}(z,\delta)$ in the forms (\ref{entropy-cond1}) and (\ref{entropy-cond2}), and condition (\ref{sup-norm-cond}) holds.
%Take $\mathcal G_j$, $\gamma_{nj}$, $\gamma_{nj,\infty}$ as described above for $j=1,\ldots,p$.
In addition, suppose that for sufficiently large $n$, $\gamma_{nj} \le w_{nj} \le 1$  and $\Gamma_n \gamma_{nj}^{1-\beta_j/2} \le 1$ for $j=1,\ldots,p$.
%Denote
%\begin{align*}
%& V_n = \max_{j=1,\ldots,p} \, w_{nj}^{\beta_j/2 -\tau_j },\quad
%W_n = \max_{j=1,\ldots,p} \,  w_{nj}^{\beta_j/2 -\tau_j + \beta_{p+1} \tau_j/2 }  ,
%\end{align*}
%where $\tau_j$ is from (\ref{sup-norm-cond}) and $\beta_{p+1} = \min_{j=1,\ldots,p} \beta_j$.
Then for any $0<\epsilon^\prime <1$ (for example, $\epsilon^\prime=\epsilon$), inequality (\ref{norm-conv-prob}) holds with $\pi = {\epsilon^\prime}^2$ and $\phi_n>0$ such that
\begin{align}
\phi_n & = O(1) \bigg\{n^{1/2}  \Gamma_n \max_j \frac{\gamma_{nj}}{\lambda_{nj}}  \max_j \frac{\tilde \gamma_{nj} w_{nj}^{\beta_{p+1} \tau_j/2}}{\lambda_{nj}} \nonumber \\
& \quad +  \max_j \frac{\tilde \gamma_{nj}}{\lambda_{nj}} \max_j \frac{ \sqrt{\log(p/\epsilon^\prime)}} {\lambda_{nj}}
+  \max_j \frac{\tilde \gamma_{nj}^2 \log(p/\epsilon^\prime)}{\lambda_{nj}^2} \bigg\}, \label{thm-Sobolev-equ}
\end{align}
where $\beta_{p+1} = \min_{j=1,\ldots,p} \beta_j$, and $O(1)$ depends only on $(C_2,C_3)$ from Lemmas~\ref{Dudley-thm} and \ref{Talagrand-thm}
and $C_4 = \max_{j=1,\ldots,p} C_{4,j}$ from condition~(\ref{sup-norm-cond}).
\end{thm}

To facilitate justification of conditions (\ref{entropy-cond1}), (\ref{entropy-cond2}), and (\ref{sup-norm-cond}),
consider the following assumption on the marginal densities of the covariates, as commonly imposed when handling random designs (e.g., Stone 1982).

\begin{ass}[Non-vanishing marginal densities] \label{density-assumption}
For $j=1,\ldots,p$, denote by $q_j(x^{(j)})$ the average marginal density function of $(X_1^{(j)},\ldots,X_n^{(j)})$, that is,
the density function associated with the probability measure $n^{-1}\sum_{i=1}^n Q_{X^{(j)}_i}$,
where $Q_{X^{(j)}_i}$ is the marginal distribution of $X^{(j)}_i$.
For some constant $0<\varrho_0 \le 1$, assume that $q_j(x^{(j)})$ is bounded from below by $\varrho_0$
simultaneously for $j=1,\ldots,p$.
\end{ass}

\begin{rem} \label{rem-Sobolev}
Conditions (\ref{entropy-cond1}), (\ref{entropy-cond2}), and (\ref{sup-norm-cond}) are satisfied under Assumption~\ref{density-assumption},
when each $\mathcal G_j$ is a Sobolev space $\mathcal W_{r_j}^{\mym_j}$ for $r_j \ge 1$ and $\mym_j\ge 1$,
or a bounded variation space $\mathcal V^{\mym_j}$ for $r_j=1$ and $\mym_j \ge 1$, on $[0,1]$.
Let $\beta_j=1/m_j$. First, (\ref{sup-norm-cond}) is implied by the interpolation inequalities for Sobolev spaces (Nirenberg 1966)
with
\begin{align}
\tau_j= (2/\beta_j+1-2/r_j)^{-1} \label{tau-sobolev}
\end{align}
and $C_{4,j}= \varrho_0 ^{-1} C_4(m_j,r_j)$ as stated in Lemma~\ref{sobolev-interpolation} of the Supplement.
Moreover, if $f_j \in {\mathcal G}^*_j(\delta)$ with $0 < \delta \le 1$, then
$ \|f_j\|_{F,j} \le 1$ and
$ \|f_j\|_Q \le \delta$, and hence $\|f_j\|_{L_{r_j}} \le\|f_j\|_\infty \le C_{4,j}$ by (\ref{sup-norm-cond}).
By rescaling the entropy estimates for Sobolev and bounded variation spaces (Lorentz et al.~1996) as in Lemmas~\ref{sobolev-entropy3} and \ref{TV-entropy} of the Supplement,
Assumptions~\ref{ass:random-entropy} and \ref{ass:supre-entropy} are satisfied
such that (\ref{entropy-cond1}) and (\ref{entropy-cond2}) hold with $B_{nj}$ independent of $n$,
and $B_{nj,\infty} =O(1) B_{nj} $ if $r_j>\beta_j$ or
$B_{nj,\infty} = O(\log^{1/2} (n)) B_{nj}$ if $r_j=\beta_j =1$.
\end{rem}

\begin{rem} \label{rem-Sobolev2}
Assumption~\ref{density-assumption} is not needed for justification of (\ref{entropy-cond1}), (\ref{entropy-cond2}), and (\ref{sup-norm-cond}),
when each class $\mathcal G_j$ is $\mathcal W_1^1$ or $\mathcal V^1$ on $[0,1]$, that is, $r_j=m_j=1$.
In this case, condition (\ref{sup-norm-cond}) directly holds with $\tau_j=1$, because $\|g_j\|_\infty \le \TV(g_j) + \|g_j\|_Q$.
Then (\ref{entropy-cond1}) and (\ref{entropy-cond2}) easily follow from the entropy estimates in Lemmas~\ref{sobolev-entropy3} and \ref{TV-entropy}.
\end{rem}

Our second result provides a sharper rate than in Theorem~\ref{thm-Sobolev}, applicable (but not limited)
to Sobolev and bounded variation spaces, provided that the following conditions hold.
For $g_j \in \mathcal G_j$, assume that $g_j(\cdot)$ can be written as $\sum_{\ell=1}^\infty \theta_{j\ell} u_{j\ell}(\cdot)$
for certain coefficients $\theta_{j\ell}$ and basis functions $u_{j\ell}(\cdot)$ on a set $\Omega$. In addition,
for certain positive constants $C_{5,1}$, $C_{5,2}$, $C_{5,3}$, $0<\tau_j<1$, and $0<w_{nj} \le 1$,
assume that for all $1\le j \le p$,
\begin{align}\label{new-cond-2}
& \sup\big\{\hbox{$\sum_{\ell=1}^k$} u_{j\ell}^2(x)/k: x\in\Omega, k \ge \ell_{j0}\big\}  \le C_{5,1},\\
\label{new-cond-1a}
& \max_{k\ge 1}\hbox{$\sum$}_{\ell_{j,k-1} < \ell \le \ell_{jk}} \theta_{j\ell}^2 \ell_{jk}^{1/\tau_j}
\le C_{5,2} (\|g_j\|_{F,j} + w_{nj}^{-1} \|g_j\|_Q)^2 ,
\\ \label{new-cond-1b}
& \hbox{$\sum$}_{\ell=1}^{\ell_{j0}} \theta_{j\ell}^2 w_{nj}^{-2}
 \le C_{5,2} (\|g_j\|_{F,j} + w_{nj}^{-1} \|g_j\|_Q)^2 ,
\end{align}
with $\ell_{jk} =\lceil (2^k/w_{nj})^{2\tau_j} \rceil$ for $k\ge 0$ and $\ell_{j,-1}=0$,
and for all $1\le j\le p$ and $k\ge 0$,
\begin{align}\label{new-cond-1c}
\sup\Big\{\|\hbox{$\sum_{\ell_{j,k-1} < \ell \le \ell_{jk}}$} \theta_{j\ell} u_{j\ell}\|_Q^2:
\hbox{$\sum_{\ell_{j,k-1} < \ell \le \ell_{jk}}$} \theta_{j\ell}^2 =1 \big\} \le C_{5,3}.
\end{align}

\vspace{.1in}
\begin{thm} \label{thm-Hilbert}
Suppose that (\ref{new-cond-2}), (\ref{new-cond-1a}), (\ref{new-cond-1b})
and (\ref{new-cond-1c}) hold as above, and %for sufficiently large $n$,
$\max_{j=1,\ldots,p} \{ e^{2/(1-\tau_j)}+ 2w_{nj}^{-\tau_j}\}\le n$.
Then
for any $0<\epsilon^\prime <1$ (for example, $\epsilon^\prime=\epsilon$), inequality (\ref{norm-conv-prob}) holds with $\pi = {\epsilon^\prime}^2$ and $\phi_n>0$ such that
\begin{align*}
& \phi_n =  O(1)
\left\{ \max_j\frac{\tilde \gamma_{nj}}{(1-\tau_j)\lambda_{nj}}  \max_j \frac{\sqrt{\log(np/\epsilon^\prime)}}{\lambda_{nj}} + \max_j \frac{\tilde \gamma_{nj}^2 \log(np/\epsilon^\prime)}{(1-\tau_j)^2\lambda_{nj}^2} \right\},
\end{align*}
where $\tilde\gamma_{nj} = n^{-1/2} w_{nj}^{-\tau_j}$%, $\mu_{p+1}= \max_{j=1,\ldots,p}\{(1-\tau_j)^{-1}\}$,
and $O(1)$ depends only on $\{C_{5,1},C_{5,2},C_{5,3}\}$.
\end{thm}

\begin{rem} \label{rem-Hilbert}
Let $\mathcal G_j$ be a Sobolev space $\mathcal W_{r_j}^{\mym_j}$
with $r_j \ge 1$, $\mym_j \ge 1$, and $(r_j\wedge 2) m_j>1$
or a bounded variation space $\mathcal V^{\mym_j}$ with $r_j=1$ and $\mym_j > 1$ (excluding $m_j=1$) on $[0,1]$.
Condition (\ref{new-cond-2}) holds for commonly used Fourier, wavelet and spline bases in $L_2$.
For any $L_2$ orthonormal bases $\{u_{j\ell},\ell\ge 1\}$,
condition (\ref{new-cond-1b}) follows from Assumption~\ref{density-assumption}
when $C_{5,2}\ge \varrho_0^{-1}$,
and condition (\ref{new-cond-1c}) is also satisfied under an additional assumption
that the average marginal density of $\{X_i^{(j)}: i=1,\ldots,n\}$ is bounded from above by $C_{5,3}$ for all $j$.
In the proof of Proposition~\ref{cor-fast0} we verify (\ref{new-cond-2}) and (\ref{new-cond-1a})
for suitable wavelet bases with $\tau_j=1/\{2m_j +1 - 2/(r_j\wedge 2)\}$, which satisfies $\tau_j<1$ because $(r_j\wedge 2) m_j>1$.
In fact, $\mathcal G_j$ is allowed to be a Besov space $\mathcal B^{m_j}_{r_j,\infty}$,
which contains $\mathcal W_{r_j}^{\mym_j}$ for $r_j\ge 1$ and $\mathcal V^{\mym_j}$
for $r_j=1$ (e.g., DeVore \& Lorentz 1993).
\end{rem}

\begin{rem} \label{rem-Hilbert-2}
The convergence rate of $\phi_n$ in Theorem \ref{thm-Hilbert} is no slower than (\ref{thm-Sobolev-equ}) in Theorem~\ref{thm-Sobolev}
if $1\le r_j\le 2$ and $(1-\tau_j)^{-1}\{\log(n)/n\}^{1/2} = O(\tilde\gamma_{nj})$,
the latter of which is valid whenever $\tau_j$ is bounded away from 1 and $\tilde \gamma_{nj} = n^{-1/2} w_{nj}^{-\tau_j/2}$ is of a slower polynomial order than $n^{-1/2}$.
However, Theorem \ref{thm-Hilbert} requires an additional side condition
(\ref{new-cond-1c}) along with
the requirement of $\tau_j<1$, which excludes for example the bounded variation space $\mathcal V^1$ on $[0,1]$.
See Equations (\ref{comparison-rate1}) and (\ref{comparison-rate2}) for implications of these rates
when used in Assumption~\ref{rate-assumption}.
\end{rem}

\subsection{Results with Sobolev and bounded variation spaces} \label{sect:specific-result}

We combine the results in Section~\ref{sect:random-design} and \ref{sect:conv-norm} (with $\epsilon^\prime=\epsilon$) to
deduce a number of concrete results on the performance of $\hat g$.
For simplicity, consider a fully homogeneous situation where each class $\mathcal G_j$ is a Sobolev space $\mathcal W_{r_0}^{\mym_0}$ for some constants $r_0 \ge 1$ and $m_0 \ge 1$
or a bounded variation space $\mathcal V^{\mym_0}$ for $r_0=1$ and $m_0\ge 1$ on $[0,1]$. Let $\beta_0=1/m_0$.
By Remark~\ref{rem-Sobolev}, if $r_0>\beta_0$, then
Assumptions~\ref{ass:random-entropy} and \ref{ass:supre-entropy} are satisfied
such that $\psi_{nj}(\delta) = B^*_0 \delta^{1-\beta_0/2}$ and
$\psi_{nj,\infty}(z,\delta)= O(1)B^*_0 z^{1-\beta_0/2}$ for $z>0$ and $0<\delta\le 1$
under Assumption~\ref{density-assumption} (non-vanishing marginal densities), where $B_0^*>0$ is a constant depending on $\varrho_0$ among others.
On the other hand, by Remark~\ref{rem-Sobolev2}, if $r_0=\beta_0=1$, then Assumptions~\ref{ass:random-entropy} and \ref{ass:supre-entropy} are satisfied such that
$\psi_{nj}(\delta) = B^*_0 \delta^{1/2}$ and
$\psi_{nj,\infty}(z,\delta)= O(\log^{1/2}(n))B^*_0 z^{1/2}$ for $z>0$ and $0<\delta\le 1$, even when Assumption~\ref{density-assumption} does not hold.
That is, $\Gamma_n$ in Theorem~\ref{thm-Sobolev} reduces to
\begin{align}
\mbox{ $\Gamma_n=O(1)$ if $r_0>\beta_0$ \;or\; $O(\log^{1/2}(n))$ if $r_0=\beta_0=1$. } \label{Gamma-n}
\end{align}

We present our results in three cases, where the underlying function $g^* = \sum_{j=1}^p g^*_j$ is assumed to satisfy (\ref{Lq-ball}) with $q=1$, $q=0$, or $0< q<1$.
As discussed in Section~\ref{sect:intro}, the parameter set (\ref{Lq-ball}) decouples sparsity and smoothness, inducing sparsity at different levels through an $L_q$ ball
in $\|\cdot\|_Q$ norm for $0 \le q \le 1$, while only enforcing smoothness through an $L_1$ ball in $\|\cdot\|_F$ norm on the components $(g^*_1,\ldots,g^*_p)$.

The first result deals with the case $q=1$ for the parameter set (\ref{Lq-ball}).

\begin{pro} \label{cor-slow}
Assume that (\ref{additive-reg}) holds and $\|g^*\|_{F,1} \le C_1 M_F$ and $\|g^*\|_{Q,1} \le C_1 M_1$ for $M_F>0$ and $M_1>0$, possibly depending on $(n,p)$.
Let $w_{nj}=1$ and $\gamma_{nj}=\gamma^*_n(1)\asymp n^{-1/2}$ by (\ref{random-rate1})--(\ref{random-rate2}).
Suppose that Assumptions~\ref{sub-gaussian-error} and \ref{density-assumption} hold, and $\log(p/\epsilon)=o(n)$.
Then for sufficiently large $n$, independently of $(M_F,M_1)$, and any $A_0> (1+\eta_0)/(1-\eta_0)$,
we have with probability at least $1-2\epsilon$,
\begin{align*}
& \|\hat g -g^*\|_n^2 + A_1 R^*_n(\hat g - g^*)
\le O(1) C_1^2 (M_F +M_1) \sqrt{\log(p/\epsilon)/n} ,
\end{align*}
where $O(1)$ depends only on $(B^*_0, A_0,\eta_0,\varrho_0)$. Moreover, we have
\begin{align*}
& \frac{1}{2} \|\hat g - g^* \|_Q^2 + A_1 R^*_n(\hat g - g^*) \\
& \le O(1) C_1^2 (M_F^2 +M_1^2) \left\{ n^{-1/2}\Gamma_n +\sqrt{\log(p/\epsilon)/n} \right\} ,
\end{align*}
with probability at least $1-2\epsilon$, where $\Gamma_n$ is from (\ref{Gamma-n}) and $O(1)$ depends only on $(B^*_0, A_0,\eta_0,\varrho_0)$ and $(C_2,$ $C_3,C_4)$ as in Theorem~\ref{thm-Sobolev}.
If $r_0=\beta_0=1$, then the results are valid even when Assumption~\ref{density-assumption} and hence $\varrho_0$ are removed.
\end{pro}

\begin{rem}[Comparison with existing results] \label{rem:cor-slow}
Proposition~\ref{cor-slow} leads to the slow rate $\{\log(p)/n\}^{1/2}$ under $L_1$-ball sparsity in $\|\cdot\|_Q$ norm,
as previously obtained for additive regression with Sobolev Hilbert spaces in Ravikumar et~al.~(2009),
except in the case where $r_0=\beta_0=1$, that is, each class $\mathcal G_j$ is $\mathcal W^1_1$ or $\mathcal V^1$.
In the latter case,
Proposition~\ref{cor-slow} shows that the convergence rate is $\{\log(np)/n\}^{1/2}$ for out-of-sample prediction,
but remains $\{\log(p)/n\}^{1/2}$ for in-sample prediction.
Previously, only the slower rate, $\{\log(np)/n\}^{1/2}$, was obtained for in-sample prediction in
additive regression with the bounded variation space $\mathcal V^1$ by Petersen et~al.~(2016).
\end{rem}

The second result deals with the case $q=0$ for the parameter set (\ref{Lq-ball}).

\begin{pro} \label{cor-fast}
Assume that (\ref{additive-reg}) holds and $\|g^*\|_{F,1} \le C_1 M_F$ and $\|g^*\|_{Q,0} \le M_0$ for $M_F>0$ and $M_0>0$,
possibly depending on $(n,p)$.
By (\ref{random-rate1})--(\ref{random-rate2}), let
\begin{align*}
& w_{nj} = w^*_n(0) = \max\left\{ {B_0^*}^{\frac{2}{2+\beta_0}} n^{ \frac{-1}{2+\beta_0} }, \, \left(\frac{\log (p/\epsilon)}{n} \right)^{1/2} \right\},\\
& \gamma_{nj} =\gamma^*_n(0) = \min\left\{{B_0^*}^{\frac{2}{2+\beta_0}} n^{ \frac{-1}{2+\beta_0} }, \, B^*_0 n^{-1/2} \left(\frac{\log (p/\epsilon)}{n}\right)^{-\frac{\beta_0}{4}}\right\}.
\end{align*}
Suppose that Assumptions~\ref{sub-gaussian-error}, \ref{compat-condition-Q}, and \ref{density-assumption} hold
with $0 < \eta_0 < (\xi_0^*-1)/(\xi^*_0+1)$ and $S = \{1 \le j \le p: \|g^*_j\|_Q > C_0^* \lambda_{nj}\}$ for some constant $C_0^* >0$,
and
\begin{align}
\Big\{ \Gamma_n {w^*_n}(0) ^{-(1-\beta_0/2) \tau_0} \gamma^*_n(0) +  {w^*_n}(0)^{-\tau_0}  \sqrt{\log(p/\epsilon)/n} \Big\}(1+M_F+M_0) = o(1), \label{cor-fast-equ}
\end{align}
where $\tau_0=1/(2/\beta_0+1-2/r_0)$.
Then for sufficiently large $n$, depending on $(M_F, M_0)$ only through the convergence rate in (\ref{cor-fast-equ}), and any $A_0> A(\xi^*_0,\eta_0)$, we have
\begin{align}
\mathcal D^*_n (\hat g, g^*) \le
O(1) C_1^2 (M_F+M_0) \left\{ n^{ \frac{-1}{2+\beta_0} }  + \sqrt{\log(p/\epsilon)/n} \right\}^{2},  \label{cor-fast-equ2}
\end{align}
with probability at least $1-2\epsilon$, where $O(1)$ depends only on $(B^*_0, A^*_0, C^*_0, \xi^*_0, \kappa^*_0, \eta_0,\varrho_0)$. % and $(C_2,C_3,C_4)$ as in Theorem~\ref{thm-Sobolev}.
If $r_0=\beta_0=1$, then the results are valid even when Assumption~\ref{density-assumption} and hence $\varrho_0$ are removed.
\end{pro}

Condition (\ref{cor-fast-equ}) is based on Theorem~\ref{thm-Sobolev} for convergence of empirical norms. By Remark~\ref{rem-Hilbert-2},
a weaker condition can be obtained using Theorem~\ref{thm-Hilbert} when $1 \le r_0 \le 2$ and $\tau_0<1$ (that is, $r_0>\beta_0$).
It is interesting to note that (\ref{cor-fast-equ}) reduces to (\ref{cor-fast-equ0}) below in the case $r_0=\beta_0=1$.

\begin{pro} \label{cor-fast0}
Proposition~\ref{cor-fast} is also valid with (\ref{cor-fast-equ}) replaced by the weaker condition
\begin{align}
\Big\{{w^*_n}(0)^{-\tau_0}  \sqrt{\log(np/\epsilon)/n} \Big\}(1+M_F+M_0) = o(1), \label{cor-fast-equ0}
\end{align}
in the case where $1 \le r_0 \le 2$, $r_0>\beta_0$, and the average marginal density of $(X_1^{(j)},\ldots,X_n^{(j)})$ is bounded from above for all $j$.
\end{pro}

\begin{rem}[Comparison with existing results] \label{rem:cor-fast1}
Propositions~\ref{cor-fast} and \ref{cor-fast0} yield the fast rate $n^{\frac{-2}{2+\beta_0}} + \log(p)/n $ under $L_0$-ball sparsity in $\|\cdot\|_Q$ norm.
Previously, the same rate was obtained for high-dimensional additive regression only with reproducing kernel Hilbert spaces (including the Sobolev space $\mathcal W^\mym_2$) by
Koltchinskii \& Yuan (2010) and Raskutti et~al.~(2012), but under more restrictive conditions.
They studied hybrid penalized estimation procedures, which involve additional constraints such that the Hilbert norms of $(g_1,\ldots,g_p)$ are bounded by {\it known} constants when minimizing
a penalized criterion.
Moreover, Koltchinskii \& Yuan (2010) assumed a constant bound on the sup-norm of possible $g^*$, whereas Raskutti et~al.~(2012) assumed the independence of the covariates $(X_i^{(1)},\ldots,X_i^{(p)})$ for each $i$.
These restrictions were relaxed in subsequent work by Suzuki \& Sugiyama (2013), but only explicitly under the assumption that the noises $\varepsilon_i$ are uniformly bounded by a constant.
Moreover, our condition (\ref{cor-fast-equ}) is much weaker than related ones in Suzuki \& Sugiyama (2013), as discussed in Remarks~\ref{rem:cor-fast2} and \ref{rem:cor-fast3} below.
See also Remark~\ref{rem-main2-cor5} for a discussion about the relationship between our results and the seemingly faster rate in Suzuki \& Sugiyama (2013).
\end{rem}

\begin{rem} \label{rem:cor-fast2}
To justify Assumptions~\ref{rate-assumption}(i)--(ii) on convergence of empirical norms, our rate condition (\ref{cor-fast-equ}) is much weaker than previous ones used.
If each class $\mathcal G_j$ is a Sobolev Hilbert space ($r_0=2$), then $\tau_0 = \beta_0/2$ and (\ref{cor-fast-equ}) becomes
\begin{align}
\left\{  n^{1/2} {w^*_n} (0)^{\beta_0^2/4} {\gamma^*_n}(0)^{2}+ \gamma^*_n (0)  \sqrt{\log(p/\epsilon)} \right\}(1+M_F+M_0) = o(1). \label{comparison-rate1}
\end{align}
%where $\gamma^*_n(0) = \min\{ \gamma_n(0), B^*_0 n^{-1/2} \nu_n^{-\beta_0/2}\}$ and $w^*_n(0) = \max\{ \gamma_n(0),\nu_n\}$ by (\ref{random-rate1})--(\ref{random-rate2}).
Moreover, by Proposition~\ref{cor-fast0}, condition (\ref{cor-fast-equ}) can be weakened to (\ref{cor-fast-equ0}), that is,
\begin{align}
\gamma^*_n (0) \sqrt{\log(np/\epsilon)} (1+M_F+M_0) = o(1), \label{comparison-rate2}
\end{align}
under an additional condition that
the average marginal density of $(X_1^{(j)},\ldots,X_n^{(j)})$ is bounded from above for all $j$. Either condition (\ref{comparison-rate1}) or (\ref{comparison-rate2}) is much weaker than
those in related analysis with reproducing kernel Hilbert spaces.
In fact, techniques based on the contraction inequality (Ledoux \& Talagrand 1991) as used in Meier et~al.~(2009) and Koltchinskii \& Yuan (2010),
lead to a rate condition such as
\begin{align}
n^{1/2} \{\gamma_n^2 (0) + \nu_n^2 \} (1+M_F+M_0) = o(1), \label{comparison-rate3}
\end{align}
where $\gamma_n(0) = {B_0^*}^{\frac{2}{2+\beta_0}} n^{ \frac{-1}{2+\beta_0} }$
and $\nu_n =  \{\log(p/\epsilon)/n\}^{1/2}$. This amounts to condition (6) assumed in Suzuki \& Sugiyama (2013), in addition to the requirement $n^{-1/2} (\log p)\le 1$.
But condition (\ref{comparison-rate3}) is even stronger than the following condition:
\begin{align}
n^{1/2}  \left\{ {\gamma_n}(0)^{2+\beta_0^2/4}+  \gamma_n (0)\nu_n  \right\} (1+M_F+M_0) = o(1),   \label{comparison-rate4}
\end{align}
because $ \Gamma_n {\gamma_n}(0)^{2+\beta_0^2/4}+  \gamma_n (0)\nu_n \ll  \gamma_n^2 (0) + \nu_n^2$
if either $\gamma_n (0) \gg \nu_n$ or $\gamma_n (0) \ll \nu_n$.
Condition (\ref{comparison-rate4}) implies (\ref{comparison-rate1}) and (\ref{comparison-rate2}), as we explain in the next remark.
\end{rem}

\begin{rem} \label{rem:cor-fast3}
Our rate condition (\ref{cor-fast-equ}) is in general weaker than
the corresponding condition with $(w^*_n(0), \gamma^*_n(0))$ replaced by $(w_n(0),\gamma_n(0))$, that is,
\begin{align}
\Big\{ \Gamma_n {\gamma_n}(0) ^{1-(1-\beta_0/2) \tau_0} +  {\gamma_n}(0)^{-\tau_0} \nu_n \Big\}(1+M_F+M_0) = o(1), \label{comparison-rate5}
\end{align}
This demonstrates the advantage of using the more careful choices $(w^*_n(0), \gamma^*_n(0))$ and also explains why (\ref{comparison-rate4}) implies (\ref{comparison-rate1}) in Remark~\ref{rem:cor-fast2}.
In fact, if $\gamma_n(0) \ge \nu_n$ then (\ref{cor-fast-equ}) and (\ref{comparison-rate5}) are identical to each other.
On the other hand, if $\gamma_n(0) < \nu_n$, then
$w^*_n(0) = \nu_n >  \gamma_n(0)$
and ${w^*_n}(0) ^{-(1-\beta_0/2) \tau_0} \gamma^*_n(0)= B^*_0 n^{-1/2} w^*_n(0)^{-(1-\beta_0/2) \tau_0 - \beta_0/2} <{\gamma_n}(0) ^{1-(1-\beta_0/2) \tau_0}$.
This also shows that if $\gamma_n(0) \ll \nu_n$, then (\ref{cor-fast-equ}) is much weaker than (\ref{comparison-rate5}).
For illustration, if $r_0=2$ and hence $\tau_0= \beta_0/2$, then
(\ref{comparison-rate5}) or equivalently (\ref{comparison-rate4}) requires at least ${\gamma_n}(0)^{-\beta_0/2} \nu_n=o(1)$, that is, $(\log p)n^{\frac{-2}{2+\beta_0}}=o(1)$,
and (\ref{comparison-rate3}) requires at least $n^{1/2} \nu_n^2 = o(1)$, that is, $\log(p) n^{-1/2} =o(1)$.
In contrast, the corresponding requirement for (\ref{cor-fast-equ}), $w^*_n(0)^{-\beta_0/2} \nu_n= o(1)$, is automatically valid as long as $\nu_n=o(1)$, that is,  $\log(p) n^{-1}=o(1)$.
\end{rem}

The following result deals with the case $0<q<1$ for the parameter set (\ref{Lq-ball}).

\begin{pro} \label{cor-medium}
Assume that (\ref{additive-reg}) holds and $\|g^*\|_{F,1} \le C_1 M_F$ and $\|g^*\|_{Q,q} \le C_1^q M_q$ for $0 < q <1$, $M_q>0$, and $M_F>0$,
possibly depending on $(n,p)$.
Let $w_{nj} = w^*_n(q)$ and $\gamma_{nj} =\gamma^*_n(q)$ by (\ref{random-rate1})--(\ref{random-rate2}).
Suppose that Assumptions~\ref{sub-gaussian-error}, \ref{compat-mono}, and \ref{density-assumption} hold
with $0 < \eta_0 < (\xi_0^*-1)/(\xi^*_0+1)$ and $S = \{1 \le j \le p: \|g^*_j\|_Q > C_0^* \lambda_{nj}\}$ for some constant $C_0^* >0$,  $\log(p/\epsilon) =o(n)$,
and
\begin{align}
\Big\{ \Gamma_n {w^*_n}(q) ^{-(1-\beta_0/2) \tau_0} \gamma^*_n(q)^{1-q} +  {w^*_n}(q)^{-\tau_0} \nu_n^{1-q} \Big\}(1+M_F+M_q) \le \eta_4, \label{cor-medium-equ}
\end{align}
for some constant $\eta_4>0$, where $\nu_n =  \{\log(p/\epsilon)/n\}^{1/2}$. Then for sufficiently large $n$, independently of $(M_F,M_q)$, and any $A_0> A(\xi^*_0,\eta_0)$, we have
\begin{align*}
\mathcal D^*_n (\hat g, g^*) \le
O(1) C_1^2 (M_F +M_q) \left\{ n^{ \frac{-1}{2+\beta_0(1-q)} }  + \sqrt{\log(p/\epsilon)/n} \right\}^{2-q},
\end{align*}
with probability at least $1-2\epsilon$, where $O(1)$ depends only on $(q, B^*_0, A^*_0, C^*_0, \xi^*_0, \kappa^*_0, \eta_0,\varrho_0,\eta_4)$ and $(C_2,C_3,C_4)$ as in Theorem~\ref{thm-Sobolev}.
If $r_0=\beta_0=1$, then the results are valid even when Assumption~\ref{density-assumption} and hence $\varrho_0$ are removed.
\end{pro}

Similarly as in Propositions~\ref{cor-fast} and \ref{cor-fast0}, condition (\ref{cor-medium-equ}) can be weakened as follows when $1 \le r_0 \le 2$ and $\tau_0<1$ (that is, $r_0>\beta_0$).
It should also be noted that (\ref{cor-medium-equ}) is equivalent to (\ref{cor-medium-equ0}) below (with different $\eta_4$ in the two equations) in the case $r_0=\beta_0=1$,
because $\gamma^*_n(q)$ with $q<1$ is of a slower polynomial order than $n^{-1/2}$ and hence $\{\log(n)/n\}^{1/2} \gamma^*_n(q)^{-1} = o(1)$.

\begin{pro} \label{cor-medium0}
Proposition~\ref{cor-medium} is also valid with (\ref{cor-fast-equ}) replaced by the weaker condition
\begin{align}
\Big\{{w^*_n}(q)^{-\tau_0} (\log(np/\epsilon)/n)^{(1-q)/2} \Big\}(1+M_F+M_0) \le \eta_4, \label{cor-medium-equ0}
\end{align}
for some constant $\eta_4>0$, in the case where $1 \le r_0 \le 2$, $r_0>\beta_0$, and the average marginal density of $(X_1^{(j)},\ldots,X_n^{(j)})$ is bounded from above for all $j$.
\end{pro}

\begin{rem} \label{rem:cor-medium}
Propositions~\ref{cor-medium} and \ref{cor-medium0} yield, under $L_q$-ball sparsity in $\|\cdot\|_Q$ norm, a convergence rate interpolating the slow and fast rates smoothly from $q=1$ to $q=0$,
similarly as in fixed designs (Section \ref{sect:fixed-design}).
However, the rate condition (\ref{cor-medium-equ}) involved does always exhibit a smooth transition to those for the slow and fast rates.
In the extreme case $q =1$, condition (\ref{cor-medium-equ}) with $q=1$ cannot be satisfied when $M_1$ is unbounded or when $M_1$ is bounded but
$\Gamma_n$ is unbounded with $r_0=\beta_0=1$.
In contrast, Proposition~\ref{cor-slow} allows for unbounded $M_1$ and the case $r_0=\beta_0=1$.
This difference is caused by the need to justify Assumption~\ref{rate-assumption}(ii) with $q \not=1$.
In the extreme case $q =0$, condition (\ref{cor-medium-equ}) with $q=0$ also differ drastically from (\ref{cor-fast-equ}) in Proposition~\ref{cor-fast}.
As seen from the proof of Corollary~\ref{main2-cor4}, this difference arises because Assumption~\ref{rate-assumption}(ii) can be justified by exploiting
the fact that $z^q\to\infty$ as $z\to\infty$ for $q>0$ (but not $q=0$).
\end{rem}

\begin{table}
\caption{Convergence rates for out-of-sample prediction under parameter set (\ref{Lq-ball}) with $(M_F,M_q)$ bounded from above} \label{rate-table}
\begin{center}
\begin{tabular}{c|c|c|c|c|c} \hline\hline
 & $r_0> \beta_0$ & \multicolumn{4}{c}{ $r_0=\beta_0=1$ } \\ \cline{3-6}
 & $0\le q \le 1$ & $q=1$ & $0<q<1$ & \multicolumn{2}{c}{$q=0$} \\ \cline{5-6}
 & & & & $\nu_n=o(\gamma_n(0))$ & otherwise  \\ \hline
scale & & & & & \\[-.1in]
adaptive & yes & yes & yes & yes & no \\ \hline
rate & $\{\gamma_n(q)+\nu_n\}^{2-q}$ & $\sqrt{\log(n)/n}+ \nu_n$ & \multicolumn{3}{c}{$\{\gamma_n(q)+\nu_n\}^{2-q}$} \\ \hline
\end{tabular}
\end{center}
{\small Note: $\gamma_n(q) \asymp n^{\frac{-1}{2+\beta_0(1-q)}}$ and $\nu_n=\{\log(p/\epsilon)/n\}^{1/2}$.
Scale-adaptiveness means the convergence rate is achieved with $(w_{nj},\gamma_{nj})$ chosen independently of $(M_F,M_q)$.}
\end{table}

For illustration, Table~\ref{rate-table} gives the convergence rates from Propositions~\ref{cor-slow}--\ref{cor-medium} in the simple situation where $(M_F,M_q)$ are bounded from above, independently of $(n,p)$.
The rate conditions (\ref{cor-fast-equ}) and (\ref{cor-medium-equ}) are easily seen to hold in all cases except that (\ref{cor-fast-equ})
is not satisfied for $q=0$ when $r_0=\beta_0=1$ but $\nu_n \not=o(\gamma_n(0))$.
In this case, we show in the following result that the convergence rate $\{\gamma_n(0)+\nu_n\}^{2}$ can still be achieved, but with
the tuning parameters $(w_{nj},\gamma_{nj})$ chosen suitably depending on the upper bound of $(M_F,M_q)$.
This is in contrast with the other cases in Table~\ref{rate-table} where the convergence rates are achieved by our penalized estimators in a scale-adaptive manner:
$(w_{nj},\gamma_{nj})=(w^*_n(q), \gamma^*_n(q))$ are chosen independently of $(M_F,M_q)$
or their upper bound.

\begin{pro} \label{cor-fast2}
Assume that $r_0=\beta_0=1$, and $M_F$ and $M_0$ are bounded from above by a constant $\overline M>0$.
Suppose that the conditions of Proposition~\ref{cor-fast} are satisfied except with (\ref{cor-fast-equ}) and Assumption~\ref{density-assumption} removed,
and Assumption~\ref{compat-mono} holds instead of Assumption~\ref{compat-condition-Q}.
Let $\hat g^\prime$ be the estimator with
$(w_{nj},\gamma_{nj})$ replaced by
$ w^\prime_{nj}= K_0 w^*_n(0)$ and $\gamma^\prime_{nj} = K_0^{-\beta_0/2} \gamma^*_n(0)$ for $K_0>0$.
Then $K_0$ can be chosen, depending on $\overline M$ but independently of $(n,p)$, such that
for sufficiently large $n$, depending on $\overline M$, and any $A_0> A(\xi^*_0,\eta_0)$,
we have
\begin{align*}
\mathcal D^*_n (\hat g^\prime, g^*) \le
O(1) C_1^2 (M_F+M_0) \left\{ n^{ \frac{-1}{2+\beta_0} }  ++\sqrt{\log(p/\epsilon)/n}\right\}^2,
\end{align*}
with probability at least $1-2\epsilon$, where $O(1)$ depends only on $(\overline M, B^*_0, A^*_0, C^*_0, \xi^*_0, \kappa^*_0, \eta_0)$ and $(C_2,C_3,C_4)$ as in Theorem~\ref{thm-Sobolev}.
\end{pro}

\section{Discussion} \label{sect:conclusion}

For additive regression with high-dimensional data, we have established new convergence results on the predictive performance of penalized estimation when
each component function can be a Sobolev space or a bounded variation space.
A number of open problems remain to be fully investigated.
First, our results provide minimax upper bounds for estimation when the component functions are restricted within an $L_1$ ball in $\|\cdot\|_{F,j}$ semi-norm and
an $L_q$ ball in $\|\cdot\|_Q$ norm. It is desirable to study whether these rates would match minimax lower bounds.
Second, while the penalized estimators have been shown under certain conditions to be adaptive to the sizes of $L_1$($\|\cdot\|_F)$ and $L_q$($\|\cdot\|_Q$) balls for fixed $q$,
we are currently investigating adaptive estimation over such balls with varying $q$ simultaneously.
Finally, it is interesting to study variable selection and inference about component functions for high-dimensional additive regression, in addition to predictive performance studied here.

\vspace{.2in}
\centerline{\bf\large References}
\begin{description}

\item Bellec, P.C. and Tsybakov, A.B. (2016) Bounds on the prediction error of penalized least squares estimators with convex penalty. {\em Festschrift in Honor of Valentin Konakov}, to appear.

\item Bickel, P., Ritov, Y., and Tsybakov, A.B. (2009) Simultaneous analysis of Lasso and Dantzig selector, {\em Annals
of Statistics}, 37, 1705--1732.

\item Bunea, F., Tsybakov, A.B., and Wegkamp, M. (2007) Sparsity oracle inequalities for the Lasso,
{\em Electronic Journal of Statistics}, 1, 169--194.

\item DeVore, R. A. and Lorentz, G.G. (1993) {\em Constructive Approximation}, Springer: New York, NY.

\item Greenshtein, E. and Ritov, Y. (2004) Persistency in high-dimensional linear predictor selection
and the virtue of over-parametrization, {\em Bernoulli} 10, 971--988.

\item Gu, C. (2002) {\em Smoothing Spline ANOVA Models}, Springer: New York, NY.

\item Hastie, T. and Tibshirani, R. (1990) {\em Generalized Additive Models}, Chapman \& Hall: New York, NY.

\item Huang, J., Horowitz, J.L., and Wei, P. (2010) Variable selection in nonparametric additive models,
{\em Annals of Statistics}, 38, 2282--2313.

\item Kim, S.-J., Koh, K., Boyd, S., and Gorinevsky, D. (2009) $\ell_1$ trend filtering, {\em SIAM Review} 51, 339--360.

\item Koltchinskii, V. and Yuan, M. (2010) Sparsity in multiple kernel learning, {\em Annals of Statistics}, 38, 3660--3695.

\item Ledoux, M. and Talagrand, M. (1991) {\em Probability in Banach Spaces: Isoperimetry and Processes}, Springer: Berlin.

\item Lin, Y. and Zhang, H.H. (2006) Component selection and smoothing in multivariate nonparametric regression,
{\em Annals of Statistics}, 34, 2272--2297

\item Lorentz, G.G., Golitschek, M.v. and Makovoz, Y. (1996) {\em Constructive Approximation: Advanced
Problems}, Springer: New York, NY.

\item Mammen, E. and van de Geer, S. (1997) Locally adaptive regression splines, {\em Annals of Statistics},
25, 387--413.

\item Meier, L., van de Geer, S., and Buhlmann, P. (2009) High-dimensional additive modeling, {\em Annals of
Statistics}, 37, 3779--3821.

\item Negahban, S.N., Ravikumar, P., Wainwright, M.J., and Yu, B. (2012) A unified framework for high-dimensional analysis of
M-estimators with decomposable regularizers, {\em Statistical Science}, 27, 538--557.

\item Nirenberg, L. (1966) An extended interpolation inequality, {\em
Annali della Scuola Normale Superiore di Pisa, Classe di Scienze}, 20, 733--737.

\item Petersen, A., Witten, D., and Simon, N. (2016) Fused Lasso additive model, {\em Journal of Computational and Graphical Statistics}, 25, 1005--1025.

\item Raskutti, G., Wainwright, M.J., and Yu, B. (2012) Minimax-optimal rates for sparse
additive models over kernel classes via convex programming, {\em Journal of Machine Learning Research} 13, 389--427.

\item Ravikumar, P., Liu, H., Lafferty, J., and Wasserman, L. (2009) SPAM: Sparse additive models, {\em Journal of
the Royal Statistical Society}, Series B, 71, 1009--1030.

\item Stone, C.J. (1982) Optimal global rates of convergence for nonparametric regression. {\em Annals of Statistics},
10, 1040--1053.

\item Stone, C.J. (1985) Additive regression and other nonparametric models. {\em Annals of Statistics}, 13, 689--705.

\item Suzuki, T. and Sugiyama, M. (2013) Fast learning rate of multiple kernel learning: Trade-off between sparsity and smoothness,
{\em Annals of Statistics}, 41, 1381--1405.

\item Tibshirani, R.J. (2014) Adaptive piecewise polynomial estimation via trend filtering, {\em Annals of
Statistics}, 42, 285--323.

\item van de Geer, S. (2000) {\em Empirical Processes in M-Estimation}, Cambridge University Press.

\item van der Vaart, A.W. and Wellner, J. (1996) {\em Weak Convergence and Empirical Processes}, Springer: New York, NY.

\item Yuan, M. and Zhou, D.-X. (2016) Minimax optimal rates of estimation in high dimensional additive models,
{\em Annals of Statistics}, 44, 2564--2593.

\item Yang, Y. and Tokdar, S.T. (2015) Minimax-optimal nonparametric regression in high dimensions, {\em Annals of Statistics},
43, 652--674.
\end{description}

\clearpage

\setcounter{page}{1}

\setcounter{section}{0}
\setcounter{equation}{0}

\renewcommand{\theequation}{S\arabic{equation}}
\renewcommand{\thesection}{S\arabic{section}}

\begin{center}
{\Large\sc Supplementary Material for ``Penalized Estimation in \\ Additive Regression with High-Dimensional Data"}

\vspace{.1in} Zhiqiang Tan \& Cun-Hui Zhang
\end{center}

\section{Proofs} \label{sect:proofs}

\subsection{Proof of Proposition~\ref{pro:BV-comput}} \label{prf:BV-comput}

Without loss of generality, assume that $k=1$ and $0 \le X^{(1)}_1 < \ldots < X^{(1)}_n \le 1$.

Consider the case $\mym=1$. For any $g = \sum_{j=1}^p g_j$ with $g_1 \in\mathcal V^1$,
define $\tilde g_1$ as a piecewise constant function:
$\tilde g_1(z) =  g_1 (X^{(1)}_1)$ for $0 \le z < X^{(1)}_2$,
$\tilde g_1(z) =  g_1(X^{(1)}_{i})$ for $X^{(1)}_{i} \le z < X^{(1)}_{i+1}$, $i=2,\ldots,n-1$, and $\tilde g_1(z) =  g_1 (X^{(1)}_n)$ for $X^{(1)}_{n} \le z \le 1$.
Let $\tilde g = \tilde g_1 +\sum_{j=2}^p  g_j$. Then $\tilde g(X_i) =g(X_i)$ for $i=1,\ldots,n$,
but $\TV(\tilde g_1) \le \TV ( g_1)$ and hence $R_n(\tilde g) \le R_n(g)$, which implies the desired result for $\mym=1$.

Consider the case $\mym=2$. For any $g = \sum_{j=1}^p g_j$ with $g_1 \in\mathcal V^2$,
define $\tilde g_1$ such that $\tilde g_1(X^{(1)}_i) = g_1(X^{(1)}_i)$, $i=1,\ldots,n$, and
$\tilde g_1(z)$ is linear in the intervals $[0, X^{(1)}_2]$, $[X^{(1)}_{i}, \, X^{(1)}_{i+1}]$, $i=2,\ldots,n-2$, and
$[X^{(1)}_{n-1},\,1]$.
Then $\TV( \tilde g_1^{(1)}) = \sum_{i=1}^{n-1} |b_{i+1} - b_i|$, where $b_i$ is the slope of $\tilde g_1$ between $[X^{(1)}_i, \,X^{(1)}_{i+1}]$.
On the other hand, by the mean-value theorem, there exists $z_i \in [X^{(1)}_i, \,X^{(1)}_{i+1}]$ such that
$g_1^{(1)}(z_i) = b_i$ for $i=1,\ldots,n-1$. Then $\TV(g_1^{(1)}) \ge\sum_{i=1}^{n-1} |b_{i+1} - b_i|$.
Let $\tilde g = \tilde g_1 +\sum_{j=2}^p  g_j$.
Then $\tilde g(X_i) =g(X_i)$ for $i=1,\ldots,n$,
but $R_n(\tilde g) \le R_n(g)$, which implies the desired result for $\mym=2$.

\subsection{Proofs of Theorem~\ref{main-thm} and corollaries} \label{sect:proof-main-thm}

We split the proof of Theorem~\ref{main-thm} and Corollary~\ref{main-cor} into five lemmas. The first one provides a probability inequality controlling the magnitude of $\langle\varepsilon, h_j\rangle_n $,
in terms of the semi-norm $\|h_j\|_{F,j}$ and the norm $\|h_j\|_n$ for all $h_j \in \mathcal G_j$ with a single $j$.

\begin{lem} \label{max-ineq2}
For fixed $j \in \{1,\ldots,p\}$, let
\begin{align*}
A_{nj} = \cup_{h_j \in \mathcal G_j}  \left\{ |\langle\varepsilon, h_j\rangle_n | /C_1 >
 \gamma_{nj,t} w_{nj} \|h_j \|_{F,j} + \gamma_{nj,t} \| h_j \|_n \right\},
\end{align*}
where $\gamma_{nj,t} = \gamma_{nj} + \sqrt{t/n}$ for $t>0$, $\gamma_{nj} = n^{-1/2} \psi_{nj}( w_{nj}) /w_{nj}$, and $w_{nj} \in (0,1]$.
Under Assumptions~\ref{sub-gaussian-error} and \ref{ass:fixed-entropy}, we have
\begin{align*}
P(A_{nj} ) \le \exp(-t).
\end{align*}
\end{lem}

\begin{prf}
In the event $A_{nj}$, we renormalize $h_j$ by letting
$ f_j = h_j / (\| h_j \|_{F,j} + \| h_j \|_n /w_{nj})$. Then $\| f_j \|_{F,j} + \| f_j \|_n /w_{nj}= 1$ and hence
$f_j \in \mathcal G_j(w_{nj})$.
By Lemma~\ref{max-ineq} with $\mathcal F_1  = \mathcal G_j(w_{nj})$ and $\delta = w_{nj}$, we have for $t>0$,
\begin{align*}
& P (A_{nj}) \le P \left\{ \sup_{f_j \in \mathcal G_j(w_{nj})  } |\langle\varepsilon, f_j\rangle_n |/C_1 > \gamma_{nj,t} w_{nj} \right\} \\
& =  P \left\{ \sup_{f_j \in \mathcal G_j(w_{nj})  } |\langle\varepsilon, f_j\rangle_n | /C_1 > n^{-1/2} \psi_{nj}( w_{nj}) + w_{nj} \sqrt{t/n}\right\}
\le \exp (-t) .
\end{align*}
\end{prf}

By Lemma~\ref{max-ineq2} and the union bound, we obtain a probability inequality controlling the magnitude of $\langle\varepsilon, h_j\rangle_n $ for $h_j \in \mathcal G_j$
simultaneously over $j=1,\ldots,p$.

\begin{lem} \label{max-ineq3}
For each $j \in \{1,\ldots,p\}$, let
\begin{align*}
A_{nj} = \cup_{h_j \in \mathcal G_j} \left\{ |\langle\varepsilon, h_j\rangle_n | >
 \lambda_{nj} w_{nj} \|h_j \|_{F,j} + \lambda_{nj} \| h_j \|_n \right\},
\end{align*}
where $\lambda_{nj}/C_1 = \gamma_{nj} + \{\log (p/ \epsilon) /n\}^{1/2} $.
Under Assumptions~\ref{sub-gaussian-error} and \ref{ass:fixed-entropy}, we have
\begin{align*}
P ( \cup_{j=1}^p A_{nj}) \le \epsilon .
\end{align*}
\end{lem}

\begin{prf}
By Lemma~\ref{max-ineq2} with $t =\log(p/\epsilon)$, we have for $j=1,\ldots,p$,
\begin{align*}
& P ( A_{nj} ) \le \exp(-t) = \frac{\epsilon}{p}.
\end{align*}
Applying the union bound yields the desired inequality.
\end{prf}

If $g^* \in \mathcal G$, then $K_n( \hat g) \le K_n( g^*)$ directly gives the basic inequality:
\begin{align}
& \frac{1}{2} \| \hat g - g^* \|_n^2  + A_0 R_n ( \hat g)
\le  \langle \varepsilon, \hat g - g^* \rangle_n +  A_0 R_n( g^*) . \label{basic-ineq}
\end{align}
By exploiting the convexity of the regularizer $R_n(\cdot)$, we provide a refinement of the basic inequality (\ref{basic-ineq}), which relates
the estimation error of $\hat g$ to that of any additive function $\bar g \in \mathcal G$ and the corresponding regularization $R_n(\bar g)$.

\begin{lem} \label{est-ineq}
The fact that $\hat g$ is a minimizer of $K_n(g)$ implies that for any function $\bar g(x) =\sum_{j=1}^p \bar g_j(x^{(j)}) \in \mathcal G$,
\begin{align}
& \frac{1}{2} \|\hat g -g^*\|_n^2 + \frac{1}{2} \| \hat g - \bar g\|^2_n + A_0 R_n ( \hat g)\nonumber \\
& \le \frac{1}{2} \|\bar g - g^*\|_n^2 + \langle\varepsilon, \hat g - \bar g \rangle_n + A_0 R_n ( \bar g). \label{basic-ineq2}
\end{align}
\end{lem}

\begin{prf}
For any $t \in (0,1]$, the fact that $K_n ( \hat g ) \le K_n( (1-t) \hat g + t  \bar g )$ implies
\begin{align*}
&\frac{t^2}{2} \|\hat g - \bar g \|_n^2 + R_n( \hat g) \le \langle Y-((1-t) \hat g + t \bar g), t(\hat g- \bar g ) \rangle_n + R_n((1-t) \hat g + t \bar g ) \\
& \le \langle y-((1-t) \hat g + t \bar g), t(\hat g- \bar g) \rangle_n + (1-t) R_n(\hat g) + t R_n(\bar g),
\end{align*}
by similar calculation leading to the basic inequality (\ref{basic-ineq}) and by the convexity of $R_n(\cdot)$: $R_n((1-t) \hat g + t \bar g ) \le (1-t) R_n(\hat g) + t R_n(\bar g)$.
Using $Y = g^*+ \varepsilon$, simple manipulation of the preceding inequality shows that for any $t \in (0,1]$,
\begin{align*}
& \langle \hat g - g^*, \hat g - \bar g\rangle_n - \frac{t}{2} \|\hat g - \bar g \|_n^2 + R_n( \hat g) \le \langle \varepsilon,\hat g- \bar g\rangle_n + R_n(\bar g),
\end{align*}
which reduces to
\begin{align*}
& \frac{1}{2} \|\hat g - g^*\|_n^2 + \frac{1-t}{2} \|\hat g - \bar g \|_n^2 + R_n( \hat g) \le \frac{1}{2} \| \bar g - g^*\|_n^2 + \langle \varepsilon,\hat g- \bar g\rangle_n + R_n(\bar g)
\end{align*}
by the fact that $2 \langle \hat g - g^*, \hat g - \bar g\rangle_n = \|\hat g - g^*\|_n^2 +  \|\hat g - \bar g \|_n^2 -\| \bar g - g^*\|_n^2$.
Letting $t \searrow0$ yields the desired inequality (\ref{basic-ineq2}).
\end{prf}

From Lemma~\ref{est-ineq}, we obtain an upper bound of the estimation error of $\hat g$ when
the magnitudes of $\langle\varepsilon, \hat g_j - \bar g_j\rangle_n $, $j=1,\ldots,p$,
are controlled by Lemma~\ref{max-ineq3}.

\begin{lem} \label{est-ineq2}
Let $A_n = \cup_{j=1}^p A_{nj}$ with $h_j = \hat g_j - \bar g_j$ in Lemma~\ref{max-ineq3}. In the event $A_n^c$,
we have for any subset $S \subset \{1,2,\ldots, p\}$,
\begin{align}
& \frac{1}{2} \|\hat g -g^*\|_n^2 + \frac{1}{2} \| \hat g - \bar g\|^2_n + (A_0-1) R_n(\hat g - \bar g) \nonumber \\
& \le \Delta_n(\bar g, S) + 2A_0 \sum_{j\in S} \lambda_{nj} \| \hat g_j- \bar g_j\|_n ,  \label{basic-ineq3}
\end{align}
where
\begin{align*}
\Delta_n(\bar g, S) = \frac{1}{2} \|\bar g - g^*\|_n^2 + 2 A_0\left( \sum_{j=1}^p  \rho_{nj}  \|\bar g_j\|_{F,j} + \sum_{j\in S^c} \lambda_{nj} \| \bar g_j \|_n\right) .
\end{align*}
\end{lem}

\begin{prf}
By the refined basic inequality (\ref{basic-ineq2}), we have in the event $A_n^c$,
\begin{align*}
& \frac{1}{2} \|\hat g -g^*\|_n^2 + \frac{1}{2} \| \hat g - \bar g\|^2_n + A_0 R_n( \hat g)  \\
& \le \frac{1}{2} \|\bar g - g^*\|_n^2  + R_n( \hat g-\bar g) + A_0 R_n( \bar g) .
\end{align*}
Applying to the preceding inequality the triangle inequalities,
\begin{align*}
\| \hat g_j\|_{F,j} & \ge \| \hat g_j - \bar g_j \|_{F,j} - \| \bar g_j \|_{F,j}, \quad j=1,\ldots,p, \\
\| \hat g_j\|_n & \ge \| \hat g_j - \bar g_j \|_n - \| \bar g_j \|_n, \quad j \in S^c, \\
\| \hat g_j\|_n & \ge \| \bar g_j \|_n - \| \hat g_j - \bar g_j \|_n , \quad j \in S,
\end{align*}
and rearranging the result leads directly to (\ref{basic-ineq3}).
\end{prf}

Taking $S=\emptyset$ in (\ref{basic-ineq3}) yields (\ref{main-thm-i}) in Corollary~\ref{main-cor}.
In general, we derive implications of (\ref{basic-ineq3}) by invoking the compatibility condition (Assumption~\ref{compat-condition}).

\begin{lem} \label{est-ineq3}
Suppose that Assumption~\ref{compat-condition} holds. If $A_0 > (\xi_0+1)/(\xi_0-1)$, then (\ref{basic-ineq3}) implies (\ref{main-thm-ii}) in Theorem~\ref{main-thm}.
\end{lem}

\begin{prf}
For the subset $S$ used in Assumption~\ref{compat-condition}, write
\begin{align*}
Z_n & = \frac{1}{2} \|\hat g -g^*\|_n^2 + \frac{1}{2} \| \hat g - \bar g\|^2_n, \\
T_{n1} & = \sum_{j=1}^p \rho_{nj} \|\hat g_j - \bar g_j \|_{F,j} + \sum_{j\in S^c} \lambda_{nj} \|\hat g_j - \bar g_j \|_n , \quad
T_{n2} = \sum_{j\in S} \lambda_{nj} \| \hat g_j- \bar g_j\|_n.
\end{align*}

Inequality (\ref{basic-ineq3}) can be expressed as
\begin{align*}
Z_n + (A_0-1) (T_{n1} + T_{n2}) \le \Delta_n(\bar g, S) + 2 A_0 T_{n2},
\end{align*}
which leads to two possible cases: either
\begin{align}
\xi_1 \{ Z_n + (A_0-1) (T_{n1} + T_{n2}) \} \le \Delta_n(\bar g, S), \label{ineq3-case1}
\end{align}
or $(1-\xi_1) \{ Z_n+ (A_0-1) (T_{n1} + T_{n2}) \} \le 2 A_0 T_{n2}$, that is,
\begin{align}
Z_n + (A_0-1) (T_{n1} + T_{n2}) \le\frac{2 A_0}{1-\xi_1}T_{n2} = (\xi_0 +1) (A_0-1) T_{n2} = \xi_2 T_{n2}, \label{ineq3-case2}
\end{align}
where $\xi_1 =1- 2A_0/\{(\xi_0 +1)(A_0-1)\} \in (0,1]$ because $A_0 > (\xi_0+1)/(\xi_0-1)$.
If (\ref{ineq3-case2}) holds, then $T_{n1} \le \xi_0 T_{n2}$, which, by Assumption~\ref{compat-condition} with $f_j = \hat g_j - \bar g_j$, implies
\begin{align}
& T_{n2} \le \kappa_0^{-1} \left( \sum_{j\in S} \lambda_{nj}^2 \right)^{1/2} \| \hat g - \bar g\|_n  . \label{ineq3-case2b}
\end{align}
Combining (\ref{ineq3-case2}) and (\ref{ineq3-case2b}) and using $\|\hat g -\bar g\|_n^2 /2\le Z_n$ yields
\begin{align}
Z_n + (A_0-1) (T_{n1} + T_{n2})\le 2 \xi_2^2 \kappa_0^{-2} \left( \sum_{j\in S} \lambda_{nj}^2 \right). \label{ineq3-case2c}
\end{align}
Therefore, inequality (\ref{basic-ineq3}), through (\ref{ineq3-case1}) and (\ref{ineq3-case2c}), implies
\begin{align*}
Z_n + (A_0-1) (T_{n1} + T_{n2})\le \xi_1^{-1} \Delta_n(\bar g, S) + 2\xi_2^2  \kappa_0^{-2} \left( \sum_{j\in S} \lambda_{nj}^2 \right) .
\end{align*}
\end{prf}

Finally, combining Lemmas~\ref{max-ineq3}, \ref{est-ineq2} and \ref{est-ineq3} completes the proof of Theorem~\ref{main-thm}.

\vspace{.1in}
\textbf{Proof of Corollary~\ref{main-cor2}.} The result follows from upper bounds of $\sum_{j\in S} \lambda_{nj}^2 $ and $\sum_{j \in S^c} \lambda_{nj} \| \bar g_j\|_n$
by the definition $S =\{1\le j \le p :\| \bar g_j \|_n > C_0 \lambda_{nj} \}$.
First, because
$\sum_{j=1}^p \lambda_{nj}^{2-q} \| \bar g_j \|_n^q \ge \sum_{j \in S} \lambda_{nj}^{2-q} (C_0^+)^q \lambda_{nj}^q $, we have
\begin{align}
\sum_{j\in S} \lambda_{nj}^2 \le (C_0^+)^{-q}  \sum_{j=1}^p \lambda_{nj}^{2-q} \| \bar g_j \|_n^q , \label{S-bound1}
\end{align}
where for $z\ge 0$, $(z^+)^q = z^q$ if $q>0$ or $=1$ if $q=0$.
Second, because $\sum_{j \in S^c} \lambda_{nj} \| \bar g_j\|_n \le \sum_{j=1}^p \lambda_{nj} (C_0\lambda_{nj})^{1-q}\| \bar g_j\|_n^q$, we have
\begin{align}
\sum_{j \in S^c} \lambda_{nj} \| \bar g_j\|_n \le  C_0^{1-q} \sum_{j=1}^p \lambda_{nj}^{2-q} \| \bar g_j\|_n^q . \label{S-bound2}
\end{align}
Inserting (\ref{S-bound1}) and (\ref{S-bound2}) into (\ref{main-thm-ii}) yields the desired inequality. \hfill $\Box$

\vspace{.1in}
\textbf{Proof of Corollary~\ref{main-cor3}.} The result follows directly from Corollary~\ref{main-cor2}, because
$\lambda_{nj}^{2-q} = C_1^{2-q} \{ \gamma_n(q) + \nu_n \}^{2-q}$
and $\rho_{nj} = C_1 \{ \gamma_n(q) + \nu_n \} \gamma_n^{1-q}(q) \le C_1 \{ \gamma_n(q) + \nu_n \}^{2-q}$, where $\nu_n=\{\log(p/\epsilon)/n\}^{1/2}$. \hfill $\Box$

\subsection{Proofs of Theorem~\ref{main2-thm} and corollaries}

Write $h_j=\hat g_j-\bar g_j$ and $h=\hat g - \bar g$ and, for the subset $S$ used in Assumption~\ref{compat-condition-Q},
\begin{align*}
Z_n & = \frac{1}{2} \|\hat g -g^*\|_n^2 + \frac{1}{2} \|h \|^2_n, \\
T^*_{n1} & = \sum_{j=1}^p \rho_{nj} \| h_j \|_{F,j} + \sum_{j\in S^c} \lambda_{nj} \| h_j \|_Q , \quad
T^*_{n2} = \sum_{j\in S} \lambda_{nj} \| h_j\|_Q.
\end{align*}
Compared with the definitions in Section~\ref{sect:proof-main-thm}, $Z_n$ is the same as before,
and $T^*_{n1}$ and $T^*_{n2}$ are similar to $T_{n1}$ and $T_{n2}$, but with $\| h_j\|_Q$ used instead of $\|h_j\|_n$.

Let
$$
\Omega_{n1} = \left\{\sup_{g \in \mathcal G} \frac{ \left| \|g\|_n^2 - \|g \|_Q^2 \right|}{ R^{*2}_n(g) } \le \phi_n \right\}.
$$
Then $P(\Omega_{n1}) \ge 1- \pi$. In the event $ \Omega_{n1}$, we have by Assumption~\ref{rate-assumption}(i),
\begin{align}
\max_{j=1,\ldots,p} \sup_{g_j \in \mathcal G_j} \frac{ \left| \|g_j\|_n - \|g_j \|_Q \right|}{ w_{nj} \|g_j\|_{F,j} + \|g_j\|_Q } \le \lambda_{n,p+1} \phi_n^{1/2} \le \eta_0. \label{norm-conv-eq1}
\end{align}
By direct calculation, (\ref{norm-conv-eq1}) implies that
if $\|g_j\|_{F,j} + \|g_j\|_n/ w_{nj} \le 1$ then $\|g_j\|_{F,j} + \|g_j\|_Q /w_{nj} \le (1-\eta_0)^{-1}$
and hence (\ref{norm-conv-eq1}) implies that
$$
H (u, \mathcal G_j(w_{nj}), \|\cdot\|_n) \le H ( (1-\eta_0)u, \mathcal G^*_j(w_{nj}), \|\cdot\|_n),
$$
and $\psi_{nj}(w_{nj})$ satisfying (\ref{random-entropy-ineq}) also satisfies (\ref{fixed-entropy-ineq}) for $\delta=w_{nj}$.
Let $\Omega_{n2} = A_n^c$ in Lemma~\ref{est-ineq2}.
Then conditionally on $X_{1:n}=(X_1,\ldots,X_n)$ for which $\Omega_{n1}$ occurs, we have $P(\Omega_{n2} |X_{1:n}) \ge 1 - \epsilon$ by Lemma~\ref{max-ineq3}.
Therefore, $P(\Omega_{n1} \cap \Omega_{n2}) \ge (1-\epsilon)(1-\pi) \ge 1-\epsilon-\pi$.

In the event $\Omega_{n2}$, recall that (\ref{basic-ineq3}) holds, that is,
\begin{align}
Z_n + (A_0-1) R_n(h) \le \Delta_n(\bar g, S) + 2A_0 \sum_{j\in S} \lambda_{nj} \|h_j \|_n . \label{basic-ineq3-re}
\end{align}
In the event $\Omega_{n1} \cap \Omega_{n2}$, simple manipulation of (\ref{basic-ineq3-re}) using (\ref{norm-conv-eq1}) shows that
\begin{align}
Z_n + A_1 R^*_n(h) \le \Delta^*_n(\bar g, S) + 2A_0 \sum_{j\in S} \lambda_{nj} \|h_j \|_Q, \label{basic-ineq3-Q}
\end{align}
where $A_1 = (A_0-1) - \eta_0(A_0+1) >0$ because $A_0 > (1+\eta_0)/(1-\eta_0)$.
In the following, we restrict to the event $\Omega_{n1} \cap \Omega_{n2}$ with probability at least $1-\epsilon-\pi$.

\vspace{.1in}
\textbf{Proof of Corollary~\ref{main2-cor}.} Taking $S = \emptyset$ in (\ref{basic-ineq3-Q}) yields (\ref{main2-thm-i1}), that is,
\begin{align*}
Z_n + A_1 R^*_n(h) \le \Delta^*_n(\bar g, g^*, \emptyset)  .
\end{align*}
As a result, $R^*_n(h) \le A_1^{-1} \Delta^*_n(\bar g, g^*, \emptyset) $
and hence $\|h\|_Q^2 \le \| h\|_n^2 + \phi_n R^{*2}_n(h) \le \| h\|_n^2 + \phi_n A_1^{-2}  \Delta^{*2}_n(\bar g, g^*, \emptyset)$.
Inequality (\ref{main2-thm-i2}) then follows from (\ref{main2-thm-i1}). \hfill $\Box$

\vspace{.1in}
\textbf{Proof of Theorem~\ref{main2-thm}.}
Inequality (\ref{basic-ineq3-Q}) can be expressed as
\begin{align*}
Z_n + A_1 (T^*_{n1} + T^*_{n2}) \le \Delta^*_n (\bar g, g^*, S) + 2 A_0 T^*_{n2},
\end{align*}
which leads to two possible cases: either
\begin{align}
\xi_1^* \{ Z_n + A_1 (T^*_{n1} + T^*_{n2})\} \le \Delta^*_n(\bar g, S), \label{ineq3-case1-Q}
\end{align}
or $(1-\xi_1^*) \{Z_n + A_1 (T^*_{n1} + T^*_{n2})\} \le 2 A_0 T^*_{n2}$, that is,
\begin{align}
Z_n+ A_1 (T^*_{n1} + T^*_{n2}) \le\frac{2 A_0}{1-\xi_1^* }T^*_{n2} = (\xi_0^* +1) A_1 T^*_{n2} = \xi_2^* T^*_{n2}, \label{ineq3-case2-Q}
\end{align}
where $ \xi_1^* = 1- 2 A_0 / \{(\xi_0^*+1)A_1\} \in (0,1]$ because
$A_0 > \{\xi^*_0+1 + \eta_0(\xi^*_0 +1)\} / \{\xi^*_0-1 - \eta_0(\xi^*_0+1)\}$.
If (\ref{ineq3-case2-Q}) holds, then $T^*_{n1} \le \xi_0^* T^*_{n2}$, which, by the theoretical compatibility condition (Assumption~\ref{compat-condition-Q}) with $f_j = \hat g_j - \bar g_j$, implies
\begin{align}
T^*_{n2} & \le \kappa_0^{*-1} \left( \sum_{j\in S} \lambda_{nj}^2 \right)^{1/2} \| h \|_Q  \label{ineq3-case2b-Q1} \\
& \le  \kappa_0^{*-1} \left( \sum_{j\in S} \lambda_{nj}^2 \right)^{1/2} \left\{ \| h \|_n + \phi_n^{1/2} (T^*_{n1}+T^*_{n2}) \right\} \label{ineq3-case2b-Q2}
\end{align}
By Assumption~\ref{rate-assumption}(ii), we have $ \phi_n\, \xi_2^{*2} \kappa_0^{*-2} (\sum_{j\in S} \lambda_{nj}^2 ) \le \eta_1^2 A_1^2$.
Combining this fact, (\ref{ineq3-case2-Q}) and (\ref{ineq3-case2b-Q2}) and simple manipulation yields
\begin{align*}
Z_n + (1-\eta_1) A_1 (T^*_{n1} + T^*_{n2}) \le \xi_2^* \kappa_0^{*-1} \left( \sum_{j\in S} \lambda_{nj}^2 \right)^{1/2} \|h\|_n,
\end{align*}
which, by the fact that $\|h\|_n^2 /2\le Z_n$, implies
\begin{align}
Z_n + (1-\eta_1) A_1 (T^*_{n1} + T^*_{n2}) \le 2 \xi_2^{*2} \kappa_0^{*-2} \left( \sum_{j\in S} \lambda_{nj}^2 \right) . \label{ineq3-case2c-Q}
\end{align}
Therefore, inequality (\ref{basic-ineq3-Q}), through (\ref{ineq3-case1-Q}) and (\ref{ineq3-case2c-Q}), implies (\ref{main2-thm-ii1}):
\begin{align*}
Z_n + (1-\eta_1) A_1 (T^*_{n1} + T^*_{n2}) \le \xi_1^{*-1} \Delta^*_n(\bar g, S) + 2\xi_2^{*2}  \kappa_0^{*-2} \left( \sum_{j\in S} \lambda_{nj}^2 \right) .
\end{align*}

To demonstrate (\ref{main2-thm-ii2}), we return to the two possible cases, (\ref{ineq3-case1-Q}) or (\ref{ineq3-case2-Q}).
On one hand, if (\ref{ineq3-case1-Q}) holds, then $A_1 R^*_n(h) = A_1( T^*_{n1} + T^*_{n2} )$ is also bounded from above by the right hand side of (\ref{ineq3-case1-Q})
and hence
\begin{align}
& \| h \|_Q^2 \le \| h\|_n^2 + \phi_n  R^{*2}_n(h)
 \le \| h\|_n^2 + \frac{\phi_n}{A_1^2} \xi_1^{*-2} \Delta^{*2}_n(\bar g, g^*, S) . \label{ineq3-case1b-Q}
\end{align}
Simple manipulation of (\ref{ineq3-case1-Q}) using (\ref{ineq3-case1b-Q}) yields
\begin{align}
& \frac{1}{2} \|\hat g -g^*\|_n^2 + \frac{1}{2} \|h\|_Q^2 + A_1 (T^*_{n1} + T^*_{n2})
\le \xi_1^{*-1} \Delta^*_n(\bar g, S)  +\frac{\phi_n}{ 2A_1^2}  \xi_1^{*-2} \Delta^{*2}_n(\bar g, g^*, S) . \label{ineq3-case1c-Q}
\end{align}
On the other hand, combining (\ref{ineq3-case2-Q}) and (\ref{ineq3-case2b-Q1}) yields
\begin{align}
Z_n + A_1 (T^*_{n1} + T^*_{n2}) \le \xi_2^* \kappa_0^{*-1} \left( \sum_{j\in S} \lambda_{nj}^2 \right)^{1/2} \|h\|_Q. \label{ineq3-case2d-Q}
\end{align}
As a result, $A_1 R^*_n(h) = A_1( T^*_{n1} + T^*_{n2} )$ is also bounded from above by the right hand side of (\ref{ineq3-case2d-Q})
and hence by Assumption~\ref{rate-assumption}(ii),
\begin{align}
& \| h \|_Q^2 \le \| h\|_n^2 + \phi_n  R^{*2}_n(h) \nonumber \\
& \le \| h\|_n^2 + \frac{\phi_n}{A_1^2} \xi_2^{*2} \kappa_0^{*-2} \left( \sum_{j\in S} \lambda_{nj}^2 \right)  \|h\|_Q^2
\le \| h\|_n^2 + \frac{1}{2} \eta_1^2 \|h\|_Q^2 . \label{ineq3-case2e-Q}
\end{align}
Simple manipulation of (\ref{ineq3-case2d-Q}) using (\ref{ineq3-case2e-Q}) yields
\begin{align*}
& \frac{1}{2} \|\hat g -g^*\|_n^2+ \frac{1-\eta_1^2}{2} \|h\|_Q^2 + A_1 (T^*_{n1} + T^*_{n2})
\le \xi_2^* \kappa_0^{*-1} \left( \sum_{j\in S} \lambda_{nj}^2 \right)^{1/2} \|h\|_Q   ,
\end{align*}
which, when squared on both sides, implies
\begin{align}
& \frac{1}{2} \|\hat g -g^*\|_n^2 + \frac{1-\eta_1^2}{2} \|h\|_Q^2 + A_1 (T^*_{n1} + T^*_{n2})
\le  \frac{2}{1-\eta_1^2} \xi_2^{*2} \kappa_0^{*-2} \left( \sum_{j\in S} \lambda_{nj}^2 \right) . \label{ineq3-case2f-Q}
\end{align}
Therefore, inequality (\ref{basic-ineq3-Q}), through (\ref{ineq3-case1c-Q}) and (\ref{ineq3-case2f-Q}), implies
\begin{align*}
& \frac{1}{2} \|\hat g -g^*\|_n^2 + \frac{1-\eta_1^2}{2} \|h\|_Q^2 + A_1 (T^*_{n1} + T^*_{n2}) \\
& \le  \xi_1^{*-1} \Delta^*_n(\bar g, S) + \frac{2}{1-\eta_1^2} \xi_2^{*2} \kappa_0^{*-2} \left( \sum_{j\in S} \lambda_{nj}^2 \right)  +
\frac{\phi_n}{ 2A_1^2} \xi_1^{*-2} \Delta^{*2}_n(\bar g, g^*, S),
\end{align*}
which yields (\ref{main2-thm-ii2}) after divided by $1-\eta_1^2$ on both sides. \hfill $\Box$

\vspace{.1in}
\textbf{Proof of Corollary~\ref{main2-cor2}.}
We use the following upper bounds, obtained from (\ref{S-bound1}) and (\ref{S-bound2}) with $S=\{1\le j\le p: \|\bar g_j\| >C_0^*\lambda_{nj}\}$,
\begin{align}
\sum_{j\in S} \lambda_{nj}^2 \le (C_0^{*+})^{-q}  \sum_{j=1}^p \lambda_{nj}^{2-q} \| \bar g_j \|_Q^q , \label{S-bound1-Q}
\end{align}
and
\begin{align}
\sum_{j \in S^c} \lambda_{nj} \| \bar g_j\|_Q \le  C_0^{* 1-q} \sum_{j=1}^p \lambda_{nj}^{2-q} \| \bar g_j\|_Q^q . \label{S-bound2-Q}
\end{align}
Equations (\ref{rate-assumption-ii}) and (\ref{main2-cor2-equ}) together imply $\phi_n \Delta^{*}_n(\bar g,,S) = O(1) + \phi_n \|\bar g-g^*\|_n^2/2$.
Inserting this into (\ref{main2-thm-ii2}) and applying (\ref{S-bound1-Q}) and (\ref{S-bound2-Q}) yields the high-probability result about $\mathcal D^*_n(\hat g, \bar g)$.
The in-probability result follows by combining the facts that $\epsilon \to 0$, $\|\bar g-g^*\|_n^2 = O_p(1) \|\bar g-g^*\|_Q^2$ by the Markov inequality,
and $\|\hat g - g^*\|_Q^2 \le 2(\|\hat g-\bar g\|_Q^2 + \|\bar g-g^*\|_Q^2)$ by the triangle inequality.
\hfill $\Box$

\vspace{.1in}
\textbf{Proof of Corollary~\ref{main2-cor3}.}
First, we show
\begin{align}
w^*_n(q) \le \{\gamma^*_n(q) + \nu_n\}^{1-q} . \label{w-gam}
\end{align}
In fact, if $\gamma_n(q) \ge \nu_n$, then $\gamma^*_n(q) = \gamma_n(q)$
and $w^*_n(q) = \gamma_n(q)^{1-q} \le \{\gamma^*_n(q) + \nu_n\}^{1-q}$.
If $\gamma_n(q) < \nu_n$, then $w^*_n(q) = \nu_n^{1-q} \le \{\gamma^*_n(q) + \nu_n\}^{1-q}$.
By (\ref{S-bound1-Q}), (\ref{S-bound2-Q}), and (\ref{w-gam}), inequality (\ref{main2-cor3-equ}) implies that for any constants $0<\eta_1<1$ and $\eta_2>0$,
(\ref{rate-assumption-ii}) and (\ref{main2-cor2-equ}) are satisfied for sufficiently large $n$.
The desired result follows from Corollary~\ref{main2-cor2} with $\bar g=g^*$, because
$\lambda_{nj}^{2-q} = C_1^{2-q} \{ \gamma^*_n(q) + \nu_n \}^{2-q} \le C_1^{2-q} \{ \gamma_n(q) + \nu_n \}^{2-q}$
and, by (\ref{w-gam}), $\rho_{nj} = C_1  w^*_n(q) \{ \gamma^*_n(q) + \nu_n \}\le C_1 \{ \gamma_n(q) + \nu_n \}^{2-q}$. \hfill $\Box$

\vspace{.1in}
\textbf{Proof of Corollary~\ref{main2-cor4}.}
For a constant $0 < \eta_1 <1$, we choose and fix ${C_0^*}^\prime \ge C_0^*$ sufficiently large, depending on $q>0$, such that
\begin{align*}
(\xi^*_0 +1)^2 \kappa_0^{*-2} \eta_3 \le ({C_0^*}^\prime)^q \eta_1^2 .
\end{align*}
Let $S^\prime = \{1\le j\le p: \|g^*_j \|_Q > {C_0^*}^\prime \lambda_{nj}\}$.
Then (\ref{rate-assumption-ii}) is satisfied with $S$ replaced by $S^\prime$, due to (\ref{main2-cor4-equ}), (\ref{S-bound1-Q}), and the definition $\lambda_{nj}=C_1\{\gamma^*_n(q)+\nu_n\}$.
Similarly, (\ref{main2-cor2-equ}) is satisfied with $S$ replaced by $S^\prime$ for $\eta_2 = \overline M_q + ({C_0^*}^\prime)^{1-q} \overline M_q$, by (\ref{S-bound2-Q}) and simple manipulation.
By Remark~\ref{mono-S}, Assumption~\ref{compat-mono} implies Assumption~\ref{compat-condition-Q} and remains valid when
$S$ is replaced by $S^\prime \subset S$.
The desired result follows from Corollary~\ref{main2-cor2} with $\bar g=g^*$.  \hfill $\Box$

\vspace{.1in}
\textbf{Proof of Corollary~\ref{main2-cor5}.}
The proof is similar to that of Corollary~\ref{main2-cor3}. First, we show
\begin{align}
w^\dag_n(q) M_F \le \{\gamma^\dag_n(q) + \nu_n\}^{1-q} M_q . \label{w-gam2}
\end{align}
In fact, if $\gamma^\prime_n(q) \ge \nu_n$, then $\gamma^\dag_n(q) = \gamma^\prime_n(q)$
and $w^\dag_n(q) M_F = {\gamma^\prime_n(q)}^{1-q} M_q \le \{\gamma^\dag_n(q) + \nu_n\}^{1-q} M_q$.
If $\gamma_n(q) < \nu_n$, then $w^\dag_n(q) M_F = \nu_n^{1-q} M_q \le \{\gamma^\dag_n(q) + \nu_n\}^{1-q} M_q$.
Then (\ref{main2-cor5-equ}) implies that for any constants $0<\eta_1<1$ and $\eta_2>0$,
(\ref{rate-assumption-ii}) and (\ref{main2-cor2-equ}) are satisfied for sufficiently large $n$.
The desired result follows from Corollary~\ref{main2-cor2} with $\bar g=g^*$, because
$\lambda_{nj}^{2-q} M_q = C_1^{2-q} \{ \gamma^\dag_n(q) + \nu_n \}^{2-q} M_q \le C_1^{2-q} \{ \gamma^\prime_n(q) + \nu_n \}^{2-q} M_q$
and, by (\ref{w-gam2}), $\rho_{nj} M_F = C_1  w^\dag_n(q) \{ \gamma^\dag_n(q) + \nu_n \} M_F \le C_1 \{ \gamma^\prime_n(q) + \nu_n \}^{2-q} M_q$. \hfill $\Box$

\subsection{Proof of Theorem~\ref{thm-Sobolev}} \label{prf:thm-Sobolev}

We split the proof into three lemmas.
First, we provide maximal inequalities on convergence of empirical inner products in functional classes
with polynomial entropies.

\begin{lem} \label{empirical-inner-product}
Let $\mathcal F_1$ and $\mathcal F_2$ be two functional classes such that
\begin{align*}
\sup_{f_j \in\mathcal F_j} \|f_j\|_Q \le \delta_j  , \quad \sup_{f_j \in\mathcal F_j}  \|f_j \|_\infty \le b_j , \quad j=1,2. %\label{norm-bounds}
\end{align*}
Suppose that for some $0<\beta_j<2$ and $B_{nj,\infty} >0$, condition (\ref{entropy-bound}) holds with
\begin{align}
\psi_{n,\infty}(z,\mathcal F_j) = B_{nj,\infty}\, z^{1-\beta_j/2}, \quad j=1,2. \label{entropy-form}
\end{align}
Then we have
\begin{align}
& E \left\{ \sup_{f_1 \in \mathcal F_1, f_2 \in \mathcal F_2} \left| \langle f_1, f_2 \rangle_n - \langle f_1, f_2 \rangle_Q \right| /C_2 \right\} \nonumber \\
& \le 2\left\{ \delta_1 + \frac{2C_2 \psi_{n,\infty}(b _1, \mathcal F_1) }{\sqrt n} \right\}^{1-\beta_1/2}
\left\{ \delta_2 + \frac{2C_2\psi_{n,\infty}(b _2, \mathcal F_2) }{\sqrt n} \right\}^{\beta_1/2}  \frac{\psi_{n,\infty} (b _2, \mathcal F_1)}{\sqrt n} \nonumber \\
& \quad + 2\left\{ \delta_2 + \frac{2C_2\psi_{n,\infty}(b _2, \mathcal F_2) }{\sqrt n} \right\}^{1-\beta_2/2}
\left\{ \delta_1 + \frac{2C_2\psi_{n,\infty}(b _1, \mathcal F_1) }{\sqrt n} \right\}^{\beta_2/2}  \frac{\psi_{n,\infty} (b _1, \mathcal F_2)}{\sqrt n}  . \label{inner-prod-exp}
\end{align}
Moreover, we have for any $t>0$,
\begin{align}
& \sup_{f_1 \in \mathcal F_1, f_2 \in \mathcal F_2} \left| \langle f_1, f_2 \rangle_n - \langle f_1, f_2 \rangle_Q \right| / C_3 \nonumber \\
& \le E \left\{ \sup_{f_1 \in \mathcal F_1, f_2 \in \mathcal F_2} \left| \langle f_1, f_2 \rangle_n - \langle f_1, f_2 \rangle_Q \right| \right\}
+ \delta_1 b_2 \sqrt{\frac{t}{n}} + b _1 b _2 \frac{t}{n} , \label{inner-prod-prob}
\end{align}
with probability at least $1-\me^{-t}$.
\end{lem}

\begin{prf}
%Let $\hat \delta_1 = \sup_{f_1 \in \mathcal F_1} \|f_1\|_n$ and  $\hat \delta_2 = \sup_{f_2 \in \mathcal F_2} \|f_2\|_n$.
For any function $f_1, f'_1 \in \mathcal F_1$ and $f_2, f'_2 \in \mathcal F_2$, we have by triangle inequalities,
\begin{align*}
\| f_1 f_2 - f'_1 f'_2 \|_n \le \hat \delta_2 \| f_1 - f'_1 \|_{n,\infty} + \hat \delta_1 \| f_2 - f'_2 \|_{n,\infty}.
\end{align*}
As a result, we have for $u>0$,
\begin{align}
H(u, \mathcal F_1 \times \mathcal F_2, \|\cdot\|_n)
\le H\{ u/(2 \hat \delta_2), \mathcal F_1, \|\cdot\|_{n,\infty} \} +  H\{ u/(2 \hat \delta_1), \mathcal F_2 , \|\cdot\|_{n,\infty} \}, \label{product-entropy}
\end{align}
where $ \mathcal F_1 \times \mathcal F_2 = \{f_1 f_2: f_1 \in \mathcal F_1, f_2 \in \mathcal F_2\}$.

By symmetrization inequality (van der Vaart \& Wellner 1996),
\begin{align*}
& E \left\{ \sup_{f_1 \in \mathcal F_1, f_2 \in \mathcal F_2} \left| \langle f_1, f_2 \rangle_n - \langle f_1, f_2 \rangle_Q \right| \right\}  \le
2 E \left\{ \sup_{f_1\in\mathcal F_1, f_2\in\mathcal F_2} | \langle \sigma, f_1f_2 \rangle_n | \right\}.
\end{align*}
Let $\hat\delta_{12} = \sup_{f_1 \in \mathcal F_1, f_2 \in \mathcal F_2} \|f_1 f_2\|_n \le \min( \hat \delta_1 b _2, \hat \delta_2 b _1)$.
By Dudley's inequality (Lemma~\ref{Dudley-thm}) conditionally on $X_{1:n}=(X_1,\ldots,X_n)$, we have
\begin{align*}
E \left\{ \sup_{f_1\in\mathcal F_1, f_2\in\mathcal F_2} | \langle \sigma, f_1f_2 \rangle_n | \Big | X_{1:n} \right\} /C_2
\le E\left\{ \int_0^{\hat\delta_{12}} H^{1/2} (u, \mathcal F_1 \times \mathcal F_2, \|\cdot\|_n) \,\dif u  \Big| X_{1:n} \right\}.
\end{align*}
Taking expectations over $X_{1:n}$, we have by (\ref{product-entropy}), (\ref{entropy-bound}), and definition of $H^*()$,
\begin{align}
& E \left\{ \sup_{f_1 \in \mathcal F_1, f_2 \in \mathcal F_2} \left| \langle f_1, f_2 \rangle_n - \langle f_1, f_2 \rangle_Q \right| /C_2 \right\} \nonumber  \\
& \le E \left[  \int_0^{\hat\delta_1 b _2} {H^*}^{1/2} \{u/(2 \hat\delta_2), \mathcal F_1, \|\cdot\|_{n,\infty}\} \,\dif u +
 \int_0^{\hat\delta_2 b _1} {H^*}^{1/2} \{u/(2 \hat\delta_1), \mathcal F_2, \|\cdot\|_{n,\infty}\} \,\dif u \right] \nonumber \\
& \le E \left[ \hat \delta_2 \psi_{n,\infty} (\hat\delta_1 b _2 /\hat \delta_2, \mathcal F_1 ) + \hat \delta_1 \psi_{n,\infty} (\hat\delta_2 b _1 /\hat \delta_1, \mathcal F_2) \right] . \label{sym-ineq}
\end{align}
By (\ref{entropy-form}) and the H\"{o}lder inequality, we have
\begin{align*}
& E \left\{ \hat \delta_2 \psi_{n,\infty} (\hat\delta_1 b _2 /\hat \delta_2, \mathcal F_1)\right\}
\le B_{n1,\infty} b _2^{1-\beta_1/2} E\left( \hat\delta_2^{\beta_1/2} \hat\delta_1 ^{1-\beta_1/2} \right) \\
& \le B_{n1,\infty} b _2^{1-\beta_1/2} E^{\beta_1/2} ( \hat\delta_2 ) E^{1-\beta_1/2} ( \hat\delta_1 )
\le  B_{n1,\infty} b _2^{1-\beta_1/2} E^{\beta_1/4} ( \hat\delta_2^2 ) E^{(2-\beta_1)/4} ( \hat\delta_1^2 ),
\end{align*}
and similarly
\begin{align*}
E \left\{\hat \delta_1 \psi_{n,\infty} (\hat\delta_2 b _1 / \hat \delta_1, \mathcal F_2 ) \right\}
\le  B_{n2,\infty} b _1^{1-\beta_2/2} E^{\beta_2/4} ( \hat\delta_1^2 ) E^{(2-\beta_2)/4} ( \hat\delta_2^2 ).
\end{align*}
Then inequality (\ref{inner-prod-exp}) follows from (\ref{sym-ineq}) and Lemma~\ref{empirical-norm}.
Moreover, inequality (\ref{inner-prod-prob}) follows from Talagrand's inequality (Lemma~\ref{Talagrand-thm}) because
$ \|f_1 f_2 \|_Q \le \delta_1 b _2$ and $\|f_1 f_2 \|_\infty \le b _1 b _2$ for $f_1 \in \mathcal F_1$ and $f_2 \in \mathcal F_2$.
\end{prf}

By application of Lemma~\ref{empirical-inner-product}, we obtain the following result on
uniform convergence of empirical inner products under conditions (\ref{entropy-cond1}), (\ref{entropy-cond2}), and (\ref{sup-norm-cond}).

\begin{lem} \label{empirical-inner-prod2}
Suppose the conditions of Theorem~\ref{thm-Sobolev} are satisfied for $j=1,2$ and $p=2$.
Let $\mathcal F_j = \mathcal G^*_j(w_{nj})$ for $j=1,2$.
Then we have
\begin{align*}
& E \left\{ \sup_{f_1 \in \mathcal F_1, f_2 \in \mathcal F_2} \left| \langle f_1, f_2 \rangle_n - \langle f_1, f_2 \rangle_Q \right| /C_2 \right\}  \\
& \le 2(1+2C_2  C_4) C_4 n^{1/2} \Gamma_n w_{n1} w_{n2}
\left( \gamma_{n1} \tilde \gamma_{n2}  w_{n2}^{ \beta_1 \tau_2/2 } +  \gamma_{n2} \tilde \gamma_{n1}  w_{n1}^{ \beta_2 \tau_1/2}  \right).
\end{align*}
where $0<\tau_j \le (2/\beta_j-1)^{-1}$ and $C_4 =\max_{j=1,2} C_{4,j}$ from condition (\ref{sup-norm-cond}), and $\tilde\gamma_{nj} = n^{-1/2} w_{nj}^{-\tau_j}$.
Moreover, we have for any $t>0$,
\begin{align*}
& \sup_{f_1 \in \mathcal F_1, f_2 \in \mathcal F_2} \left| \langle f_1, f_2 \rangle_n - \langle f_1, f_2 \rangle_Q \right| /C_3 \\
& \le E \left\{ \sup_{f_1 \in \mathcal F_1, f_2 \in \mathcal F_2} \left| \langle f_1, f_2 \rangle_n - \langle f_1, f_2 \rangle_Q \right| \right\}
+ w_{n1}  w_{n2} \left(C_4 \, t^{1/2}\tilde \gamma_{n2}  + C_4^2 \, t \tilde \gamma_{n1}\tilde \gamma_{n2}  \right),
\end{align*}
with probability at least $1-\me^{-t}$.
\end{lem}

\begin{prf}
For $f_j \in \mathcal F_j$ with $w_{nj}\le 1$, we have $ \|f_j\|_{F,j} \le 1$ and
$ \|f_j\|_Q \le w_{nj}$, and hence $\|f_j\|_\infty \le C_4 w_{nj}^{1-\tau_j}  $ by (\ref{sup-norm-cond}).
Let $\psi_{n,\infty}(\cdot,\mathcal F_j) =\psi_{nj,\infty}(\cdot, w_{nj})$ from (\ref{entropy-cond2}), that is, in the form (\ref{entropy-form}) such that (\ref{entropy-bound}) is satisfied.
We apply Lemma~\ref{empirical-inner-prod2} with $\delta_j= w_{nj}$ and $b_j = C_4 w_{nj}^{1-\tau_j}$.
By simple manipulation, we have
\begin{align*}
& n^{-1/2} \psi_{n,\infty}(b _j, \mathcal F_j) = n^{-1/2} \psi_{nj,\infty}( C_4 w_{nj}^{1-\tau_j}, w_{nj}) \\
& \le C_4 B_{nj,\infty} n^{-1/2} w_{nj}^{-\beta_j/2} w_{nj}^{1-(1-\beta_j/2)\tau_j } \le C_4 \Gamma_n \gamma_{nj} w_{nj}^{1-\beta_j/2} \le C_4 w_{nj},
\end{align*}
where $C_4\ge 1$ is used in the second step, $B_{nj,\infty}\le \Gamma_n B_{nj}$ and $(1-\beta_j/2)\tau_j \le \beta_j/2$ in the third step,
and $\gamma_{nj} \le w_{nj}$ and $\Gamma_n \gamma_{nj} w_{nj}^{-\beta_j/2} \le \Gamma_n \gamma_{nj}^{1-\beta_j/2} \le 1$
in the fourth step.
Therefore, inequality (\ref{inner-prod-exp}) yields
\begin{align*}
& E \left\{ \sup_{f_1 \in \mathcal F_1, f_2 \in \mathcal F_2} \left| \langle f_1, f_2 \rangle_n - \langle f_1, f_2 \rangle_Q \right|/C_2 \right\} \nonumber \\
& \le 2(1+2C_2  C_4) n^{-1/2} w_{n1}^{1-\beta_1/2} w_{n2}^{\beta_1/2} \psi_{n1,\infty}(C_4 w_{n2}^{1-\tau_2}, w_{n1})\\
& \quad + 2(1+2C_2  C_4) n^{-1/2} w_{n2}^{1-\beta_2/2}  w_{n1}^{\beta_2/2} \psi_{n2,\infty}(C_4 w_{n1}^{1-\tau_1}, w_{n2}) \\
& \le 2(1+2C_2  C_4) C_4 n^{-1/2} w_{n1}^{1-\beta_1/2} B_{n1,\infty} w_{n2} w_{n2}^{-\tau_2 + \beta_1 \tau_2/2 } \\
& \quad + 2(1+2 C_2 C_4)C_4 n^{-1/2} w_{n2}^{1-\beta_2/2} B_{n2,\infty} w_{n1} w_{n1}^{-\tau_1 + \beta_2 \tau_1/2 },
\end{align*}
which leads to the first desired inequality because $B_{nj,\infty}\le \Gamma_n B_{nj}$.
Moreover, simple manipulation gives
\begin{align*}
&\delta_1 b_2 \sqrt{\frac{t}{n}} = C_4 w_{n1} w_{n2}^{1-\tau_2 } \sqrt{\frac{t}{n}}
= C_4 \, t^{1/2} w_{n1} w_{n2} \tilde \gamma_{n2} ,\\
&b_1 b_2 \frac{t}{n} =  C_4^2 w_{n1}^{1-\tau_1} w_{n2}^{1-\tau_2} \frac{t}{n}
= C_4^2\, t w_{n1} \tilde \gamma_{n1} w_{n2} \tilde \gamma_{n2} .
\end{align*}
The second desired inequality follows from (\ref{inner-prod-prob}).
\end{prf}

The following result concludes the proof of Theorem~\ref{thm-Sobolev}.

\begin{lem}
In the setting of Theorem~\ref{thm-Sobolev}, let
\begin{align*}
& \phi_n = 4 C_2 C_3 (1+ 2C_2  C_4)C_4 n^{1/2} \Gamma_n \max_j \frac{\gamma_{nj}}{\lambda_{nj}} \max_j \frac{\tilde\gamma_{nj} w_{nj}^{\beta_{p+1} \tau_j/2}}{\lambda_{nj}} \\
& \quad + \sqrt 2 C_3  C_4 \max_j \frac{\tilde \gamma_{nj}}{\lambda_{nj}} \max_j \frac{ \sqrt{\log(p/\epsilon^\prime)}} {\lambda_{nj}}
+ 2 C_3  C_4^2 \max_j \frac{\tilde \gamma_{nj}^2 \log(p/\epsilon^\prime)}{\lambda_{nj}^2} ,
\end{align*}
where $\tilde\gamma_{nj} = n^{-1/2} w_{nj}^{-\tau_j}$ and $\beta_{p+1}=\min_{j=1,\ldots,p} \beta_j$. Then
\begin{align*}
P\left\{ \sup_{g\in \mathcal G} \frac{ \left| \|g\|_n^2 - \|g \|_Q^2 \right|}{ R^{*2}_n(g) } > \phi_n \right\} \le {\epsilon^\prime}^2 .
\end{align*}
\end{lem}

\begin{prf}
For $j=1,\ldots,p$, let $r^*_{nj} ( g_j) = \|g_j\|_{F,j} + \|g_j\|_Q /w_{nj}$ and $f_j = g_j /r^*_j ( g_j) $.
Then $\|f_j\|_{F,j} + \|f_j\|_Q /w_{nj}=1$ and hence $f_j \in {\mathcal G}^*_j( w_{nj} )$.
By the decomposition $ \|g\|_n^2 = \sum_{j,k} \langle g_j, g_k \rangle_n $, $ \|g\|_Q^2 = \sum_{j,k} \langle g_j, g_k \rangle_Q $, and
the triangle inequality, we have
\begin{align*}
& \left| \|g\|_n^2 - \|g \|_Q^2 \right| \le \sum_{j,k} \left| \langle g_j, g_k \rangle_n - \langle g_j, g_k \rangle_Q \right| \nonumber \\
& = \sum_{j,k} r^*_{nj}(g_j) r^*_{nk}(g_k) \left| \langle f_j, f_k \rangle_n - \langle f_j, f_k \rangle_Q \right| .
\end{align*}
Because $R^{*2}(g)= \sum_{j,k} r^*_{nj}(g_j) r^*_{nk}(g_k) w_{nj} \lambda_{nj} w_{nk} \lambda_{nk}  $, we have
\begin{align*}
& \left\{ \sup_{g=\sum_{j=1}^p g_j} \frac{ \left| \|g\|_n^2 - \|g \|_Q^2 \right|}{ R^{*2}(g) } > \phi_n \right\}
= \bigcup_{g=\sum_{j=1}^p g_j} \Big\{ \left| \|g\|_n^2 - \|g \|_Q^2 \right| > \phi_n R^{*2}(g) \Big\} \\
& \subset \bigcup_{j,k} \left\{ \sup_{f_j \in {\mathcal G}^*(w_{nj}) , f_k\in {\mathcal G}^*(w_{nk})}
\left| \langle f_j, f_k \rangle_n - \langle f_j, f_k \rangle_Q \right|
> \phi_n \, w_{nj} \lambda_{nj} w_{nk} \lambda_{nk} \right\}
\end{align*}
By Lemma~\ref{empirical-inner-prod2} with $\mathcal F_1 = {\mathcal G}^*_j ( w_{nj})$,  $\mathcal F_2 = {\mathcal G}^*_k( w_{nk})$,
and $t=\log (p^2 /{\epsilon^\prime}^2)$, we have with probability no greater than ${\epsilon^\prime}^2 /p^2$,
\begin{align*}
& \sup_{f_j \in {\mathcal G}^*(w_{nj}) , f_k\in {\mathcal G}^*(w_{nk}) }
\left| \langle f_j, f_k \rangle_n - \langle f_j, f_k \rangle_Q \right| /C_3 \\
& > 4 C_2 (1+2 C_2 C_4)C_4 n^{1/2} \Gamma_n W_n w_{nj} \gamma_{nj} w_{nk} \gamma_{nk} \\
& \quad +  C_4 n^{1/2} V_n w_{nj} w_{nk} \gamma_{nk} \sqrt{\log(p^2 /{\epsilon^\prime}^2)/n} + C_4^2 V_n^2 \log(p^2 /{\epsilon^\prime}^2) w_{nj} \gamma_{nj} w_{nk} \gamma_{nk} .
\end{align*}
Therefore, we have by the definition of $\phi_n$,
\begin{align*}
& P\left( \sup_{f_j \in {\mathcal G}^*(w_{nj}) , f_k\in {\mathcal G}^*(w_{nk})}
\left| \langle f_j, f_k \rangle_n - \langle f_j, f_k \rangle_Q \right|
> \phi_n \, w_{nj} \lambda_{nj} w_{nk} \lambda_{nk} \right)
\le \frac{{\epsilon^\prime}^2}{p^2} .
\end{align*}
The desired result follows from the union bound.
\end{prf}

\subsection{Proofs of Propositions \ref{cor-slow}, \ref{cor-fast}, \ref{cor-medium}, and \ref{cor-fast2}}

Denote $w_{nj} =w_{n,p+1}$ and $\gamma_{nj} =\gamma_{n,p+1}$ for $j=1,\ldots,p$. By direct calculation, (\ref{thm-Sobolev-equ}) implies
that for any $0\le q \le 1$,
\begin{align}
 \phi_n(\gamma_{n,p+1}+\nu_n)^{2-q} \le & \,O(1) \Big\{ n^{1/2} \Gamma_{n} W_n \gamma^{2-q}_{n,p+1} + n^{1/2} V_n \min\left( \gamma_{n,p+1} \nu_n^{1-q},\, \gamma_{n,p+1}^{1-q}\nu_n\right) \nonumber \\
& \quad + n V_n^2 \min\left(\gamma^2_{n,p+1} \nu_n^{2-q}, \gamma^{2-q}_{n,p+1} \nu_n^2\right) \Big\} , \label{phi-lam-rate}
\end{align}
where
\begin{align}
V_n = w_{n,p+1}^{\beta_0/2-\tau_0} , \quad W_n = w_{n,p+1}^{\beta_0/2-\tau_0+\beta_0\tau_0/2} \label{Vn-equ}.
\end{align}

We verify that the technical conditions hold as needed for Theorem~\ref{thm-Sobolev}, with $w_{nj}=w_n^*(q)$ and $\gamma_{nj}=\gamma_n^*(q)$ for $0\le q \le 1$.
First, we verify $\gamma_{nj} \le w_{nj}$ for sufficiently large $n$. It suffices to show that $\gamma^*_n(q) \le w^*_n(q)$ whenever $\gamma_n(q) \le 1$ and $\nu_n\le 1$.
In fact, if $\gamma_n(q) \ge \nu_n$, then $w^*_n(q) = \gamma_n(q)^{1-q}$ and $\gamma^*_n(q) =\gamma_n(q) \le \gamma_n(q)^{1-q}$ provided $\gamma_n(q) \le 1$.
If $\gamma_n(q) < \nu_n$, then $w^*_n(q) = \nu_n^{1-q}$ and $\gamma^*_n(q) = B^*_0 n^{-1/2} \nu_n^{-(1-q)\beta_0/2} \le \nu_n \le \nu_n^{1-q}$ provided $\nu_n\le 1$.
Moreover, we have $\Gamma_n \gamma^*_n(q)^{1-\beta_0/2}\le 1$ for sufficiently large $n$, because
$\Gamma_n$ is no greater than $O(\log^{1/2}( n))$ and
$ \gamma^*_n(q)^{1-\beta_0/2} \le \gamma_n(q)^{1-\beta_0/2}$ decreases polynomially in $n^{-1}$ for $0<\beta_0<2$.

\vspace{.1in}
\textbf{Proof of Proposition~\ref{cor-slow}.} For $w_{nj}=1$ and $\gamma_{nj}=\gamma^*_n(1)\asymp n^{-1/2}$, inequality (\ref{phi-lam-rate}) with $q=0$
and $\nu_n = o(1)$ gives
\begin{align*}
\phi_n \{\gamma^*_n(1)+\nu_n\}^2 &\le O(1) \left\{ n^{1/2} \Gamma_n \gamma^{*2}_n(1) + n^{1/2} \gamma^*_n(1) \nu_n + n \gamma^{*2}_n(1) \nu_n^2 \right\} \\
& =  O(1) \left( n^{-1/2} \Gamma_n + \nu_n \right),
\end{align*}
Assumption~\ref{rate-assumption}(i) holds because $\Gamma_n$ is no greater than $O(\log^{1/2}( n))$.
Inserting the above inequality into (\ref{main2-thm-i2}) in Corollary~\ref{main2-cor} yields the out-of-sample prediction result.
The in-sample prediction result follows directly from Corollary~\ref{main2-cor}.
\hfill $\Box$

\vspace{.1in}
\textbf{Proof of Proposition~\ref{cor-fast}.}
For $\gamma_{nj}=\gamma^*_n(0)$, inequality (\ref{phi-lam-rate}) with $q=0$ gives
\begin{align}
\phi_n \{\gamma^*_n(0)+\nu_n\}^2 &\le O(1) \left\{ n^{1/2} \Gamma_n W_n \gamma^{*2}_n(0) + n^{1/2} V_n \gamma^*_n(0) \nu_n + n V_n^2 \gamma^{*2}_n(0) \nu_n^2 \right\} . \label{prf-cor-fast-equ1}
\end{align}
By (\ref{Vn-equ}) and $\gamma^*_n(0)=B^*_0 n^{-1/2} w^*_n(0)^{-\beta_0/2}$, simple manipulation gives
%\begin{align*}
%& n^{1/2} V_n \gamma^*_n(0) \\
%&= {B_0^*}^{-(1/\mym_0-1/2)2 \tau_0}\{ n^{1/2} \gamma^*_n(0)\}^{2\tau_0/\beta_0}  \\
%&={B_0^*}^{-(1/\mym_0-1/2)2 \tau_0} [B_0^* \min\{ \gamma_n(0)^{-\beta_0/2}, \nu_n^{-\beta_0/2}\} ]^{2\tau_0/\beta_0} \\
%&={B_0^*}^{-(1/\mym_0-1/2)2 \tau_0}{B_0^*}^{2\tau_0 /\beta_0} \min\{ \gamma_n(0)^{-\tau_0}, \nu_n^{-\tau_0} \} \\
%&=B_0^*  \min\{ \gamma_n(0)^{-\tau_0}, \nu_n^{-\tau_0} \},
%\end{align*}
%and hence
%\begin{align}
%& n^{1/2} V_n \gamma^*_n(0) \nu_n
%=B_0^*  \min\{ \nu_n \gamma_n(0)^{-\tau_0}, \nu_n^{1-\tau_0} \}, \label{prf-cor-fast-equ2}
%\end{align}
\begin{align}
& n^{1/2} V_n \gamma^*_n(0) \nu_n = B^*_0 w^*_n(0)^{-\tau_0} \nu_n, \label{prf-cor-fast-equ2}  \\
& n^{1/2} W_n \gamma^*_n(0)^2 = B^*_0 w^*_n(0)^{-(1-\beta_0/2)\tau_0}\gamma^*_n(0). \nonumber
\end{align}
Then (\ref{cor-fast-equ}) and (\ref{prf-cor-fast-equ1}) directly imply that Assumption~\ref{rate-assumption}(i) holds for sufficiently large $n$ and also (\ref{main2-cor3-equ}) holds.
The desired result follows from Corollary~\ref{main2-cor3} with $q=0$. \hfill $\Box$

\vspace{.1in}
\textbf{Proof of Proposition~\ref{cor-medium}.}
For $\gamma_{nj}=\gamma^*_n(q)$, inequality (\ref{phi-lam-rate}) with $q=0$ gives
\begin{align}
\phi_n \{\gamma^*_n(q)+\nu_n\}^2 &\le O(1) \left\{ n^{1/2} \Gamma_n W_n \gamma^{*2}_n(q) + n^{1/2} V_n \gamma^*_n(q) \nu_n + n V_n^2 \gamma^{*2}_n(q) \nu_n^2 \right\} . \label{prf-cor-medium-equ1}
\end{align}
By (\ref{Vn-equ}) and $\gamma^*_n(q)=B^*_0 n^{-1/2} w^*_n(q)^{-\beta_0/2}$, simple manipulation gives
%\begin{align*}
%& n^{1/2} V_n \gamma^*_n(q)  \\
%&={B_0^*}^{-(1/\mym_0-1/2)2\tau_0}  \{ n^{1/2} \gamma^*_n(q)\}^{2\tau_0/\beta_0} \\
%&={B_0^*}^{-(1/\mym_0-1/2)2\tau_0} [B_0^* \min\{ \gamma_n(q)^{-(1-q)\beta_0/2}, \nu_n^{-(1-q)\beta_0/2}\} ]^{2\tau_0/\beta_0}\\
%&={B_0^*}^{-(1/\mym_0-1/2)2\tau_0} (B_0^*)^{2\tau_0 /\beta_0} \min\{ \gamma_n(q)^{-(1-q)\tau_0}, \nu_n^{-(1-q)\tau_0} \}\\
%&= B_0^* \min\{ \gamma_n(q)^{-(1-q)\tau_0}, \nu_n^{-(1-q)\tau_0} \}
%\end{align*}
%and hence
%\begin{align*}
%n^{1/2} V_n \gamma^*_n(q) \nu_n^{1-q} = B_0^* \min\{\nu_n^{1-q}\gamma_n(q)^{-(1-q)\tau_0}, \nu_n^{(1-q)(1-\tau_0)} \}.
%\end{align*}
\begin{align*}
& n^{1/2} V_n \gamma^*_n(q) \nu_n^{1-q} = B^*_0 w^*_n(q)^{-\tau_0} \nu_n^{1-q}, \\
& n^{1/2} W_n \gamma^*_n(q)^{2-q} = B^*_0 w^*_n(q)^{-(1-\beta_0/2)\tau_0}\gamma^*_n(q)^{1-q}.
\end{align*}
Then (\ref{cor-medium-equ}) and (\ref{prf-cor-medium-equ1}) imply that Assumption~\ref{rate-assumption}(i) holds for sufficiently large $n$, along with the fact that $\nu_n=o(1)$, $\gamma_n(q)=o(1)$, and $q>0$.
Moreover, (\ref{cor-medium-equ}) and (\ref{phi-lam-rate}) with $\gamma_{nj}=\gamma^*_n(q)$ directly yield (\ref{main2-cor4-equ}).
The desired result follows from Corollary~\ref{main2-cor4}.\hfill $\Box$

\vspace{.1in}
\textbf{Proof of Proposition~\ref{cor-fast2}.}
Denote by $\gamma^\prime_{n,p+1}$, $V_n^\prime$, $W_n^\prime$, etc., the corresponding quantities based on $(w^\prime_{nj}, \gamma^\prime_{nj})$.
By (\ref{Vn-equ}) and (\ref{prf-cor-fast-equ2}) with $\tau_0=1$, we have
$n^{1/2} V^\prime_n \gamma^\prime_{n,p+1} \nu_n = K_0^{-1} (n^{1/2} V_n \gamma_{n,p+1} \nu_n)$ and
$n^{1/2} V_n \gamma_{n,p+1} \nu_n =B_0^* \min\{\nu_n \gamma_n^{-1}(0), 1\} \le B_0^*$.
Moreover, we have
$n^{1/2} \Gamma_{n} W^\prime_n {\gamma^{\prime}}^2_{n,p+1} = o(1)$ for a constant $K_0$, because
$W^\prime_n = 1$, $\Gamma_n$ is no greater than $O(\log^{1/2}( n))$, and
$n^{1/2} \gamma_n^2(0)$ decreases polynomially in $n^{-1}$.
For a constant $0<\eta_1<1$, we choose and fix $K_0 \ge 1$ sufficiently large, depending on $\overline M$ but independently of $(n,p)$, such that
Assumptions~\ref{rate-assumption}(i)--(ii) are satisfied,
with $(w_{nj}, \gamma_{nj})$ replaced by $(w^\prime_{nj}, \gamma^\prime_{nj})$,
for sufficiently large $n$, due to (\ref{S-bound1-Q}), (\ref{prf-cor-fast-equ1}), and the definition $\lambda^\prime_{nj} = C_1(\gamma^\prime_{nj}+\nu_n)$.
Moreover, by (\ref{w-gam}), $\rho^\prime_{nj} = \lambda^\prime_{nj} w^\prime_{nj} \le K_0^{1-\beta_0/2} \lambda_{nj} w_{nj}
\le K_0^{1-\beta_0/2} C_1 \{\gamma^*_n(0)+\nu_n\}^2$, which together with (\ref{S-bound2-Q}) implies that
(\ref{main2-cor2-equ}) is satisfied for some constant $\eta_2>0$.
Assumption~\ref{compat-mono} is also satisfied with $C^*_0$ replaced by $C_0^{*\prime} = C^*_0 K_0^{\beta_0/2}$ and $S$ replaced by
$\{1\le j\le p: \|g^*_j \|_Q > C_0^{*\prime} \lambda^\prime_{nj}\}\subset S$ for $K_0 \ge 1$
due to monotonicity in $S$ for the validity of Assumption~\ref{compat-mono} by Remark~\ref{mono-S}, and
with $(w_{nj},\gamma_{nj})$ replaced by $(w^\prime_{nj},\gamma^\prime_{nj})$ because (\ref{compat-Q-equ1}) after the modification implies (\ref{compat-Q-equ1}) itself,
with $w^\prime_{nj} \ge w_{nj}$ for $K_0\ge 1$ and $\lambda_{nj}$ constant in $j$.
The desired result follows from Corollary~\ref{main2-cor2} with $\bar g=g^*$. \hfill $\Box$

\subsection{Proof of Theorem \ref{thm-Hilbert}}

We use the non-commutative Bernstein inequality (Lemma~\ref{Bernstein-ineq}) to prove Theorem \ref{thm-Hilbert}.
Suppose that $(X_1,\ldots,X_n)$ are independent variables in a set $\Omega$.
First, consider finite-dimensional functional classes $\mathcal F_j$ with elements of the form
\begin{align} \label{group}
f_j(x) = u_j^\T(x) \theta_j, \quad \forall\, \theta_j \in \mathbb R^{d_j}, j=1,2,
\end{align}
where $u_j(x)$ is a vector of basis functions from $\Omega$ to $\mathbb R^{d_j}$,
and $\theta_j$ is a coefficient vector.
Let $U_j = \{u_j(X_1),\ldots,u_j(X_n)\}^\T $, and
$\Sigma_{jj^\prime}= E \big(U_{j}^\T U_{j^\prime}/n\big)\in {\mathbb R}^{d_{j}\times d_{j^\prime}}$.
The population inner product is
$\langle f_j, f_j^\prime \rangle_Q
= \theta_{j}^\T \Sigma_{jj^\prime}\theta_{j^\prime}$, $j, j^\prime=1,2$.
The difference between the sample and population inner products can be written as
\begin{align*}
&\sup_{\|\theta_{j}\|=\|\theta_{j^\prime}\|=1} \Big|\big\langle f_{j},f_{j^\prime}\big\rangle_n
- \big\langle f_{j},f_{j^\prime}\big\rangle_{Q}\Big|
= \sup_{\|\theta_{j}\|=\|\theta_{j^\prime}\|=1}
| \theta_j^\T (U_{j}^\T U_{j^\prime}/n - \Sigma_{jj^\prime} ) \theta_j^\prime | \\
%= \sup_{\|\theta_{j^\prime}\|=1} \| (U_{j}^\T U_{j^\prime}/n - \Sigma_{jj^\prime} ) \theta_j^\prime \|
&  = \| U_{j}^\T U_{j^\prime}/n - \Sigma_{jj^\prime}\|_S.
\end{align*}
%where the second step holds by the Cauchy--Schwartz inequality.

\begin{lem}\label{lm-a}
Let $f_j$ be as in (\ref{group}).
Assume that for a constant $C_{5,1}$,
\begin{align*} %\label{group-cond}
\sup_{x\in \Omega} \|u_j(x)\|^2 \le C_{5,1} \ell_j , \quad \forall j=1,2.
\end{align*}
Then for all $t>0$,
\begin{align*}
\| U_{j}^\T U_{j^\prime}/n - \Sigma_{jj^\prime}\|_S
> \sqrt{(\ell_{j}\|\Sigma_{j^\prime j^\prime}\|_{S}) \vee (\ell_{j^\prime}\|\Sigma_{jj}\|_{S})} \sqrt{\frac{2C_{5,1}t}{n}}
+ C_{5,1}\sqrt{\ell_{j}\ell_{j^\prime}} \frac{4t}{3n}
\end{align*}
with probability at least $1-(d_{j}+d_{j^\prime}) \me^{-t}$.
\end{lem}

\begin{prf}
Let $M_i = u_{j}(X_i)u_{j^\prime}^\T (X_i) - E \{u_{j}(X_i)u_{j^\prime}^\T (X_i)\}$.
Because $u_{j}(X_i)u_{j^\prime}^\T (X_i)$ is of rank 1,
$\|M_i\|_S\le 2\sup_{x\in\Omega} \{ \|u_{j}(x)\|\|u_{j^\prime}^\T (x)\| \} \le 2C_{5,1}\sqrt{\ell_{j}\ell_{j^\prime}}$.
Hence we set $s_0 = 2C_{5,1}\sqrt{\ell_{j}\ell_{j^\prime}}$ in Lemma~\ref{Bernstein-ineq}.
Similarly, $W_{\rm col} \le C_{5,1}\ell_{j^\prime} \|\Sigma_{jj}\|_{S}$ because
\begin{align*}
E  (M_iM_i^\T ) \le E \{u_{j}(X_i)u_{j^\prime}^\T (X_i)u_{j^\prime}(X_i)u_{j}^\T (X_i)\}
\le C_{5,1}\ell_{j^\prime}\,E \{u_{j}(X_i)u_{j}^\T (X_i)\},
\end{align*}
and $W_{\rm row} \le C_{5,1}\ell_{j} \|\Sigma_{j^\prime j^\prime}\|_{S}$.
Thus, (\ref{matrix-exp-2}) gives the desired result.
\end{prf}

Now consider functional classes $\mathcal F_j$ such that $f_j \in \mathcal F_j$ admits an expansion
$$
f_j(\cdot) =\sum_{\ell=1}^\infty \theta_{j\ell} u_{j\ell}(\cdot) ,
$$
where $\{u_{j\ell}(\cdot): \ell=1,2,\ldots\}$ are basis functions and $\{\theta_{j\ell}: \ell=1,2,\ldots\}$ are the associated coefficients.

\begin{lem}\label{lm-c}
Let $0< \tau_j < 1$, $0 < w_{nj} \le 1$ and
\[
B_j=\Big\{f_j: \hbox{$\sum$}_{k/4 < \ell \le k} \theta_{j\ell}^2 \le k^{-1/\tau_j}\ \forall\, k\ge (1/w_{nj})^{2\tau_j}, \
\hbox{$\sum$}_{0\le \ell^{1/\tau_j} w_{nj}^2 < 1} \theta_{j,\ell+1}^2 \le w_{nj}^{2}\Big\}
\]
Suppose that (\ref{new-cond-2}) and (\ref{new-cond-1c}) hold with certain positive
constants $C_{5,1}$, $C_{5,3}$.
Then, for a certain constant $C_{5,4}$ depending on
$\{C_{5,1},C_{5,3}\}$ only,
\begin{align*}
& \sup_{f_j\in B_j ,f_{j^\prime}\in B_{j^\prime} }
\Big|\langle f_j,f_{j^\prime}\rangle_n - \langle f_j,f_{j^\prime}\rangle_{Q}\Big| \\
&\le C_{5,4}w_{nj} w_{nj^\prime}\Big[
(\mu_j w_{nj}^{-\tau_j} + \mu_{j^\prime}w_{nj^\prime}^{-\tau_{j^\prime}})
\sqrt{\big\{ \mu_j+ \mu_{j^\prime}+\log(w_{nj}^{-\tau_j}+w_{nj^\prime}^{-\tau_{j^\prime}})
+t\big\}/n} \nonumber \\
& \qquad\qquad\qquad
+\big\{ \mu_j+ \mu_{j^\prime}+\log(w_{nj}^{-\tau_j}+w_{nj^\prime}^{-\tau_{j^\prime}})  +t\big\}
(\mu_j w_{nj}^{-\tau_j})(\mu_{j^\prime}w_{nj^\prime}^{-\tau_{j^\prime}})/n \Big]
\end{align*}
with at least probability $1-e^{-t}$ for all $t>0$,
where $ \mu_j=1/(1-\tau_j)$. % and $ \mu_{j^\prime}=1/(1-\tau_{j^\prime})$.
\end{lem}

\begin{prf}
Let $\ell_{jk} =\lceil (2^k/w_{nj})^{2\tau_j} \rceil$.
We group the basis and coefficients as follows:
\begin{align*}
u_{j,G_{jk}}(x) = (u_{j\ell}(x), \ell\in G_{jk})^\T ,\quad \theta_{j,G_{jk}} = (\theta_{j\ell},\ell\in G_{jk})^\T ,\quad
k = 0,1,\ldots
\end{align*}
where $G_{j0}=\{1,\ldots,\ell_{j0}\}$ of size $|G_{j0}| = \ell_{j0}$ and
$G_{jk}=\{\ell_{j,k-1}+1,\ldots,\ell_{jk}\}$ of
size $|G_{jk}| = \ell_{jk} - \ell_{j,k-1}\le (2^k/w_{nj})^{2\tau_j}$ for $k \ge 1$.
Define $\tilde\theta_j$, a rescaled version of $\theta_j$, by
\begin{align*}
\tilde\theta_{j,G_{jk}} = (\tilde\theta_{j\ell},\ell\in G_{jk}) = 2^k w_{nj}^{-1}\theta_{j,G_{jk}} .
\end{align*}
%Because $\ell_{jk}^{1/(2\tau_j)}\ge 2^k/w_{nj}$ and
%$\ell_{jk} \le  1+ 2^{2\tau_j}\ell_{j,k-1}$ with $\tau_j<1$, we have $\ell_{jk} \le 4\ell_{j,k-1}$ and
It follows directly from (\ref{new-cond-1a}) and (\ref{new-cond-1b}) that
\begin{align*}
\|\tilde\theta_{j,G_{j,0}}\|_2 \le 1,\quad
\|\tilde\theta_{j,G_{jk}}\|_{r_j} \le (2^k/w_{nj})/\ell_{jk}^{-1/(2\tau_j)}\le 1\ \forall\ k\ge 1,\ \forall\, f_j\in B_j.
\end{align*}
Let $U_{jk}=\{ u_{j,G_{jk}}(X_1),\ldots,u_{j,G_{jk}}(X_n) \}^\T \in {\mathbb R}^{n\times |G_{jk}|}$. We have
\begin{align} \label{new-pf-1}
& \sup_{f_j\in B_j,f_{j^\prime}\in B_{j^\prime}}
\Big|\langle f_j,f_{j^\prime}\rangle_n - \langle f_j,f_{j^\prime}\rangle_{L_2}\Big| \nonumber \\
&= \sup_{f_j\in B_j,f_{j^\prime}\in B_{j^\prime}}
\Bigg|\sum_{k=0}^\infty \sum_{\ell=0}^\infty
\theta_{j,G_{jk}}^\T \bigg(U_{jk}^\T U_{j^\prime,\ell}/n - E \,U_{jk}^\T U_{j^\prime,\ell}/n\bigg)\theta_{j^\prime,G_{j^\prime,\ell}}\Bigg| \nonumber \\
&\le \max_{\|\tilde\theta_j\| \vee \|\tilde\theta_{j^\prime}\|\le 1}
\Bigg|\sum_{k=0}^\infty \sum_{\ell=0}^\infty
\tilde\theta_{j,G_{jk}}^\T \bigg(\frac{U_{jk}^\T U_{j^\prime,\ell}/n - E \,U_{jk}^\T U_{j^\prime,\ell}/n}
{2^k w_{nj}^{-1} 2^\ell w_{nj^\prime}^{-1} }\bigg)\tilde\theta_{j^\prime,G_{j^\prime,\ell}}\Bigg| \nonumber \\
&\le  w_{nj} w_{nj^\prime} \sum_{k=0}^\infty \sum_{\ell=0}^\infty
\bigg\|\frac{U_{jk}^\T U_{j^\prime,\ell}/n - E \,U_{jk}^\T U_{j^\prime,\ell}/n}{2^{k} 2^{\ell}}\bigg\|_S.
\end{align}
Let $a_k = 1/\{(k+1)(k+2)\}$. By (\ref{new-cond-2}),
$\sup_{x \in \Omega} \| u_{j,G_{jk}}(x) \|^2 \le \sup_{x\in\Omega}\sum_{\ell=1}^{\ell_{jk}} u_{j\ell}^2(x)
\le C_{5,1}\ell_{jk}$ for $k \ge 0$.
By (\ref{new-cond-1c}), $\| E U^\T_{jk}  U_{jk}/n \|_S \le C_{5,3}$.
Because $|G_{jk}|\le \ell_{j,k}$, it follows from Lemma \ref{lm-a} that
\begin{align} \label{new-pf-2}
& \|U_{jk}^\T U_{j^\prime,\ell}/n - E \,U_{jk}^\T U_{j^\prime,\ell}/n\|_S \nonumber \\
\le &\sqrt{ \big\{\log(\ell_{jk}+\ell_{j^\prime,\ell}) - \log(a_ka_\ell)+t \big\}
2C_{5,1}C_{5,3} (\ell_{jk} \vee \ell_{j^\prime,\ell})/n}\nonumber \\
& + \big\{\log(\ell_{jk}+\ell_{j^\prime,\ell}) - \log(a_ka_\ell)+t \big\}
(4/3)C_{5,1}\sqrt{\ell_{jk}\ell_{j^\prime,\ell}}/n
\end{align}
with probability at least $1-a_ka_\ell \me^{-t}$ for any fixed $k\ge 0$ and $\ell \ge 0$.
By the union bound and the fact that $\sum_{k=0}^\infty a_k=1$, inequality (\ref{new-pf-2}) holds simultaneously
for all $k\ge 0$ and $\ell\ge 0$ with probability at least $1- \me^{-t}$.
Because $\ell_{jk} =\lceil (2^k/w_{nj})^{2\tau_j} \rceil$, we rewrite (\ref{new-pf-2}) as
\begin{align} \label{new-pf-3}
& \|U_{jk}^\T U_{j^\prime,\ell}/n - E \,U_{jk}^\T U_{j^\prime,\ell}/n\|_S \nonumber \\
\le & C_{5,4}\Big[(2^{\tau_j k} w_{nj}^{-\tau_j}
+ 2^{\tau_{j^\prime}\ell} w_{nj^\prime}^{-\tau_{j^\prime}})
\sqrt{ \big\{k+\ell+\log(w_{nj}^{-\tau_j}+w_{nj^\prime}^{-\tau_{j^\prime}})  +t\big\}/n} \nonumber \\
& \qquad
+ \big\{k+\ell+\log(w_{nj}^{-\tau_j}+w_{nj^\prime}^{-\tau_{j^\prime}})  +t\big\}
(2^{\tau_j k}w_{nj}^{-\tau_j})(2^{\tau_{j^\prime}\ell}w_{nj^\prime}^{-\tau_{j^\prime}})/n\Big].
\end{align}
where $C_{5,4}$ is a constant depending only on $\{C_{5,1},C_{5,3}\}$.
For any $\alpha \ge 0$, $\sum_{k=0}^\infty k^\alpha 2^{-k(1-\tau_j)}
\le C_\alpha  \mu_j ^{\alpha+1}$, where $C_\alpha$ is a numerical constant and $ \mu_j = 1/(1-\tau_j)$.
Using this fact and inserting (\ref{new-pf-3}) into (\ref{new-pf-1}) yields the desired result.
\end{prf}

Finally, the following result concludes the proof of Theorem~\ref{thm-Hilbert}.

\begin{lem}
In the setting of Theorem~\ref{thm-Hilbert}, let
\begin{align*}
& \phi_n = C_{5,2} C_{5,4}
\left\{ \max_j \frac{\sqrt{2\log(np/\epsilon^\prime)}}{\lambda_{nj}} \max_j \frac{\mu_j \tilde \gamma_{nj}}{\lambda_{nj}} + \max_j \frac{2\log(np/\epsilon^\prime) \mu_j^2 \tilde \gamma_{nj}^2}{\lambda_{nj}^2} \right\},
\end{align*}
where $\tilde\gamma_{nj} = n^{-1/2} w_{nj}^{-\tau_j}$, $\mu_j =1/(1-\tau_j)^{-1} $, and $C_{5,4}$ is a constant depending only on $\{C_{5,1},C_{5,3}\}$ as in
Lemma \ref{lm-c}.
Then
\begin{align*}
P\left\{ \sup_{g\in \mathcal G} \frac{ \left| \|g\|_n^2 - \|g \|_Q^2 \right|}{ R^{*2}_n(g) } > \phi_n \right\} \le {\epsilon^\prime}^2 .
\end{align*}
\end{lem}

\begin{prf}
Recall that $\ell_{jk} =\lceil (2^k/w_{nj})^{2\tau_j} \rceil$.
For $g_j =\sum_{\ell=1}^\infty \theta_{j\ell} u_{j\ell}$, define $r_{nj}(g_j)$ by
\[
r^2_{nj}(g_j) = \bigg(\sum_{\ell=1}^{\ell_{j0}} \theta_{j\ell}^2/w_{nj}^{2}\bigg)\vee
\bigg( \max_{k\ge 1}\sum_{\ell_{j,k-1} < \ell \le \ell_{jk}} \theta_{j\ell}^2 \ell_{jk}^{1/\tau_j} \bigg).
\]
Let $f_j = g_j/r_{nj}(g_j)$ and $ \mu_j=1/(1-\tau_j)$.
Then $f_j \in B_j$ as in Lemma~\ref{lm-c} and
\begin{align*}
& \left| \|g\|_n^2 - \|g \|_Q^2 \right|
\le \sum_{j=1}^p\sum_{j'=1}^p\Big|\langle g_j,g_{j^\prime}\rangle_n - \langle g_j,g_{j^\prime}\rangle_{Q}\Big|  \\
& = \sum_{j=1}^p\sum_{j'=1}^p r_{nj}(g_j) r_{nj^\prime} (g_{j^\prime})
\Big|\langle f_j,f_{j^\prime}\rangle_n - \langle f_j,f_{j^\prime}\rangle_{Q}\Big|.
\end{align*}
Because $\sum_{i=1}^p w_{nj} \lambda_{nj} r_{nj}(g_j)
\le \sum_{j=1}^p C_{5,2}^{1/2}\lambda_{nj}   (w_{nj}\|g_j\|_{F,j} + \|g_j\|_{Q})
= C_{5,2}^{1/2}R^*_n(g)$ by (\ref{new-cond-1a}),
\begin{align} \label{pf-th4-1}
& \left\{\sup_{g\in \mathcal G} \frac{ \left| \|g\|_n^2 - \|g \|_Q^2 \right|}{ R^{*2}_n(g) } > \phi_n \right\} \nonumber \\
& \subset \bigcup_{j,j^\prime} \left\{ \sup_{f_j \in B_j, f_{j^\prime} \in B_{j^\prime} }
\left| \langle f_j, f_k \rangle_n - \langle f_j, f_k \rangle_Q \right|
> C_{5,2}^{-1} \phi_n \, w_{nj} \lambda_{nj} w_{nj^\prime} \lambda_{nj^\prime} \right\}.
\end{align}
By Lemma \ref{lm-c} with $t=\log(p^2/{\epsilon^\prime}^2)$ and $e^{2 \mu_j}+ 2w_{nj}^{-\tau_j}\le n$,
we have
\begin{align*}
&\sup_{f_j\in B_j ,f_{j^\prime}\in B_{j^\prime}}
\Big|\langle f_j,f_{j^\prime}\rangle_n - \langle f_j,f_{j^\prime}\rangle_{Q}\Big| \\
&\le C_{5,4}w_{nj} w_{nj^\prime}\Big[
(\mu_j w_{nj}^{-\tau_j} + \mu_{j^\prime}w_{nj^\prime}^{-\tau_{j^\prime}})
\sqrt{\big\{ \mu_j+ \mu_{j^\prime}+\log(w_{nj}^{-\tau_j}+w_{nj^\prime}^{-\tau_{j^\prime}})
+\log(p^2/{\epsilon^\prime}^2)\big\}/n} \nonumber \\
& \qquad\qquad\qquad
+\big\{ \mu_j+ \mu_{j^\prime}+\log(w_{nj}^{-\tau_j}+w_{nj^\prime}^{-\tau_{j^\prime}})
+\log(p^2/{\epsilon^\prime}^2)\big\}
(\mu_j w_{nj}^{-\tau_j})(\mu_{j^\prime}w_{nj^\prime}^{-\tau_{j^\prime}})/n \Big]\\
&\le C_{5,4}w_{nj} w_{nj^\prime}\Big\{
(\mu_j w_{nj}^{-\tau_j} + \mu_{j^\prime}w_{nj^\prime}^{-\tau_{j^\prime}})
\sqrt{2\log(np/{\epsilon^\prime})/n}
+ 2\log(p/{\epsilon^\prime})(\mu_j w_{nj}^{-\tau_j})(\mu_{j^\prime}w_{nj^\prime}^{-\tau_{j^\prime}})/n \Big\},
\end{align*}
with probability at least $1-{\epsilon^\prime}^2/p^2$.
By the definition of $\phi_n$, we have
\begin{align*}
P \left\{ \sup_{f_j\in B_j ,f_{j^\prime}\in B_{j^\prime}}
\Big|\langle f_j,f_{j^\prime}\rangle_n - \langle f_j,f_{j^\prime}\rangle_{Q}\Big|\le C_{5,2}^{-1} \phi_n w_{nj} w_{nj^\prime} \lambda_{nj} \lambda_{nj^\prime} \right\} \ge 1-\frac{{\epsilon^\prime}^2}{p^2}.
\end{align*}
The conclusion follows from the union bound using (\ref{pf-th4-1}).
\end{prf}

\subsection{Proof of Proposition \ref{cor-fast0}}

Here we verify explicitly conditions of Theorem \ref{thm-Hilbert}
for Sobolev spaces ${\mathcal W}_{r_i}^{m_i}$ and bounded variation spaces ${\mathcal{V}}^{m_i}$
with $r_i=1$ on $[0,1]$ in the case of $\tau_j<1$, where $\tau_i = 1/(2m_i+1-2/(r_i\wedge 2))$.
Because conditions (\ref{new-cond-2}), (\ref{new-cond-1a}),
(\ref{new-cond-1b}) and (\ref{new-cond-1c}) depend on $(m_j,r_j)$ only through $\tau_j$,
we assume without loss of generality $1\le r_j\le 2$.
When the average marginal density of $\{X_i^{(j)}: i=1,\ldots,n\}$ is uniformly bounded away from $0$ and $\infty$,
the norms $\|g_j\|_Q$ and $\|g_j\|_{L_2}$ are equivalent, so that
condition (\ref{new-cond-1b}) and (\ref{new-cond-1c}) hold for any $L_2$-orthonormal bases
$\{u_{j\ell}: \ell\ge 1\}$.
Let $u_0(x)$ be a mother wavelet with $m$ vanishing moments, e.g.,
$u_0(x)=0$ for $|x|>c_0$, $\int u_0^2(x)dx=1$, $\int  x^m u_0(x) dx = 0$ for $m = 0,\ldots,\max_j m_j$,
and
$\{u_{0,k\ell}(x) = \sqrt{2^k}u_0(2^k(x-j)): \ell =1,\ldots,2^k, k=0,1,\ldots\}$ is $L_2$-orthonormal.
We shall identify $\{u_{j\ell}: \ell\ge 1\}$ as
$\{u_{0,11}, u_{0,21}, u_{0,22}, u_{0,31}, \ldots\}$.
Because $\#\{\ell: u_{0,k\ell}(x)\neq 0\}\le 2c_0k\ \forall x$,
\[
\sum_{\ell=2^k}^{2^{k+1}-1} u_{j\ell}^2(x) = \sum_{\ell=1}^{2^k}u^2_{0,k\ell}(x)
\le 2c_0 2^k\|u_0\|_\infty, \ \forall x,
\]
so that (\ref{new-cond-2}) holds.
Suppose $g_j(x) = \sum_{\ell=1}^\infty \theta_{j\ell}u_{j\ell}(x)
= \sum_{k=0}^\infty \sum_{\ell=1}^{2^k} \theta_{jk\ell} u_{0,k\ell}(x)$.
Define $u_0^{(-m)}(x)$ as the $m$-th integral of $u_0$,
$u_0^{(-m)}(x) = \int_{-\infty}^x u_0^{(-m+1)}(t)dt$, and $g_j^{(m)}(x)=(d/dx)^mg_j(x)$.
Because $u_0$ has vanishing moments, $\int u_0^{(m)}(x)dx=0$ for $m=0,\ldots,\max_j m_j$,
so that $u_0^{(m)}(x)=0$ for $|x|>c_0$.
Due to the orthonormality of the basis functions, for $1 \le \ell \le 2^k$, we have
\[
2^{m_jk}\theta_{jk\ell} = 2^{m_jk}\int g_j(x)u_{0,k\ell}(x)dx = (-1)^m\int g_j^{(m_j)}(x)u_{0,m_jk\ell}(x)dx
\]
with $u_{0,mk\ell}(x) = \sqrt{2^k}u_0^{(-m)}(2^k(x-j))$.
By the H\"older inequality,
\begin{align*}
\sum_{\ell=2^{k-1}}^{2^k-1}\big| 2^{m_jk}\theta_{j\ell}\big|^{r_j}
&\le \sum_{\ell=2^{k-1}}^{2^k-1}\int \Big| g_j^{(m_j)}(x)\Big|^{r_j}\Big|u_{0,m_jk\ell}(x)\Big|^{r_j(1-2(1-/r_j))} dx
\Big\|u_{0,m_jk\ell}\Big\|_{L_2}^{r_j(2(1-1/r_j))}
\cr &\le \Big\| g_j^{(m_j)}\Big\|^{r_j}_{L_{r_j}} 2c_0\Big\|u_0^{(m_j)}\Big\|_\infty^{r_j(1-2(1-/r_j))}
2^{(k/2)r_j(1-2(1-/r_j))}
\Big\|u_0^{(m_j)}\Big\|_{L_2}2^{r_j(2(1-1/r_j))}.
\end{align*}
Because $2^{m_jk - (k/2)(1-2(1-/r_j))} = 2^{k(m_j+1/2 - 1/r_j)} = 2^{k/(2\tau_j)}$ and $1\le r_j\le 2$, we have
\begin{align*}
\bigg\{\sum_{\ell=2^{k-1}}^{2^k-1}2^{k/\tau_j}\theta_{j\ell}^2\bigg\}^{1/2}
&\le \bigg\{\sum_{\ell=2^{k-1}}^{2^k-1}\big| 2^{k/(2\tau_j)}\theta_{j\ell}\big|^{r_j}\bigg\}^{1/r_j}
\cr &\le \Big\|g_j^{(m_j)}\Big\|_{L_{r_j}} (2c_0)^{1/r_j}\Big\|u_0^{(m_j)}\Big\|_\infty^{2/r_j-1}
\Big\|u_0^{(m_j)}\Big\|_{L_2}^{2-2/r_j} .
\end{align*}
%which implies (\ref{new-cond-1a}).
Because $\ell_{jk}^{1/(2\tau_j)}\ge 2^k/w_{nj}$ and
$\ell_{jk} \le  1+ 2^{2\tau_j}\ell_{j,k-1}$ with $\tau_j<1$, we have $\ell_{jk} \le 4\ell_{j,k-1}$,
so that $\{\ell_{j,k-1}+1,\ldots,\ell_{j,k}\}$ involves at most three resolution levels.
Thus, condition (\ref{new-cond-1a}) follows from the above inequality.
For the bounded variation class, we have
\[
2^{m_jk}\theta_{jk\ell} = 2^{m_jk} \int g_j(x)u_{0,k\ell}(x)dx
= (-1)^m \int u_{0,m_jk\ell}(x)dg_j^{(m_j-1)}(x),
\]
so that (\ref{new-cond-1a}) follows from the same proof with $r_j=1$.

\section{Technical tools} \label{sect:tech-tools}

\subsection{Sub-gaussian maximal inequalities} \label{sect:preparation}

The following maximal inequality can be obtained from van de Geer (2000, Corollary~8.3),
or directly derived using Dudley's inequality for sub-gaussian variables and Chernoff's tail bound (see Proposition 9.2, Bellec et~al.~2016).

\begin{lem} \label{max-ineq}
%Let $(\varepsilon_1,\ldots,\varepsilon_n)$ be independent sub-Gaussian variables under Assumption~\ref{sub-gaussian-error}
For $\delta>0$, let $\mathcal F_1$ be a functional class such that $\sup_{f_1 \in \mathcal F_1 } \|f_1\|_n  \le \delta$,
and
\begin{align}
\psi_n (\delta, \mathcal F_1) \ge \int_0^\delta H^{1/2} (u, \mathcal F_1, \|\cdot\|_n) \,\dif u . \label{entropy-integral}
\end{align}
Let $(\varepsilon_1,\ldots,\varepsilon_n)$ be independent variables. Under Assumption~\ref{sub-gaussian-error}, we have for any $t>0$,
\begin{align*}
P \left\{ \sup_{f_1 \in \mathcal F_1 } |\langle\varepsilon, f_1\rangle_n | /C_1 > n^{-1/2} \psi_n(\delta, \mathcal F_1) + \delta \sqrt{t/n} \right\}
\le \exp( -t),
\end{align*}
where $C_1 =C_1(D_0,D_1)>0$ is a constant, depending only on $(D_0,D_1)$.
\end{lem}

\subsection{Dudley and Talagrand inequalities}

The following inequalities are due to Dudley (1967) and Talagrand (1996).

\begin{lem} \label{Dudley-thm}
For $\delta>0$, let $\mathcal F_1 $ be a functional class such that $\sup_{f_1 \in \mathcal F_1 } \|f_1\|_n  \le \delta$ and
(\ref{entropy-integral}) holds.
Let $(\sigma_1,\ldots,\sigma_n)$ be independent Rademacher variables, that is, $P(\sigma_i=1)=P(\sigma_i=-1)=1/2$.
Then for a universal constant $C_2>0$,
\begin{align*}
E \left\{ \sup_{f_1 \in \mathcal F_1 } |\langle\sigma, f_1\rangle_n | /C_2  \right\} \le n^{-1/2} \psi(\delta,\mathcal F_1) .
\end{align*}
\end{lem}

\begin{lem}\label{Talagrand-thm}
For $\delta>0$ and $b>0$, let $(X_1,\ldots,X_n)$ be independent variables, and
$\mathcal F$ be a functional class such that $\sup_{f \in \mathcal F} \|f\|_Q \le \delta$ and $\sup_{f \in \mathcal F} \|f\|_\infty \le b $.
Define
$$
Z_n = \sup_{f \in \mathcal F} \left| \frac{1}{n} \sum_{i=1}^n \{ f(X_i) - E f(X_i) \} \right|.
$$
Then for a universal constant $C_3>0$, we have
\begin{align*}
P \left\{ Z_n /C_3 > E(Z_n) + \delta \sqrt{\frac{t}{n}} + b \frac{t}{n} \right\} \le \exp(-t), \quad t>0.
\end{align*}
\end{lem}

\subsection{Non-commutative Bernstein inequality}

We state the non-commutative Bernstein inequality
%(Oliveira, 2010;
(Troop, 2011) as follows.

\begin{lem} \label{Bernstein-ineq}
Let $\{M_i: i=1,\ldots, n\}$ be independent random matrices in ${\mathbb R}^{d_1\times d_2}$
such that $E (M_i)=0$ and $P \{\| M_i\|_S \le s_0\}=1$, $i=1,\ldots,n$, for a constant $s_0>0$, where $\|\cdot\|_S$
denotes the spectrum norm of a matrix.
Let $\Sigma_{\rm col} = \sum_{i=1}^n E  (M_i M_i^\T)/n$
and $\Sigma_{\rm row} = \sum_{i=1}^n E  (M_i^\T M_i)/n$. Then, for all $t>0$,
\begin{align} \label{matrix-exp}
P \bigg( \bigg\|\frac{1}{n}\sum_{i=1}^n M_i\bigg\|_S > t\bigg)
\le (d_1+d_2) \exp\bigg( \frac{ - n t^2/2}{\|\Sigma_{\rm col}\|_{S}\vee\|\Sigma_{\rm row}\|_{S} + s_0 t/3}\bigg).
\end{align}
Consequently, for all $t>0$,
\begin{align} \label{matrix-exp-2}
P \bigg\{ \bigg\|\frac{1}{n}\sum_{i=1}^n M_i\bigg\|_S
> \sqrt{\|\Sigma_{\rm col}\|_{S}\vee\|\Sigma_{\rm row}\|_{S}}\sqrt{2t/n}
+ ( s_0/3)2t/n \bigg\}
\le (d_1+d_2)\me^{-t}.
\end{align}
\end{lem}

\subsection{Convergence of empirical norms}

For $\delta>0$ and $b>0$, let $\mathcal F_1$ be a functional class such that
\begin{align*}
\sup_{f_1 \in\mathcal F_1} \|f_1\|_Q \le \delta  , \quad \sup_{f_1 \in\mathcal F_1}  \|f_1 \|_\infty \le b ,
\end{align*}
and let $\psi_{n,\infty}(\cdot, \mathcal F_1)$ be an upper envelope of the entropy integral:
\begin{align}
\psi_{n,\infty} (z, \mathcal F_1) \ge \int_0^z {H^*}^{1/2} (u/2, \mathcal F_1, \|\cdot\|_{n,\infty} )\,\dif u, \quad z>0 , \label{entropy-bound}
\end{align}
where $H^* (u, \mathcal F_1, \|\cdot\|_{n,\infty}) = \sup_{(X_1^{(1)},\ldots, X_n^{(1)})} H (u, \mathcal F_1, \|\cdot\|_{n,\infty})$.
Let $\hat \delta = \sup_{f_1 \in \mathcal F_1} \|f_1\|_n$. The following result can be obtained from Guedon et~al.~(2007) and, in its present form, van de Geer (2014), Theorem 2.1.

\begin{lem} \label{empirical-norm}
For the universal constant $C_2$ in Lemma~\ref{Dudley-thm}, we have
\begin{align*}
E \left\{ \sup_{f_1 \in \mathcal F_1} \left| \|f_1\|_n^2 - \|f_1\|_Q^2 \right| \right\} \le
\frac{2 \delta C_2 \psi_{n,\infty}(b, \mathcal F_1)}{\sqrt n} + \frac{4 C_2^{2} \psi^2_{n,\infty}(b, \mathcal F_1) }{n}.
\end{align*}
Moreover, we have
\begin{align*}
\sqrt{ E( \hat \delta^2 ) } \le \delta + \frac{2 C_2 \psi_{n,\infty} (b, \mathcal F_1)}{\sqrt n} .
\end{align*}
\end{lem}

\subsection{Metric entropies} \label{sect:entropy}

For $r\ge 1$ and $\mym>0$ (possibly non-integral), let $\overline{\mathcal W}_r^\mym = \{f :\|f\|_{L_r} + \| f^{(\mym)}\|_{L_r}\le 1\}$.
The following result is taken from Theorem 5.2, Birman \& Solomjak (1967).

\begin{lem} \label{sobolev-entropy}
If $r \mym >1$ and $1\le q \le \infty$, then
\begin{align*}
H( u , \overline{\mathcal W}_r^\mym, \|\cdot\|_{L_q} ) \le B_1 u ^{-1/\mym}, \quad u  >0,
\end{align*}
where $B_1 =B_1(\mym,r)>0$ is a constant depending only on $(\mym, r)$. If $r \mym \le 1$ and $1\le q< r/(1-r \mym)$, then
\begin{align*}
H( u , \overline{\mathcal W}_r^\mym, \|\cdot\|_{L_q} ) \le B_2 u ^{-1/\mym}, \quad u  >0,
\end{align*}
where $B_2 =B_2(\mym,r,q)>0$ is a constant depending only on $(\mym,r,q)$.
\end{lem}

For $\mym \ge 1$, let $\overline{\mathcal{V}}^\mym = \{f: \|f\|_{L_1} + \TV(f^{(\mym-1)})\le 1\}$.
The following result can be obtained from Theorem 15.6.1, Lorentz et~al.~(1996), on
the metric entropy of the ball $\{f: \|f\|_{L_r} + [f]_{\mbox{\scriptsize Lip}(\mym,L_r)} \le 1\}$,
where $[f]_{\mbox{\scriptsize Lip}(\mym,L_r)}$ is a semi-norm in the Lipschitz space $\mbox{Lip}(\mym,L_r)$.
By Theorem 9.9.3, DeVore \& Lorentz (1993), the space $\mbox{Lip}(\mym,L_1)$ is equivalent to
$\mathcal V^\mym$, with the semi-norm $[f]_{\mbox{\scriptsize Lip}(\mym,L_r)}$ equal to $\TV(f)$, up to suitable
modification of function values at (countable) discontinuity points.
However, it should be noted that the entropy of $\overline{\mathcal{V}}^1$ endowed with the norm $\|\cdot\|_{L_\infty}$ is infinite.

\begin{lem} \label{sobolev-entropy2}
If $ \mym \ge 2$ and $1\le q \le \infty$, then
\begin{align*}
H( u , \overline{\mathcal V}^\mym, \|\cdot\|_{L_q} ) \le B_3 u ^{-1/\mym}, \quad u  >0,
\end{align*}
where $B_3 =B_3(\mym)>0$ is a constant depending only on $\mym$. If $1\le q < \infty$, then
\begin{align*}
H( u , \overline{\mathcal V}^1, \|\cdot\|_{L_q} ) \le B_4 u ^{-1}, \quad u  >0,
\end{align*}
where $B_4=B_4(r) >0$ is a constant depending only on $r$.
\end{lem}

By the continuity of functions in $\mathcal W_r^\mym$ for $\mym \ge 1$ and $\mathcal V^\mym$ for $\mym \ge 2$,
the maximum entropies of these spaces in $\|\cdot\|_{n,\infty}$ and $\|\cdot\|_n$ norms over all possible design points can be derived from Lemmas~\ref{sobolev-entropy} and \ref{sobolev-entropy2}.

\begin{lem} \label{sobolev-entropy3}
If $r \mym >1$, then for $B_1=B_1(\mym,r)$,
\begin{align*}
H^*( u , \overline{\mathcal W}_r^\mym, \|\cdot\|_{n} ) \le H^*( u , \overline{\mathcal W}_r^\mym, \|\cdot\|_{n,\infty} ) \le B_1 u ^{-1/\mym}, \quad u  >0,
\end{align*}
and hence (\ref{entropy-integral}) and (\ref{entropy-bound}) hold with $\psi_n(z, \overline{\mathcal W}_r^\mym)\asymp \psi_{n,\infty}(z, \overline{\mathcal W}_r^\mym) \asymp z^{1-1/(2\mym)}$.
If $\mym \ge 2$, then for $B_3=B_3(\mym)$,
\begin{align*}
H^*( u , \overline{\mathcal V}^\mym, \|\cdot\|_{n} ) \le H^*( u , \overline{\mathcal V}^\mym, \|\cdot\|_{n,\infty} ) \le B_3 u ^{-1/\mym}, \quad u  >0,
\end{align*}
and hence (\ref{entropy-integral}) and (\ref{entropy-bound}) hold with $\psi_n(z, \overline{\mathcal V}^\mym)
\asymp \psi_{n,\infty}(z, \overline{\mathcal V}^\mym)
\asymp z^{1-1/(2\mym)}$.
\end{lem}

The maximum entropies of $\overline{\mathcal{V}}^1$ over all possible design points
can be obtained from Section 5, Mammen (1991) for the norm $\|\cdot\|_n$ and Lemma 2.2, van de Geer (2000) for the norm $\|\cdot\|_{n,\infty}$.
In fact, the proof of van de Geer shows that for $\mathcal F$ the class of nondecreasing functions
$f: [0,1]\to [0,1]$, $H^*( u , \mathcal F, \|\cdot\|_{n,\infty} ) \le n\log(n+u^{-1})$ if $u\le n^{-1}$ or $\le u^{-1}\log(n+u^{-1})$ if $u >n^{-1}$.
But if $u\le n^{-1}$, then $ n\log(n+u^{-1}) \le n(\log n + n^{-1} u^{-1}) \le (1+\log n) u^{-1}$.
If $u > n^{-1}$, then $u^{-1}\log(n+u^{-1}) \le u^{-1}\log(2n)$. Combining the two cases gives the stated result about $H^*( u , \overline{\mathcal{V}}^1, \|\cdot\|_{n,\infty} )$,
because each function in $\overline{\mathcal{V}}^1$ can be expressed as a difference two nondecreasing functions.

\begin{lem} \label{TV-entropy}
For a universal constant $B_5>0$, we have
\begin{align*}
H^*( u , \overline{\mathcal W}_{1}^1, \|\cdot\|_{n} ) \le H^*( u , \overline{\mathcal{V}}^1, \|\cdot\|_n) \le B_5 u ^{-1}, \quad u >0,
\end{align*}
and hence (\ref{entropy-integral}) holds with $\psi_n(z, \overline{\mathcal W}_{1}^1) \asymp \psi_n(z, \overline{\mathcal{V}}^1) \asymp z^{1/2}$.
Moreover, for a universal constant $B_6>0$, we have
\begin{align*}
H^*( u , \overline{\mathcal W}_{1}^1, \|\cdot\|_{n,\infty} ) \le H^*( u , \overline{\mathcal{V}}^1, \|\cdot\|_{n,\infty} ) \le B_6 \frac{1+\log n}{u }, \quad u >0,
\end{align*}
and hence (\ref{entropy-bound}) holds with $\psi_{n,\infty}(z, \overline{\mathcal W}_{1}^1) \asymp \psi_{n,\infty}(z, \overline{\mathcal{V}}^1)
\asymp (1+\log n)^{1/2} (z/2)^{1/2}$.
\end{lem}

\subsection{Interpolation inequalities}

The following inequality (\ref{inter-ineq1}) can be derived from the Gagliardo-Nirenberg inequality for Sobolev spaces (Theorem 1, Nirenberg 1966). Inequality (\ref{inter-ineq2}) can be
shown by approximating $f \in \mathcal V^\mym$ by functions in $\mathcal W_1^\mym$.

\begin{lem} \label{sobolev-interpolation}
For $r \ge 1$ and $\mym \ge 1$, we have
for any $f \in \mathcal W_r^\mym$,
\begin{align}
\| f \|_{\infty} \le (C_4/2) \left\{ \|f^{(\mym)}\|_{L_r} + \|f\|_{L_2} \right\}^{\tau} \|f\|_{L_2}^{1-\tau} , \label{inter-ineq1}
\end{align}
where $\tau=(2\mym+1-2/r)^{-1} \le 1$ and $C_4=C_4(\mym,r) \ge 1$ is a constant depending only on $(\mym,r)$.
In addition, we have for any $f \in \mathcal V^\mym$,
\begin{align}
\| f \|_{\infty} \le (C_4/2) \left\{ \TV(f^{(\mym-1)}) + \|f\|_{L_2} \right\}^{\tau} \|f\|_{L_2}^{1-\tau} . \label{inter-ineq2}
\end{align}
\end{lem}

From this result, $\|f\|_\infty$ can be bounded in terms of $\|f\|_{L_2}$ and $\| f^{(\mym)} \|_{L_r}$ or $\TV(f^{(\mym-1)})$ in a convenient manner.
For $f \in \mathcal W_r^\mym$ and $0< \delta \le 1$,
if $\|f\|_{L_2} \le \delta$ and $\| f^{(\mym)} \|_{L_r} \le 1$, then $\|f\|_\infty \le C_4 \delta^{1-1/(2\mym+1-2/r)}$.
Similarly, for $f \in \mathcal V^\mym$ and $0< \delta \le 1$,
if $\|f\|_{L_2} \le \delta$ and $\TV(f^{(\mym-1)})\le 1$, then $\|f\|_\infty \le C_4 \delta^{1-1/(2\mym-1)}$.

\vspace{.2in}
\centerline{\bf\large References}
\begin{description}

\item Bellec, P.C., Lecue, G., Tsybakov, A.B. (2016) Slope meets Lasso: Improved oracle bounds and optimality, arXiv:1605.08651.

\item Birman, M.ˇS. and Solomjak, M.Z. (1967) Piecewise-polynomial approximations
of functions of the classes $W_p^\alpha$, {\em Mathematics of the USSR--Sbornik}, 2, 295--317.

\item Dudley, R.M. (1967) The sizes of compact subsets of Hilbert space and continuity of
Gaussian processes, {\em Journal of Functional Analysis}, 1, 290--330.

\item Guedon, O., Mendelson, S., Pajor, A., and Tomczak-Jaegermann, N. (2007) Subspaces and orthogonal decompositions generated by
bounded orthogonal systems, {\em Positivity}, 11, 269--283.

\item Mammen, E. (1991) Nonparametric regression under qualitative smoothness assumptions, {\em Annals of Statistics}, 19, 741--759.

\item Talagrand, M. (1996) New concentration inequalities in product spaces. {\em Inventiones Mathematicae} 126, 505--563.

\item Tropp, J.A. (2011) Freedman's inequality for matrix martingales,
{\em Electronic Communications in Probability}, 16, 262--270.

\item van de Geer, S. (2014) On the uniform convergence of
empirical norms and inner products, with application to causal inference, {\em Electronic Journal of Statistics},
8, 543--574.
\end{description}

\end{document}